\documentclass[runningheads]{llncs}
\usepackage{amssymb,amsmath}
\usepackage{graphicx}
\usepackage{enumerate}
\usepackage{comment}
\usepackage{epstopdf}
\usepackage[bookmarksnumbered=true,hidelinks]{hyperref}
\usepackage[numbers]{natbib}
\renewcommand\bibname{References}
\usepackage{etoolbox}
\usepackage{lipsum}%
\usepackage{adjustbox}

\usepackage{algorithm}
\usepackage{algorithmicx}
\usepackage{algpseudocode}
\algdef{SE}[DOWHILE]{Do}{doWhile}{\algorithmicdo}[1]{\algorithmicwhile\ #1}%

\patchcmd{\thebibliography}{\chapter*}{\section*}{}{}
\usepackage{todonotes,bm}
\usepackage{url}
\usepackage{cleveref}

\usepackage{subfigure}
\usepackage[subtle, indent=normal, mathspacing=normal]{savetrees}
\usepackage{multirow}
\usepackage{booktabs}
\usepackage{float}
\usepackage{pdflscape}
\usepackage{longtable}
\usepackage{makecell}

\def\Re{\mathbb{R}}

\def\I{\mathcal{I}}

\def\hat{\widehat}

\def \V{\mathrm{V}}

\def\O{{\mathcal O}}

\def\Re{{\mathbb R}}

\def\O{{\mathcal O}}

\newcommand{\exclude}[1]{}

\DeclareMathOperator{\conv}{conv}

\DeclareMathOperator{\Diag}{Diag}

\DeclareMathOperator{\sign}{sign}

\allowdisplaybreaks
 
\usepackage{thm-restate}

\begin{document}

\mainmatter  %

\title{Relax-and-Sparsify Approximation Algorithms for Sparse Ridge Regression}

\titlerunning{Relax-and-Sparsify Approximation Algorithms for Sparse Ridge Regression}

\author{Shixuan Yi  and Yongchun Li }

\authorrunning{S. Yi and Y. Li}

\institute{ School of Data Science, The Chinese University of Hong Kong, Shenzhen\\
\email{shixuanyi@link.cuhk.edu.cn} \\
\email{yongchunli@cuhk.edu.cn}}

\maketitle

\begin{abstract}
This paper studies sparse ridge regression in the low-rank, high-dimensional setting that often arises when data collection is expensive. By exploiting the low-rank structure, we develop approximation algorithms within a relax-and-sparsify framework: solve a convex relaxation and then sparsify its solution to satisfy the original sparsity constraint. Specifically, we establish a  tight rank-dependent sparsity bound for the convex relaxation
by decomposing its support into two parts. This decomposition in turn motivates two approximation strategies, parameter reduction and binary
redistribution, each exploiting a different part. 
Experiments on synthetic and real
datasets show that the binary redistribution algorithms
outperform the benchmarks in most
instances.
\end{abstract}

\section{Introduction} \label{sec:intro}
Sparse ridge regression is a widely used statistical technique for feature selection that balances prediction accuracy and model interpretability \cite{liu2023okridge,tan2025screening}. It has been applied in various fields, including biomedical analysis \cite{bica2025learning,li2025cell}, clinical prediction \cite{chen2024fast,grzeszczyk2025regscore}, and signal processing \cite{saha2023harfe}. This paper studies the sparse ridge regression problem of the form:
\begin{equation}\label{eq:srr}
\V_{\mathrm{opt}}(k):=   \min_{ \beta \in \Re^p} \; \left\{ \frac{1}{n}\| y-  X \beta\|_2^2 + \lambda\| \beta\|_2^2: \, \| \beta\|_0 \le k\right\},
\end{equation}
where $ X \in \Re^{n\times p}$ denotes the matrix of features, $ y\in \Re^n$ denotes the response vector, $k\le p$ is a positive integer,  $\lambda>0$ is a Tikhonov regularization coefficient to address collinearity, and $\|\beta\|_0$ denotes the number of nonzero entries in $\beta$. Problem \eqref{eq:srr} aims to select the best subset of at most $k$ features from a total of $p$ features to minimize the least squares error with a squared $\ell_2$ penalty.

Let $r$ denote the rank of $X$. In this paper, we are interested in instances with $r\ll p$. Such instances arise naturally in high-dimensional regression with limited observations, particularly when  data collection is costly. For example, gene-expression datasets often contain thousands of candidate predictors but only dozens of observations \cite{buhlmann2014high}. Such an imbalance is also common in clinical studies, since the high costs of recruiting
patients, processing biological samples, and collecting outcomes limit the sample size \cite{bentley2019conducting,myers2019recruitment}. Our numerical experiments
in Subsection \ref{subsec:real} include several healthcare datasets with this
$n\ll p$ structure. Since $r\le n$,  $n\ll p$ implies directly $r\ll p$.

Problem \eqref{eq:srr} is NP-hard in general \cite{chen2019approximation}. We demonstrate that it remains NP-hard for every fixed rank $r\ge 2$, so low rank alone does not make problem \eqref{eq:srr} tractable. Nevertheless, we exploit this low-rank structure to develop effective approximation algorithms with provable performance guarantees. Our algorithms adopt a  relax-and-sparsify strategy: they first solve the convex
relaxation \eqref{eq:rel} of \eqref{eq:srr} and then reduce the support of its solution to satisfy
the sparsity constraint in \eqref{eq:srr}. The key is that the rank $r$ controls how much  sparsification is required.

\subsection{Literature review}

In this subsection, we review relevant literature on three aspects: convex relaxations, sparsity bounds, and approximation algorithms.

\noindent\emph{Convex Relaxations.} The lower bounds of \eqref{eq:srr} have been derived in various ways. A widely used
relaxation is based on the squared $k$-support norm
\cite{argyriou2012sparse}. Subsequent work has developed alternative
representations of this norm
\cite{mcdonald2014spectral,sander2023fast}. Equivalent relaxations have also been derived using the Woodbury matrix identity \cite{bertsimas2020sparse}, perspective reformulation \cite{xie2020scalable}, and majorization theory \cite{li2025partial}.  In addition, a variety of semidefinite programming (SDP) relaxations have also been proposed by \cite{atamturk2025rank,cifuentes2026solving,hou2026sparsity,pilanci2015sparse}. Although SDP relaxations can provide tighter lower bounds, the
computational cost can
limit their scalability. We therefore focus on the squared
$k$-support norm relaxation, which can be solved efficiently using
first-order methods (see, e.g., \cite{argyriou2012sparse,ijcai2016p252,mcdonald2014spectral}). 

\vspace{\baselineskip}
\noindent\emph{Sparsity Bounds.} The literature on sparsity bounds for  the sparse ridge regression problem \eqref{eq:srr} remains limited. To the best of our knowledge, two particularly relevant works are \cite{askari2022approximation,li2025partial}. Specifically, 
\cite{askari2022approximation} derived the bound $r+k+2$ using the Shapley--Folkman Theorem, which can be attained by solving an auxiliary linear program. \cite{li2025partial} instead used majorization theory to derive the bound $r+k$. Their analysis applies to a broad class of problems with permutation-invariant feasible sets. Moreover, they proposed a general algorithm to attain the bound, but no polynomial-time complexity guarantee was provided.

Although the prior works are built on different convex relaxation formulations of \eqref{eq:srr}, each formulation can be equivalently converted into the squared $k$-support norm relaxation (see \cite{li2025partial,mcdonald2014spectral,xie2020scalable}). Consequently, the differences in the resulting bounds stem from the proof techniques employed. In this paper, we exploit an invariance property of the squared $k$-support norm and the subgradient optimality condition to further sharpen the sparsity bound to $r+k-1$.

\vspace{\baselineskip}
\noindent\emph{Approximation Algorithms.} We review approximation algorithms for \eqref{eq:srr} that admit
provable performance guarantees.
Iterative hard thresholding (IHT) is a commonly used technique for developing approximation algorithms of cardinality-constrained optimization \cite{blumensath2009iterative,jain2014iterative}. For \eqref{eq:srr}, \cite{yuan2020dual} proposed a Dual IHT (DIHT) algorithm  and established convergence bounds under regularity assumptions. 
Combinatorial algorithms that search directly over supports, including greedy and local search, have been developed for various sparse optimization problems (see, e.g., \cite{JMLR:v14:bahmani13a,hazimeh2020fast,li2025exact}). Specifically, \cite{xie2020scalable} proposed an efficient greedy algorithm for  \eqref{eq:srr}. More recently, \cite{liu2023okridge} developed a beam search variant that retains multiple candidate supports at each step. 
Another thread of work focuses on randomized rounding algorithms based on convex relaxations of \eqref{eq:srr}, including \cite{pilanci2015sparse} and \cite{xie2020scalable}.
However, these algorithms do not exploit the rank-dependent sparsity bound. We take a different angle and make this bound central to the development of our
approximation algorithms. 

\subsection{Contributions and Organization}
The main contributions of this paper are summarized as follows.
\begin{enumerate}[(i)]
    \item \textbf{Tight sparsity bound.} We decompose the support of an optimal solution to the squared $k$-support norm relaxation into two parts. 
    We show that these two parts contain at most $k-1$ and $r$ nonzero entries, respectively. This yields a new sparsity bound $r+k-1$, which improves the existing bounds in the literature. We further prove that this bound is tight and develop a polynomial-time algorithm to find an optimal relaxation solution  satisfying the bound;
    \item \textbf{Approximation algorithms.} Motivated by the support decomposition and sparsity bound, we propose two complementary algorithmic strategies, \emph{Parameter Reduction} and \emph{Binary Redistribution},  for constructing feasible solutions to \eqref{eq:srr}. We also establish
their approximation guarantees using the invariance and monotonicity
properties of the squared $k$-support norm.
    
Specifically, for $k\ge r$,  the  parameter reduction  algorithm achieves an approximation ratio of $k/(k-r+1)$ that depends only on $k$ and $r$.
Binary redistribution constructs a $k$-sparse solution to \eqref{eq:srr} by
redistributing at most $r$ coefficients
while keeping all other coefficients fixed.
We develop an algorithm based on binary quadratic optimization
to find the best redistribution and a scalable variant using deterministic pipage rounding. Both admit the same approximation ratio; and
    \item \textbf{Numerical study.}  
Experiments on synthetic and real datasets show that the binary
redistribution algorithms outperform the baselines in most tested instances. Notably, the pipage rounding variant substantially improves computational efficiency than the optimization-based algorithm while achieving comparable solution quality. 
\end{enumerate}

\noindent\textit{Organization:} The remainder of the paper is organized as follows. \Cref{sec:bound} studies the computational complexity of problem \eqref{eq:srr} and derives the sparsity bound. \Cref{sec:algo} develops  approximation algorithms and establishes their performance guarantees. \Cref{sec:num} presents the numerical study. Finally, \Cref{sec:conclude} concludes the paper.\\
\\
\noindent\textit{Notation:} We let $\Re^n$ denote the set of all the $n$-dimensional vectors. Given two positive integers $q\le p$, we let $[p]=\{1,2,\dots,p\}$, and let $[q,p]=\{q,q+1,\dots,p\}$. We let $I_n$ be the $n\times n$ identity matrix, and let $1_n$ and $0_n$ be the $n$-dimensional all-ones and zero vectors, respectively. For a vector $ x \in \Re^p$, we let $\Diag(x)$ denote a diagonal matrix with diagonal entries being $x$, and let $\|x\|_1$ and $\|x\|_2$ denote its $\ell_1$-norm and Euclidean norm, respectively. We let $|x|$ denote a vector containing absolute entries of $x$, and let $|x|_i^{\downarrow}$ denote the $i$-th largest entry of $|x|$. For a subset $S\subseteq [p]$, let $x_S$ denote the subvector of $x$ indexed by $S$. For a matrix $X\in\Re^{n\times p}$, we let $X_i\in\Re^n$ denote its $i$-th column vector for each $i\in[p]$, and for a subset $S\subseteq [p]$, let $X_S$ denote the submatrix consisting of the columns of $X$ indexed by $S$. For a symmetric matrix $X\in\Re^{n\times n}$, we let $\lambda_{\max}(X)$ denote its largest eigenvalue. Additional notation will be introduced later as needed. 

\section{Convex Relaxation and Sparsity Bound}\label{sec:bound}
In this section, we first show that \eqref{eq:srr} remains NP-hard for
every fixed rank $r\ge2$. We then present the  widely used convex relaxation proposed in \cite{argyriou2012sparse} and derive a tight sparsity bound for its optimal solution.

\subsection{Complexity and Rank Reduction}
We demonstrate that
the sparse ridge regression problem \eqref{eq:srr} remains NP-hard even when the rank of $X$ is fixed and small. This motivates the development of approximation
algorithms in the low-rank setting.
Moreover, we present a dimensionality reduction that reduces the number of
rows of $X$ from $n$ to $r$ without changing the optimal solutions.

\begin{restatable}{theorem}{themNP}\label{them:NP}
For any fixed rank $r\ge 2$, problem \eqref{eq:srr} is NP-hard.     
\end{restatable}
\begin{proof}
    See \Cref{proof:themNP}. \qed
\end{proof}

Although low rank does not resolve the computational
difficulty of \eqref{eq:srr}, it allows a useful
dimensional reduction.
Let
$X=U\Sigma V^{\top}$ be a thin singular value decomposition of $X$, where $U
\in \Re^{n\times r}, \Sigma\in \Re^{r\times r}$, and $V\in \Re^{p\times r}$. 
\begin{lemma}\label{lem:dim}
For any $\beta \in \Re^p$,
we have
\[
 \|y-X\beta\|_2^2
 =\|U^{\top}y-\Sigma V^{\top}\beta\|_2^2
  +\|(I-U U^{\top})y\|_2^2.
\]
\end{lemma}
\begin{proof}
We can decompose
\[
y-X\beta
=
U(U^\top y-\Sigma V^\top\beta)
+(I-UU^\top)y.
\]
The two terms are orthogonal because
$U^\top(I-UU^\top)=0$.
Together with the fact that $U^\top U=I_r$,
taking squared Euclidean norms on both sides above yields the identity.
\qed
\end{proof}
\begin{remark}
  Since the last term in \Cref{lem:dim} is independent of $\beta$, replacing
$(X,y)$ by $(\Sigma V^{\top},U^{\top}y)$ and dropping the constant term preserves the form of \eqref{eq:srr} and its optimal solutions.
Hence, we continue to work with \eqref{eq:srr} throughout this paper, and assume, without loss of generality, that $X\in\mathbb{R}^{r\times p}$ has
full row rank and $y\in \Re^r$. We keep the factor $1/n$ in problem \eqref{eq:srr}, where $n$ denotes the original sample size.

This reduction simplifies our subsequent analysis and may reduce
computational costs when $r< n$.
Moreover, since the dropped constant is
nonnegative, the approximation ratios established in this paper may become tighter when transferred to the original problem.
\end{remark}

\subsection{Convex Relaxation}
Our approximation algorithms for \eqref{eq:srr} leverage the convex relaxation based on the squared $k$-support norm \cite{argyriou2012sparse}.  
Let us first introduce the squared $k$-support norm.
\begin{definition}[Squared $k$-support norm]\label{def:knorm}
For a vector $ \beta \in \Re^p$ and an integer $k\in [p]$, suppose that $\iota(\beta)$, $0\le \iota(\beta) \le k-1$, is  an integer such that  $$|\beta|^{\downarrow}_{\iota(\beta)} >  \Delta(\beta):=\frac{1}{k-\iota(\beta)} \sum_{i\in[\iota(\beta)+1,p]} |\beta|^{\downarrow}_{i} \ge |\beta|^{\downarrow}_{\iota(\beta)+1},$$ with the convention $|\beta|^{\downarrow}_0=\infty$. Then the square of the $k$-support norm is
\begin{equation}\label{eq:phi}
  \phi( \beta; k):=\sum_{i\in [\iota(\beta)]} \left(|\beta|^{\downarrow}_i\right)^2 + (k-\iota(\beta))  \left(\Delta(\beta)\right)^2.
  \end{equation}
 $\sqrt{\phi(\beta; k)}$ is the $k$-support norm.
\end{definition}
Note that the integer $\iota(\beta)$  is unique and always exists, a technical result from \citet[Lemma 14]{nikolov2015randomized}.

As established by \cite{argyriou2012sparse}, the function $\phi$ exactly characterizes the convex hull of the intersection between an $\ell_2$-ball and an $\ell_0$-norm constraint set:
\[\conv\left(\{ \beta \in \Re^p: \| \beta\|_2^2\le 1, \| \beta\|_0\le k \}\right) = \{ \beta \in \Re^p:   \phi( \beta; k)\le 1 \}. \] 
This naturally leads to the following convex relaxation for \eqref{eq:srr}:
\begin{equation}\label{eq:rel}
\V_{\mathrm{rel}}(k):=   \min_{ \beta \in \Re^p} \; \frac{1}{n}\| y-  X \beta\|_2^2 + \lambda\phi( \beta; k) \le \V_{\mathrm{opt}}(k). 
\end{equation}
The relaxation \eqref{eq:rel} can be solved efficiently using
first-order methods, including proximal gradient and
Frank--Wolfe algorithms (see, e.g., \cite{argyriou2012sparse,ijcai2016p252,mcdonald2014spectral}).

\subsection{Support Decomposition and Sparsity Bound}
In this subsection, we establish a tight sparsity bound for
an optimal solution of \eqref{eq:rel}. Our approach is to split the support of the optimal solution into two parts and bound the size of each part separately.

By \Cref{def:knorm},  the entries of $\beta$ whose absolute values exceed $\Delta(\beta)$ contribute individually to $\phi(\beta; k)$, whereas the remaining entries contribute only through the sum of their absolute values. 
This structure allows us to redistribute the latter entries while keeping $\phi(\beta; k)$ unchanged, under the conditions in \Cref{lem:inv}. 

To facilitate our analysis of the sparsity bound, we restrict redistribution to nonzero entries of $\beta$ whose absolute values are strictly below $\Delta(\beta)$ and keep all other entries fixed. 
This leads to the following decomposition of the support  of $\beta$. Specifically, \Cref{def:decom} introduces notation for the nonzero indices below $\Delta(\beta)$ and at least $\Delta(\beta)$, respectively.
\begin{definition}[Support decomposition] \label{def:decom}
For a vector $\beta\in\Re^p$ and an integer $k\in[p]$, let $\Delta(\beta)$ be as in  \Cref{def:knorm}. Define the \emph{redistribution index set}:
\[
    \I(\beta):=\{i\in[p]\ :\ 0<|\beta_i|<\Delta(\beta)\}
\] 
and the number of nonzero entries of $\beta$ outside $\I(\beta)$:
\begin{align*}
\upsilon(\beta)
 :=
 \left|\left\{i\in[p] \ : \
 |\beta_i|> 0,\ 
 |\beta_i|\ge \Delta(\beta)
 \right\}\right|.
\end{align*}
\end{definition}

\begin{remark}
The support decomposition yields the following identity
\begin{equation}\label{eq:partition}
\|\beta\|_0=\upsilon(\beta)+|\I(\beta)|.
\end{equation}
This identity reduces sparsity bound of $\beta$ to bounding $\upsilon(\beta)$ and $\I(\beta)$, respectively. 
\end{remark}

We first establish properties of $\upsilon(\beta)$, including
an upper bound on $\upsilon(\beta)$ whenever $\beta$ has more
than $k$ nonzero entries.
\begin{lemma}\label{lem:upsilon} For any $ \beta \in \Re^p$ and $k\in [p]$,  the following hold:
\begin{enumerate}[(i)]
    \item For any $i\in[\iota(\beta)+1, \upsilon(\beta)]$, 
$|\beta|^\downarrow_i=\Delta(\beta)$, where $\iota(\beta)$ and $\Delta(\beta)$ follow from \Cref{def:knorm}; and
\item If $\|\beta\|_0\ge k+1$, then  $\upsilon(\beta)\le k-1$. 
\end{enumerate}   
\end{lemma}
\begin{proof}
By \Cref{def:knorm}, we have $|\beta|^\downarrow_i>\Delta(\beta)$ if $i\le \iota(\beta)$ and $|\beta|^\downarrow_i\le\Delta(\beta)$ if $i\ge\iota(\beta)+1$. Furthermore, \Cref{def:decom} gives $|\beta|^\downarrow_i\ge \Delta(\beta)$ for all $i\in [\upsilon(\beta)]$. Combining these inequalities yields \[ |\beta|^\downarrow_i=\Delta(\beta), ~ \forall i\in[\iota(\beta)+1, \upsilon(\beta)].\]

When $\|\beta\|_0 \geq k+1$, we can prove $\upsilon(\beta) \leq k-1$ by contradiction. Suppose that $\upsilon(\beta) \geq k$. Since $\upsilon(\beta)\ge k\ge\iota(\beta)+1$,
Part (i) implies that
        \[ \sum_{i\in [\iota(\beta)+1, k]} |\beta|^{\downarrow}_i= (k-\iota(\beta)) \cdot \Delta(\beta) = (k-\iota(\beta)) \cdot \frac{\sum_{i=\iota(\beta)+1}^{p} |\beta|^{\downarrow}_i}{k-\iota(\beta)} = \sum_{i\in [\iota(\beta)+1, p]}|\beta|^{\downarrow}_i,
        \]
 which contradicts $\|\beta\|_{0} \ge k+1$. Thus, $\upsilon(\beta) \leq k-1$ must hold when $\|\beta\|_0\ge k+1$. \qed
\end{proof}

We next turn to the redistribution index set $\I(\beta)$. We investigate an invariance property of the function $\phi(\beta; k)$ based on $\I(\beta)$.
The following lemma gives explicit conditions under which modifying
the entries indexed by $\I(\beta)$ preserves $\phi(\beta;k)$. This plays a key role in deriving the sparsity bound  and developing approximation algorithms.

\begin{lemma}\label{lem:inv}
For any $\beta \in \Re^p$ and  $k\in [p]$, let 
$\I:=\I(\beta).$
 If $\bar\beta\in\Re^p$ satisfies
\[\bar \beta_i=\beta_i, \ \ \forall i\notin \I, \quad |\bar \beta_i|\le \Delta(\beta), \ \ \forall i\in \I, \quad \|\bar\beta_{\I}\|_1=\|\beta_{\I}\|_1, \]
then we have $\phi(\bar \beta; k)=\phi(\beta; k)$. 
\end{lemma}
\begin{proof}
For notational simplicity, let $\iota:=\iota(\beta)$ and $\Delta:=\Delta(\beta)$. By \Cref{def:knorm}, the first $\iota$ largest absolute entries of $\beta$ strictly exceed $\Delta$ and are therefore outside $\I$. They remains unchanged in $\bar\beta$, which all other entries of $\bar\beta$ have absolute value at most $\Delta$. It follows that
\[ |\bar\beta|^\downarrow_i= |\beta|^\downarrow_i, \ \forall i\in [\iota].\]

Since $\bar\beta_i=\beta_i$ for all $i\notin\I$ and
$\|\bar\beta_{\I}\|_1=\|\beta_{\I}\|_1$, we have
$\|\bar\beta\|_1=\|\beta\|_1$.
    Therefore, we can obtain that
    \[\sum_{i\in[\iota+1,p]}|\bar\beta|^\downarrow_i=\|\bar\beta\|_1-\sum_{i\in[\iota]}|\bar\beta|^\downarrow_i=\|\beta\|_1-\sum_{i\in[\iota]}|\beta|^\downarrow_i=\sum_{i\in[\iota+1,p]}|\beta|^\downarrow_i.
    \]
    
According to the results above, we can show that
    \[
        |\bar\beta|^\downarrow_{\iota}=|\beta|^\downarrow_{\iota}>\Delta=\frac{1}{k-\iota}\sum_{i\in[\iota+1,p]}|\beta|^\downarrow_i=\frac{1}{k-\iota}\sum_{i\in[\iota+1,p]}|\bar\beta|^\downarrow_i\geq|\bar\beta|^\downarrow_{\iota+1},
    \]
which, by \Cref{def:knorm}, implies that the same values of $\iota$ and $\Delta$ remain valid for $\bar\beta$. Consequently,
    \[\phi(\bar\beta;k)=\sum_{i\in[\iota]} \left(|\bar\beta|^\downarrow_i\right)^2+(k-\iota)\Delta^2=\sum_{i\in[\iota]}\left(|\beta|^\downarrow_i\right)^2+(k-\iota)\Delta^2=\phi(\beta; k).
    \]
    We thus complete the proof. \qed
\end{proof}

We now demonstrate that the second term $|\I(\beta)|$ in  \eqref{eq:partition} is at most $r$ for some
optimal solution of \eqref{eq:rel}. If an optimal solution $\beta^*$ satisfies $|\I(\beta^*)|>r$, then, by leveraging 
 \Cref{lem:inv} and the subgradient
optimality condition, we can construct a direction that strictly reduces the size of the redistribution index set without
hurting optimality. Repeating this procedure yields
the desired bound.
\begin{theorem}\label{them:index_bound}
The convex relaxation \eqref{eq:rel} admits an optimal solution $\beta^*$  satisfying
\[|\I(\beta^*)| \le   r. \]
 \end{theorem}

\begin{proof}
Let $\beta^*$ be any optimal solution of \eqref{eq:rel}, and
let $\I:=\I(\beta^*)$ and $\Delta:=\Delta(\beta^*)$.
If $|\I|\le r$, the result follows immediately. Otherwise,
we construct another optimal solution with a strictly
smaller redistribution index set.

Since the submatrix $X_{\I}\in \Re^{r\times |\I|}$ has $r$ rows and more than $r$ columns, there exists a nonzero vector $d_{\I}\in \Re^{|\I|}$ such that
\begin{align}\label{eq:data-zero}
 X_{\I} d_{\I}=0.
\end{align}
We extend $d_{\I}$ to a vector $d\in \Re^p$ by setting $d_i=0$ for all $i\notin \I$. 

Next, we define a vector $\sigma \in \Re^p$ as
\[\sigma_i:=\sign(\beta^*_i), \; \forall i\in \I, \quad \sigma_i :=0, \;\forall i\notin \I. \]
Since $0<|\beta^*_i|<\Delta$ for all $i\in \I$, it follows from \cite[Eq.~(60)]{richard2014tight} that every subgradient $g\in \Re^p$ of the function
$\phi$ at $\beta^*$ satisfies
\begin{align}\label{eq:subgrad}
  g_{\I}= 2\Delta \sigma_{\I} \in \Re^{|\I|}.
\end{align}
By the subgradient optimality condition of \eqref{eq:rel} at $\beta^*$ (see \cite[Eq.~(7.1)]{rockafellar1974conjugate}), there exists a subgradient $ g$ of $\phi(\beta^*, k)$  such that
\[ \frac{2}{n}X^{\top} (X\beta^*- y) + \lambda g=0.\]
Restricting the above equation to $\I$ and using \eqref{eq:subgrad}, we obtain
\[\frac{2}{n} X_{\I}^{\top} (X\beta^*- y) + 2\lambda \Delta \sigma_{\I}=0 \Longrightarrow  \frac{2}{n} d_{\I}^{\top}X_{\I}^{\top} (X\beta^*- y) + 2\lambda \Delta d_{\I}^{\top}\sigma_{\I}=0,\]
which, together with  \eqref{eq:data-zero}, yields
\begin{align}\label{eq:sign-zero}
    \sigma_{\I}^{\top}d_{\I}=0.
\end{align}

We now perturb $\beta^*$ along the direction $d$. For $0<\delta\le \bar \delta$, define 
 $\hat{ {\beta}} := 
       \beta^*
        + \delta  d$, where
\begin{equation}\label{eq:delta}
\bar \delta :=\max_{t} \left\{t \; : \;   0 \le \sigma_i(\beta_i^*+ t d_i) \le \Delta, \;\forall i\in \I\right\}.
  \end{equation}
We show that the function $\phi$ remains unchanged along this perturbation.
Since $d_i=0$ for all $i\notin \I$, we have $\hat \beta_i=\beta^*_i$ for all $i\notin \I$. By \eqref{eq:delta}, the entries  of $|\hat \beta|$ indexed by $\I$ are at most $\Delta$. It follows that
\[
    \sum_{i\in \I} |\hat \beta_i| =  \sum_{i\in \I} \sigma_i(\beta_i^*+\delta d_i) = \sum_{i\in \I} | \beta^*_i| + \delta \sigma_{\I}^{\top} d_{\I}= \sum_{i\in \I} | \beta^*_i|,
\]
where the last equality is because of \eqref{eq:sign-zero}. Hence, all three conditions of \Cref{lem:inv} hold, and we obtain $\phi(\beta^*,k)= \phi(\hat \beta,k)$.

In addition, \eqref{eq:data-zero} implies that
\[X\hat \beta = X \beta^* + \delta X d =X \beta^* + \delta X_{\I} d_{\I}= X \beta^*. \]
Combining this with $\phi(\beta^*,k)=\phi(\hat \beta; k)$, we can conclude that the objective values at $\beta^*$ and $\hat \beta$ coincide for any $0<\delta\le \bar \delta$. Thus, $\hat \beta$ is also optimal for \eqref{eq:rel}. 

When $\delta=\bar \delta$, at least one of the inequalities in \eqref{eq:delta} is active. Hence, at least one index
$\ell\in\I$ satisfies $|\hat\beta_\ell|\in\{0,\Delta\}$. The proof of \Cref{lem:inv} also gives
$\Delta(\hat\beta)=\Delta$. Therefore,  the redistribution
index set strictly shrinks, that is, $\I(\hat\beta) \subsetneq\I$. Replacing $\beta^*$ with $\hat\beta$ and repeating this
procedure whenever $|\I(\beta^*)|>r$ yields the desired bound.
We thus complete the proof.
 \qed
\end{proof}

According to the support decomposition \eqref{eq:partition}, we can combine \Cref{them:index_bound} with Part (ii) of \Cref{lem:upsilon} to obtain the following 
sparsity bound.
\begin{theorem}\label{cor:bound}
There is an optimal solution $\beta^*$ to \eqref{eq:rel} satisfying
\[\|\beta^*\|_0 \le  r+k-1. \]
 \end{theorem}
 \begin{proof}
 According to \Cref{them:index_bound}, \eqref{eq:rel} admits an optimal solution $\beta^*$ satisfying $|\I(\beta^*)|\le r$.
 If $\|\beta^*\|_0\le k$, the result is trivial.
If $\|\beta^*\|_0\ge k+1$, then we have $\upsilon(\beta^*)\le k-1$  based on Part (ii) of \Cref{lem:upsilon}. Then it follows from \eqref{eq:partition} that $$\|\beta^*\|_0= \upsilon(\beta^*) + |\I(\beta^*)|\le r+k-1. $$
We complete the proof.
\qed
 \end{proof}

 \begin{remark}
    The sparsity bound in \Cref{cor:bound} is a uniform worst-case guarantee. For a given instance, the support decomposition also yields the instance-dependent sparsity bound $\upsilon(\beta^*)+r$ for an optimal solution $\beta^*$ to \eqref{eq:rel}, which can be tighter than $r+k-1$. 
    Indeed, our numerical experiments show that $\upsilon(\beta^*)$ is often substantially smaller than $k-1$. 
    We also develop an algorithm to achieve the tighter bound.
\end{remark}

\begin{corollary}\label{cor:exact}
    When $r=1$, the relaxation \eqref{eq:rel} and the original problem \eqref{eq:srr} have the same optimal value.
\end{corollary}

\begin{restatable}{proposition}{proptight}\label{prop:tight}
There exists an instance of \eqref{eq:rel} that admits a unique optimal solution $\beta^*$ satisfying 
\[ \|\beta^*\|_0=r+k-1. \] 
\end{restatable}
\begin{proof}
    See \Cref{proof:proptight}. \qed
\end{proof}

\subsection{Optimality-Preserving Sparsification Algorithm}
Although \eqref{eq:rel} can be solved efficiently using first-order methods,  these methods may return a dense optimal solution $\beta^*$ that violates the sparsity bound in \Cref{cor:bound}. 
To close this gap, we develop an optimality-preserving sparsification algorithm to attain this bound.  At each iteration, the algorithm constructs a nonzero direction  and computes the largest admissible
step size as in the proof of \Cref{them:index_bound}. The updated solution strictly reduces the  size of the redistribution index set. Repeating this procedure yields an optimal solution that satisfies the sparsity bound in \Cref{cor:bound}. The detailed implementation is presented in
\Cref{alg:sparsity_reduction}.

\begin{algorithm}[H]
    \caption{Optimality-Preserving Sparsification Algorithm}
    \label{alg:sparsity_reduction}
    
    \begin{algorithmic}[1]
        \State \textbf{Input:} An optimal solution $ {\beta}^*$ to \eqref{eq:rel}, matrix $ {X}$, and integer $k$
        
\item Compute $\Delta=\Delta(\beta^*)$ and let $\I=\I(\beta^*)$ 
\While{$|\I| > r$ and $\|\beta^*\|_0\ge k+1$}
\State Find a nonzero vector $d\in \Re^p$ satisfying
\[ X_{\I} d_{\I}=0, \quad d_i=0, \; \forall i\notin \I. \]
\State Compute the largest step size
            $$ \bar \delta :=\max_{t} \left\{t \; : \;   0 \le \sign(\beta_i^*)(\beta_i^*+ t d_i) \le \Delta, \;\forall i\in \I\right\}. $$
           
            \State Update the solution $\beta^*\leftarrow\beta^*+\bar \delta d$ 
            \State Update the set $\I\leftarrow \I\setminus\left\{i\in \I:|\beta^*_i|\in\left\{0,\Delta\right\}\right\}$
        \EndWhile
        
        \State \textbf{return} \(\beta^*\)
    \end{algorithmic}
\end{algorithm}

\begin{proposition}\label{prop:algo_sparse}
\Cref{alg:sparsity_reduction} runs in $\O(p\log p+(p-r)r^3)$ time and  returns an optimal solution $\beta^*$ to \eqref{eq:rel} satisfying
\[\|\beta^*\|_0 \le  \upsilon(\beta^*) +r.\]
\end{proposition}
\begin{proof} 
The sparsity bound for $\beta^*$ follows directly from \eqref{eq:partition} and  \Cref{them:index_bound}; we thus omit  the details.

Next, we derive the time complexity of \Cref{alg:sparsity_reduction}. In Step 2, we can compute $\Delta$ and $\upsilon$ by sorting the $p$ entries of $|\beta^*|$, which runs in $\O(p\log p)$ time (see \cite{cormen1990introduction}). 

The loop continues only while $|\I|\ge r+1$. Thus, we can implement Step 5
by selecting any subset $S\subseteq \I$ with $|S|=r+1$ and finding
a nonzero vector $d_S$ satisfying $X_Sd_S=0$. Since
$X_S\in\Re^{r\times(r+1)}$, such a vector $d_S$ always exists and can be
computed in $\O(r^3)$ operations (see \cite{golub2013matrix}). Setting
$d_i=0$ for all $i\notin S$ gives a valid direction $d$.
Because the entries of $\beta^*$ outside $S$ remain unchanged, Steps 6--7
need only update the entries indexed by $S$, which requires $\O(r)$ operations in total.

By the definition of $\bar \delta$ in Step 5,
at least one updated entry
indexed by $\I$ has absolute value $0$ or $\Delta$. Hence, $|\I|$ decreases by at least one in each
iteration. Since $r+1\le |\I|\le p$, \Cref{alg:sparsity_reduction} terminates after at most
$p-r$ iterations. Together with the $\O(p\log p)$ initialization cost, the total time complexity of \Cref{alg:sparsity_reduction} is
$\O(p\log p+(p-r)r^3)$. \qed
\end{proof}

\section{Relax-and-Sparsify Approximation Algorithms}
\label{sec:algo}
In this section, we develop approximation algorithms for \eqref{eq:srr} using the properties of $\phi$ and the sparsity bound of \eqref{eq:rel}. They follow a common
relax-and-sparsify framework: solve the convex relaxation \eqref{eq:rel} and then sparsify its optimal solution to satisfy the constraint $\|\beta\|_0\le k$ in \eqref{eq:srr}.
The support
decomposition \eqref{eq:partition} motivates two complementary algorithmic strategies to achieving this sparsity.

Specifically, according to the support decomposition, $\|\beta\|_0$ splits into two terms $\upsilon(\beta)$ and $|\I(\beta)|$. To achieve $k$-sparsity, the \emph{Parameter Reduction} and \emph{Binary Redistribution}
algorithms exploit $\upsilon(\beta)$ and $|\I(\beta)|$,
respectively. Parameter reduction decreases the bound on $\upsilon(\beta)$ by
reducing the relaxation parameter, while binary redistribution aims to reduce the number of nonzero entries within the redistribution index set $\I(\beta)$.

\subsection{Properties of The Squared $k$-Support Norm}
In this subsection, we derive properties of the squared $k$-support norm $\phi$, which are key for deriving approximation ratios of approximation algorithms. 

First, we show that $\phi(\beta;k)$ coincides with
$\|\beta\|_2^2$ whenever 
 $\beta$ has at most $k$ nonzero entries.
\begin{lemma}\label{lem:phi_k}
   For any $k\in [p]$, let $ \beta \in \Re^p$ be a vector such that
   $\| \beta\|_0\le k$. Then
   \[\phi( \beta; k) = \| \beta\|_2^2. \]
\end{lemma}
\begin{proof}
By \Cref{def:knorm}, we have
    \[
|\beta|^\downarrow_{\iota(\beta)+1}\leq \Delta(\beta)= \frac{1}{k-\iota(\beta)} \sum_{i\in [\iota(\beta)+1,k]}|\beta|^\downarrow_i\leq \frac{1}{k-\iota(\beta)}  (k-\iota(\beta)) |\beta|^\downarrow_{\iota(\beta)+1}=|\beta|^\downarrow_{\iota(\beta)+1},
    \]  
    where the first equality is because $|\beta|^\downarrow_i=0$ for all $i\in [k+1,p]$. It is clear that both inequalities must hold with equality, which implies
\[
|\beta|^\downarrow_{\iota(\beta)+1}
=\cdots=|\beta|^\downarrow_k
=\Delta(\beta).
\]
    Plugging the above equalities into \eqref{eq:phi} yields
    \[
        \phi(\beta; k)=\sum_{i\in [\iota(\beta)]} \left(|\beta|^\downarrow_i\right)^2+\sum_{i\in [\iota(\beta)+1,k]} \left(|\beta|^\downarrow_i\right)^2=\|\beta\|_2^2,
    \]
    where the last equality again uses $|\beta|^\downarrow_i=0$ for all $i\in [k+1,p]$.
    \qed
\end{proof}

By leveraging \Cref{lem:phi_k} and \cite[Theorem 4.1]{mcdonald2014spectral}, we show that for fixed $\beta$, the function $\phi(\beta; k)$ is monotonically non-increasing with respect to $k$. Furthermore,  we provide a tight upper bound on how much its value can increase when $k$ decreases.
\begin{lemma}\label{lem:phi-mont}
For any vector $ \beta \in \Re^p$ and integers $k_1, k_2\in [p]$ with $k_1\le k_2$, we have 
\begin{equation*}
  \phi( \beta; k_2)  \le \phi( \beta; k_1) \le \frac{k_2}{k_1} \phi( \beta; k_2),
\end{equation*}
and both inequalities are tight. 
\end{lemma}
\begin{proof}
The proof includes two parts. We first establish the inequalities and then prove that equality can be attained in each inequality. If $k_1=k_2$, the result immediately follows. We therefore focus on the nontrivial case where $k_1<k_2$.
\begin{itemize}
    \item[(i)] By \cite[Theorem 4.1]{mcdonald2014spectral}, we can represent $\phi(\beta;k_1)$ and $\phi(\beta;k_2)$ as
\begin{equation}\label{eq:phi_min_k1}
\phi(\beta;k_1)= \inf_{ z\in [0,1]^p} \; \bigg\{
            \sum_{i\in [p]} \frac{\beta_i^2}{z_i} \ : \;  \sum_{i\in [p]} z_i \le k_1
            \bigg\}
        \end{equation}
        and
\begin{equation}\label{eq:phi_min_k2}
\phi(\beta;k_2)= \inf_{ z\in [0,1]^p} \; \bigg\{
            \sum_{i\in [p]} \frac{\beta_i^2}{z_i} \ : \;  \sum_{i\in [p]} z_i \le k_2
            \bigg\}.
        \end{equation}
Thus, it suffices to compare these two minimization problems, which
    differ only in the cardinality constraint. Since $k_1<k_2$, the feasible region of \eqref{eq:phi_min_k1} is contained in that of \eqref{eq:phi_min_k2}. It follows that $\phi(\beta;k_2)\leq\phi(\beta;k_1).$

  To derive the second inequality, let $z^*$ be an optimal solution to problem \eqref{eq:phi_min_k2}, and define
$  \tilde{z}:={k_1}/{k_2}z^*$.
        Since $z^*\in [0,1]^p$ and $k_1< k_2$, we have
    $ \tilde{z}\in[0,1]^p$.
        In addition, given $\sum_{i\in [p]} z^*_i\leq k_2$, we can show that 
        \[
            \sum_{i\in [p]} \tilde z_i =\frac{k_1}{k_2}\sum_{i\in [p]}z_i^*\leq \frac{k_1}{k_2}k_2=k_1.
        \]
        Therefore, $\tilde{z}$ is feasible for \eqref{eq:phi_min_k1} with the  objective value:
        \[
\sum_{i\in [p]}\frac{\beta_i^2}{\tilde{z}_i}=\frac{k_2}{k_1}\sum_{i\in [p]}\frac{\beta_i^2}{z^*_i}=\frac{k_2}{k_1}\phi(\beta;k_2)\ge \phi(\beta; k_1),
        \]
        where the second equality is because $z^*$ is optimal for \eqref{eq:phi_min_k2}, and the last inequality is from the fact that $\phi(\beta; k_1)$ is the optimal value of \eqref{eq:phi_min_k1}. Hence, we have
            \[
        \phi(\beta;k_2)
        \le \phi(\beta;k_1)
        \le \frac{k_2}{k_1}\phi(\beta;k_2).
    \]
\item[(ii)] We next prove that both inequalities are tight.

Consider any nonzero vector $\beta$ satisfying $\|\beta\|_0\leq k_1\le k_2$.  \Cref{lem:phi_k} implies that
        \[
            \phi(\beta;k_1)=\|\beta\|_2^2, \quad \text{and} \quad \phi(\beta;k_2) = \|\beta\|_2^2.
        \]
Hence, the first inequality $   \phi(\beta;k_2)
        \le \phi(\beta;k_1)$ can hold with equality.

For the second inequality, let us define the vector $\beta$ as
        \[\beta_i=a, ~ \forall i \in [k_2], \quad \beta_i=0, ~\forall i \in [k_2+1, p],
        \]
        where $a\neq 0$ is a scalar.
To compute $\phi(\beta; k_1)$, we apply \Cref{def:knorm} with $k=k_1$ in the following. It is easy to verify that
\[ |\beta|^\downarrow_0 > \frac{1}{k_1-0}\sum_{i\in[p]}|\beta|^\downarrow_i=|a|\frac{k_2}{k_1} \ge |a|= |\beta|^\downarrow_1, \]
where $|\beta|^\downarrow_0=\infty$ by convention.  By the conditions and the uniqueness of $\iota(\beta)$ in \Cref{def:knorm}, we must have $\iota(\beta)=0$. Consequently, $\Delta(\beta)=|a|\frac{k_2}{k_1}$.
Plugging these values into \eqref{eq:phi} yields
        \[
 \phi(\beta;k_1)=k_1\left(|a|\frac{k_2}{k_1}\right)^2=\frac{k_2}{k_1}(k_2a^2)=\frac{k_2}{k_1}\|\beta\|_2^2=\frac{k_2}{k_1}\phi(\beta;k_2),
        \]
        where the last equality follows from \Cref{lem:phi_k}, since $\|\beta\|_0=k_2$. Therefore, the second inequality $\phi(\beta;k_1)
        \le \frac{k_2}{k_1}\phi(\beta;k_2)$ holds with equality for this choice of $\beta$.    \qed
\end{itemize}
\end{proof}

\subsection{Parameter Reduction  Approximation Algorithm}\label{subsec:algo2}
According to  \Cref{lem:upsilon}, the upper bound on $\upsilon(\beta)$ is $k-1$, which directly determines the  sparsity bound $k+r-1$ in \Cref{cor:bound}. Inspired by this,
the \emph{Parameter Reduction} algorithm applies the sparsity
bound with a smaller relaxation parameter than $k$. 
Indeed, an optimal solution of \eqref{eq:rel} with parameter $q$ can be sparsified to have at most $q+r-1$ nonzero entries.
Thus, assume  $k\ge r$ and set
\[
    \tilde k:=k-r+1\in [k].
\]
By replacing $k$ with $\tilde k$ in \eqref{eq:rel}, the resulting sparsity bound becomes
$\tilde k+r-1=k$.

Specifically, the parameter reduction algorithm first applies \Cref{alg:sparsity_reduction}   to find an optimal solution of \eqref{eq:rel}. If the resulting solution is $k$-sparse, it is also optimal for \eqref{eq:srr}. Otherwise,  we solve problem \eqref{eq:rel} with the reduced parameter $\tilde k$:
\begin{equation}\label{eq:rel_reduced}
\V_{\mathrm{rel}}(\tilde k):=   \min_{ \beta \in \Re^p} \; \frac{1}{n}\| y-  X \beta\|_2^2 + \lambda\phi( \beta; \tilde k). 
\end{equation}
By \Cref{cor:bound}, problem \eqref{eq:rel_reduced} has an optimal solution $\tilde \beta$ satisfying $\|\tilde \beta\|_0 \le \tilde k+r-1=k.$
Thus, $\tilde \beta$ is  feasible for the original problem \eqref{eq:srr}. 
The detailed implementation is presented in \Cref{alg:reduce}.

\begin{algorithm}[H]
\caption{Parameter Reduction Approximation Algorithm}
\label{alg:reduce}
\begin{algorithmic}[1]
\State \textbf{Input:} An optimal solution $\beta^*$ of \eqref{eq:rel} returned by \Cref{alg:sparsity_reduction} and an integer $k\ge r$
\State Let $\tilde\beta=\beta^*$
\If{\(\|\beta^*\|_0\geq k+1\)}
\State Set $\tilde k=k-r +1$
\State Solve \eqref{eq:rel_reduced}  and let $\hat \beta$ be an optimal solution
\State Apply \Cref{alg:sparsity_reduction} to $\hat\beta$
with parameter $\tilde k$, and denote the output by $\tilde \beta$

\EndIf
\State \textbf{return} \(\tilde \beta\)
\end{algorithmic}
\end{algorithm}

By leveraging \Cref{lem:phi_k,lem:phi-mont}, we derive the theoretical guarantees of \Cref{alg:reduce}. 
\begin{theorem}\label{thm:reduce_approx_ratio}
Assume $k\ge r$ and let $\tilde k=k-r+1$. Let $\tilde \beta \in \Re^p$ be returned by \Cref{alg:reduce}. Then $\tilde \beta$ is feasible for \eqref{eq:srr} and  satisfies
    \begin{align*}
{\frac{1}{n}\|y-X\tilde\beta\|_2^2+\lambda\|\tilde\beta\|_2^2}
 \leq \frac{k}{\tilde k} \V_{\mathrm{opt}}(k).
\end{align*}
Moreover, if the input solution $\beta^*$ in \Cref{alg:reduce} satisfies $\|\beta^*\|_1< \tilde k \|\beta^*\|_{\infty}$, then
    \begin{align*}
{\frac{1}{n}\|y-X\tilde\beta\|_2^2+\lambda\|\tilde\beta\|_2^2}
 \leq \min\left\{\tilde k,
 \frac{k}{\tilde k} \right\}\V_{\mathrm{opt}}(k).  
\end{align*}
\end{theorem}
\begin{proof}
 If $\|\beta^*\|_0\le k$, the algorithm
returns $\beta^*$, which is optimal for \eqref{eq:srr} by
\Cref{lem:phi_k}. Both ratios follow immediately.
We therefore consider the case $\|\beta^*\|_0\ge k+1$.
In this case, we have $\|\tilde\beta\|_0\le r+\tilde k-1=k$
by \Cref{cor:bound}. Thus, $\tilde\beta$ is feasible for
\eqref{eq:srr}.

To prove the approximation ratios, we bound the increase in $\phi$ at $\beta^*$ when the parameter
is reduced from $k$ to $\tilde k$. Suppose there exists a constant
$\rho\ge1$ satisfying
\begin{equation}\label{eq:reduce-bound}
    \phi(\beta^*;\tilde k)\le\rho\,\phi(\beta^*;k).
\end{equation}

Since $\|\tilde\beta\|_0\le k$ and $\tilde k\le k$, we have
\[\|\tilde\beta\|_2^2
    =\phi(\tilde\beta;k)
    \le\phi(\tilde\beta;\tilde k),\]
where the equality is due to \Cref{lem:phi_k} and the inequality is from \Cref{lem:phi-mont}. Leveraging the above results, the objective value of \eqref{eq:srr} at $\tilde \beta$ satisfies
\begin{align*}
    &\frac{1}{n}\|y-X\tilde\beta\|_2^2
    +\lambda\|\tilde\beta\|_2^2
    \le
    \frac{1}{n}\|y-X\tilde\beta\|_2^2+\lambda\phi(\tilde\beta;\tilde k)\le
    \frac{1}{n}\|y-X\beta^*\|_2^2
    +\lambda\phi(\beta^*;\tilde k)\\ & \le
    \rho\left(
        \frac{1}{n}\|y-X\beta^*\|_2^2
        +\lambda\phi(\beta^*;k)
\right)=\rho\,\V_{\mathrm{rel}}(k)
    \le\rho\,\V_{\mathrm{opt}}(k),
\end{align*}
where the second inequality is because $\beta^*$ is a feasible solution for \eqref{eq:rel_reduced}, and the
third one uses \eqref{eq:reduce-bound}, $\rho\ge1$,
and the nonnegativity of the least squares loss. Therefore, \Cref{alg:reduce} admits an approximation ratio $\rho$. 

We next establish two choices of $\rho$ in \eqref{eq:reduce-bound}:
\begin{enumerate}[(i)]
    \item By \Cref{lem:phi-mont}, \eqref{eq:reduce-bound}
    holds with $\rho=k/\tilde k$, giving the first
    approximation ratio.

    \item Suppose that
    $\|\beta^*\|_1<\tilde k\|\beta^*\|_{\infty}$.
    Let $\iota$ and $\Delta$ denote $\iota(\beta^*)$ and
    $\Delta(\beta^*)$, respectively, defined as in \Cref{def:knorm} with  $k$ replaced by $\tilde k$. If $\iota=0$, \Cref{def:knorm} would require $
     \frac{\|\beta^*\|_1}{\tilde k}
    \ge|\beta^*|^\downarrow_1$,
    contradicting the condition. Thus,
    $1\le\iota\le\tilde k-1$.

    By \Cref{def:knorm}, we have
    \begin{align*}
        \phi(\beta^*;\tilde k)
        &=
        \sum_{i\in[\iota]}(|\beta^*|^\downarrow_i)^2
        +(\tilde k-\iota)\Delta^2\le
        \sum_{i\in[\iota]}(|\beta^*|^\downarrow_i)^2
        +(\tilde k-\iota)(|\beta^*|^\downarrow_\iota)^2\\
        &\le
        \frac{\tilde k}{\iota}
        \sum_{i\in[\iota]}(|\beta^*|^\downarrow_i)^2\le
        \frac{\tilde k}{\iota}\|\beta^*\|_2^2
        =
        \frac{\tilde k}{\iota}\phi(\beta^*;p)
        \le
        \frac{\tilde k}{\iota}\phi(\beta^*;k),
    \end{align*}
    where the first inequality uses
    $\Delta<|\beta^*|^\downarrow_\iota$, the second inequality uses
    $(|\beta^*|^\downarrow_\iota)^2
        \le
        \frac{1}{\iota}
        \sum_{i\in[\iota]}(|\beta^*|^\downarrow_i)^2
   $, and the last equality and inequality follow from
    \Cref{lem:phi_k,lem:phi-mont}, respectively.
    Since $\iota\ge1$, \eqref{eq:reduce-bound} also holds
    with $\rho=\tilde k$ under the stated condition. Combining these two choices of $\rho$ proves the strengthened guarantee.

\end{enumerate}
We thus complete the proof. \qed
\end{proof}

As a byproduct of the proof of \Cref{thm:reduce_approx_ratio}, we directly obtain the following bounds  on the gap between the relaxation \eqref{eq:rel}
and the original problem \eqref{eq:srr}.
\begin{corollary}\label{cor:opt_gap}
Assume $k\ge r$ and let $\tilde k=k-r+1$. Then $
\V_{\mathrm{opt}}(k)
 \leq \frac{k}{\tilde k} \V_{\mathrm{rel}}(k)$.
Moreover, if the input solution $\beta^*$ in \Cref{alg:reduce} satisfies $\|\beta^*\|_1< \tilde k \|\beta^*\|_{\infty}$, then
$
\V_{\mathrm{opt}}(k)
 \leq \min\left\{\tilde k,
 \frac{k}{\tilde k} \right\}\V_{\mathrm{rel}}(k).
$
\end{corollary}

To the best of our knowledge, \Cref{thm:reduce_approx_ratio}
provides the first-known  approximation ratios for \eqref{eq:srr} that do not rely on data $X$, $y$ and parameter $\lambda$. However, the assumption $k\ge r$ limits its practical applicability, since $k$ is often
small and may fall below $r$.
This limitation motivates the development of the binary redistribution algorithm in Subsection \ref{subsec:algo1}, which does not require $k\ge r$.

\subsection{Binary Redistribution Approximation Algorithms}
\label{subsec:algo1}
The \emph{Binary Redistribution} approximation algorithms introduce binary redistribution of entries within $\I(\beta^*)$ to 
convert an optimal solution $\beta^*$ of \eqref{eq:rel} into a feasible one of \eqref{eq:srr} while preserving $\phi(\beta^*; k)$.
We first solve a binary quadratic optimization
problem to select the best redistribution. To improve efficiency, we further use
deterministic pipage rounding to yield a binary redistribution in
polynomial time. We establish the same approximation guarantee for both variants.

\subsubsection{Binary Redistribution via Optimization}
Let $\beta^*$ be an optimal solution of \eqref{eq:rel}
returned by \Cref{alg:sparsity_reduction}. If $\|\beta^*\|_0\ge k+1$, let
\[
\I:=\I(\beta^*),
\quad
\upsilon:=\upsilon(\beta^*),
\quad \Delta:=\Delta(\beta^*), \quad \sigma_i:=\sign(\beta_i^*), \ \forall i\in [p].
\]

We first express the $\ell_{1}$ norm of entries within $\I$ explicitly.
\begin{lemma}\label{lem:sum}
We have     $\|\beta^*_{\I}\|_1= (k-\upsilon)\Delta.$ 
\end{lemma}
\begin{proof}
First, \Cref{lem:upsilon} shows that $|\beta|_i^{\downarrow}= \Delta$ for any $i\in [\iota+1, \upsilon]$. Besides, $|\beta|_i^{\downarrow}< \Delta$ for any $i\in [\upsilon+1, p]$. The positive entries in the latter range correspond exactly to the entries of $\beta^*$ indexed by $\I$.
    Thus, by the definition of $\Delta$,
    \begin{align*}
  (k-\iota)  \Delta  =  \sum_{i\in [\iota+1, \upsilon]} |\beta|_i^{\downarrow} + \sum_{i\in [\upsilon+1, p]} |\beta|_i^{\downarrow} =    (\upsilon-\iota)  \Delta + \sum_{i\in \I} |\beta_i|.
    \end{align*}
The result immediately follows from rearranging the above equality.
    \qed
\end{proof}

Binary distribution 
aims to convert $\beta^*$ into a feasible solution $\bar \beta$ for \eqref{eq:srr} from the perspective of the redistribution index set $\I$. Besides, this redistribution preserves $\phi(\beta^*;k)$.
By leveraging \Cref{lem:inv,lem:sum}, we construct such a solution $\bar \beta$ in three steps.
\begin{enumerate}
    \item \textit{Determine the sparsity budget within $\I$.}

    By the support decomposition \eqref{eq:partition},  the entries outside $\I$ contribute exactly $\upsilon$ nonzero
entries. We keep all entries
outside $\I$ fixed. That is
\[\bar \beta_i = \beta_i^*, \ \forall i\notin \I. \]
Thus,  $\bar\beta$ is feasible
for \eqref{eq:srr} if and only if $\|\bar\beta_{\I}\|_0\le k-\upsilon$. By Part (ii) of \Cref{lem:upsilon}, we get $1\le k-\upsilon\le |\I|$. 

\vspace{0.7em}

\item \textit{Construct binary redistributions that preserve
$\phi(\beta^*;k)$.}

Since $\beta^*_i=\bar \beta_i$ for any $i\notin \I$, to guarantee $\phi(\beta^*;k)=\phi(\bar \beta;k)$,
\Cref{lem:inv} shows that it suffices to keep $|\bar \beta_i| \le \Delta$ for each $i\in \I$ and preserve the $\ell_1$ norm over $\I$. By \Cref{lem:sum}, this norm is exactly $(k-\upsilon)\Delta$. Therefore we require
\[|\bar \beta_i|\le \Delta, \ \ \forall i\in \I, \quad \|\bar\beta_{\I}\|_1=(k-\upsilon)\Delta.\]
To meet the sparsity budget for $\bar \beta_{\I}$, these conditions force each  $|\bar \beta_i|$, $i\in \I$, to be either 0 or $\Delta$. Consequently, we have $\|\bar \beta_{\I}\|_0=k-\upsilon$.

We let $z_i\in\{0,1\}$ indicate whether entry
$i\in\I$ is selected. Then $|\bar\beta_i|=\Delta z_i$ for each $i\in \I$ and $\sum_{i\in\I}z_i=k-\upsilon$.  

\vspace{0.7em}

\item \textit{Fix the signs of the entries within $\I$.}

In the previous step, the binary variables determine the redistributed magnitudes. We restrict attention to redistributions that preserve the signs of $\beta_{\I}^*$. Thus, we set $$\bar \beta_i=\sigma_i  \Delta z_i, \ \forall i\in \I.$$
\end{enumerate}
Every binary vector $z=(z_i)_{i\in \I}$ with
$\sum_{i\in\I}z_i=k-\upsilon$ specifies which $k-\upsilon$ entries within $\I$ remain nonzero. Moreover, each such selection uniquely determines the redistributed vector $\bar \beta$. The following corollary formalizes the expression of $\bar \beta$ as a function of $z$ and shows its properties. 
\begin{corollary}
    \label{cor:inv} 
Let  $z=(z_i)_{i\in \I}$ satisfy $z_i\in\{0,1\}$ for all $i\in \I$ and $\sum_{i\in \I} z_i=k-\upsilon$.  Define $\beta(z)\in \Re^p$ by
\begin{equation}\label{eq:beta-round}
    \beta_i(z) :=   \beta_i^*,\  \forall i \notin \I, \quad      \beta_i(z) :=   \sigma_i\Delta z_i, \ \forall i \in \I.
\end{equation}
Then 
\[\|\beta(z)\|_0=k, \quad \text{and} \quad \phi(\beta^*; k)=\phi(\beta(z); k)=\|\beta(z)\|_2^2.\]
\end{corollary}
\begin{proof}
The construction gives $\|\beta(z)\|_0=\upsilon+\sum_{i\in\I}z_i=k$,
and \Cref{lem:inv,lem:sum} yield $\phi(\beta(z);k)=\phi(\beta^*;k)$. Since $\|\beta(z)\|_0=k$, we have $\phi(\beta(z);k)=\|\beta(z)\|_2^2$ based on \Cref{lem:phi_k}.
\qed
\end{proof}

According to \Cref{cor:inv}, $\beta(z)$ is $k$-sparse for any feasible $z$, and the objective value of \eqref{eq:srr} at $\beta(z)$ is $\frac{1}{n}\|y-X\beta(z)\|_2^2+\lambda\phi(\beta^*;k)$. Because the second term and the positive scaling factor $\frac{1}{n}$ are constants, the best binary redistribution is obtained by minimizing the least squares loss. By substituting the definition of $\beta(z)$ in \eqref{eq:beta-round}, this loss can be written as the following quadratic function over $z$:
\[ \|y-X\beta(z)\|_2^2 =\bigg\|y-\sum_{i\notin \I} X_i\beta_i^* -\sum_{i\in \I} \sigma_i\Delta X_i z_i\bigg\|_2^2.\]
Once an optimal binary vector is obtained, we can construct a feasible solution for \eqref{eq:srr}. The detailed implementation is presented in \Cref{alg:round}. 

\begin{algorithm}[H]
\caption{Binary Redistribution via  Optimization}
\label{alg:round}
\begin{algorithmic}[1]
\State \textbf{Input:} An optimal solution $\beta^*$ returned by \Cref{alg:sparsity_reduction}, $X$, $y$, and $k$
\State Let $\hat \beta =\beta^*$
\If{\(\|\beta^*\|_0\geq k+1\)}
\State Compute $\Delta = \Delta(\beta^*)$, $\I=\I(\beta^*)$ and $\upsilon=\upsilon(\beta^*)$
\State Let $\sigma_i=\sign(\beta^*_i), \ \forall i\in\I$

\State Solve the following binary quadratic optimization problem
\begin{equation}\label{eq:sign-miqp}
 \min_{(z_i)_{i\in \I}} \bigg\{\bigg\|y-\sum_{i\notin \I} X_i\beta_i^* -\sum_{i\in \I} \sigma_i\Delta X_i z_i\bigg\|_2^2 ~ : ~\sum_{i\in \I} z_i=k-\upsilon, ~  
    z_i\in \{0,1\}, ~  \forall i\in \I\bigg\} 
\end{equation}
and let $z^{\mathrm{IP}}=(z_i^{\mathrm{IP}})_{i\in \I}$ denote an optimal solution 

\For{$i\in \I$}
\State Update $\hat \beta_i\gets \sigma_i\Delta z_i^{\mathrm{IP}}$
\EndFor
\EndIf
\State \textbf{return} \(\hat \beta\)
\end{algorithmic}
\end{algorithm}

Note that  problem \eqref{eq:sign-miqp} can be  solved directly by off-the-shelf solvers, such as Gurobi and MOSEK. By \Cref{them:index_bound}, $|\I|\le r$. Thus, \eqref{eq:sign-miqp} aims to select $k-\upsilon$ ones from at most $r$ candidates, whereas the original problem \eqref{eq:srr} seeks to select at most $k$ features from $p$ candidates. In high-dimensional and low-rank settings where $r\ll p$, the reduced problem \eqref{eq:sign-miqp} can be much easier to solve than \eqref{eq:srr}.
Unfortunately, as an integer programming problem, \eqref{eq:sign-miqp} is still hard to scale. This motivates us to provide an efficient approximate variant of \Cref{alg:round} in the following subsection.

\subsubsection{Efficient Binary Redistribution via Pipage Rounding} 
Pipage rounding is an efficient, commonly used method that rounds a fractional solution into a binary one \cite{ageev2004pipage,calinescu2011maximizing,harvey2014pipage}. To avoid solving \eqref{eq:sign-miqp} to optimality, we apply pipage rounding to an optimal solution of its continuous relaxation to obtain a feasible binary solution.
According to \cite[Theorem 4.1]{mcdonald2014spectral}, the continuous relaxation of \eqref{eq:sign-miqp} admits an optimal solution that can be expressed in closed-form in terms of $\beta^*$.
\begin{remark}\label{rem:cont-sol}
The continuous relaxation of \eqref{eq:sign-miqp} admits an optimal solution $z^{\mathrm{rel}}$:
\begin{equation}\label{eq:z-rel}
    0< z^{\mathrm{rel}}_i=\frac{|\beta^*_i|}{\Delta}<1, \quad \forall i\in \I.
\end{equation}
\end{remark}

Given the fractional solution $z^{\mathrm{rel}}$, pipage rounding iteratively adjusts two fractional coordinates while preserving the cardinality $k-v$, until all coordinates become binary. At each iteration, we consider all possible swap directions and choose the endpoint with the smallest objective value of \eqref{eq:sign-miqp}. The detail swap rounding procedure can be found in \cite{harvey2014pipage}. Each iteration makes at least one fractional coordinate binary. Since the cardinality $k-\upsilon$ is integral, the procedure performs at most $|\I|-1\le r-1$ iterations. The resulting \Cref{alg:pipage} differs from \Cref{alg:round} only in how the binary vector is obtained.

\begin{algorithm}[H]
\caption{Binary Redistribution via  Pipage Rounding}
\label{alg:pipage}
\begin{algorithmic}[1]
\State \textbf{Input:} An optimal solution $\beta^*$ returned by \Cref{alg:sparsity_reduction}, $X$, $y$, $k$, and $z^{\mathrm{rel}}$ 
\State Let $\hat \beta =\beta^*$
\If{\(\|\beta^*\|_0\geq k+1\)}
\State Compute $\Delta = \Delta(\beta^*)$ and $\I=\I(\beta^*)$
\State Let $\sigma_i=\sign(\beta^*_i), \ \forall i\in\I$

\State Apply pipage rounding to $z^{\mathrm{rel}}$ and denote the output by $z^{\mathrm{PR}}=(z_i^{\mathrm{PR}})_{i\in \I}$ 

\For{$i\in \I$}
\State Update $\hat \beta_i\gets \sigma_i\Delta z_i^{\mathrm{PR}}$
\EndFor
\EndIf
\State \textbf{return} \(\hat \beta\)
\end{algorithmic}
\end{algorithm}

\subsubsection{Performance Guarantee}
Since \Cref{alg:round} dominates \Cref{alg:pipage}, it suffices to establish the approximation ratio for \Cref{alg:pipage}. A common approach to analyzing pipage rounding is to construct an auxiliary function that is concave under swaps \cite{calinescu2011maximizing,harvey2014pipage}.  Next, we introduce such a function tailored to  \eqref{eq:sign-miqp}.

\begin{lemma}\label{lem:concave-swap}
For any $z=(z_i)_{i\in\I}\in[0,1]^{|\I|}$, define
\begin{align*}
    F(z)&:=\bigg\|y-\sum_{i\notin \I} X_i\beta_i^* -\sum_{i\in \I} \Delta\sigma_i X_i z_i\bigg\|_2^2+\Delta^2\lambda_{\max}(X^{\top}X)\sum_{i\in \I}z_i(1-z_i).
\end{align*}
The function $F(z)$ is concave under swaps.
\end{lemma}
\begin{proof}
For any $z\in[0,1]^{|\I|}$ and any distinct $i,j\in \I$, let
\[
    s:=\begin{bmatrix}
        \sigma_i \\
        -\sigma_j
    \end{bmatrix}, \quad
    A:=\begin{bmatrix}
       X_i \; X_j 
    \end{bmatrix}.
\]
Also, we let $z(\delta)$ be obtained from $z$ by replacing $(z_i,z_j)$ with $(z_i+\delta,z_j-\delta)$, and set $h(\delta):=F(z(\delta))$. It is easy to verify that $h(\delta)$ is concave in $\delta$ because the second-order derivative
\begin{align*}
    &\frac{d^2h(\delta)}{d\delta^2}=2\Delta^2\left[\|\sigma_i X_i-\sigma_j X_j\|_2^2-2\lambda_{\max}(X^{\top}X)\right]=2\Delta^2\left[s^\top A^\top As-2\lambda_{\max}(X^{\top}X)\right] \\
    &\le4\Delta^2\left[\lambda_{\max}(A^\top A)-\lambda_{\max}(X^{\top}X)\right]\le4\Delta^2\left[\lambda_{\max}(X^\top X)-\lambda_{\max}(X^{\top}X)\right]=0,
\end{align*}
where the first inequality uses Rayleigh quotient inequality and $\sigma_i,\sigma_j\in\{-1,1\}$, and the last inequality follows from Theorem 8.1.7 in \cite{golub2013matrix}.
\qed
\end{proof}

By leveraging \Cref{lem:phi-mont,lem:concave-swap} and Theorem 3.2 in \cite{harvey2014pipage}, we derive the same approximation ratio of \Cref{alg:round,alg:pipage}.

\begin{theorem}\label{them:binary-ratio}
Let $\hat \beta \in \Re^p$ be returned by \Cref{alg:round} or \Cref{alg:pipage}. Then $\hat\beta$ is a feasible solution for \eqref{eq:srr} and  satisfies 
\begin{align*}
{\frac{1}{n}\|y-X\hat \beta\|_2^2+\lambda\|\hat\beta\|_2^2}
 \leq \left[1+\frac{\lambda_{\max}(X^{\top}X)}{\lambda n}
     \min\left\{\frac{r-1}{k+r-1},\frac{r}{4k}\right\}\right] \V_{\mathrm{opt}}(k).
\end{align*}
\end{theorem}
\begin{proof}
If $\|\beta^*\|_0\le k$, the algorithm
returns $\beta^*$, which is optimal for \eqref{eq:srr} by
\Cref{lem:phi_k}. The result follows immediately.
We therefore consider the case $\|\beta^*\|_0\ge k+1$. According to \Cref{cor:inv}, $\hat \beta$ is feasible and equal to $\beta(z^{\mathrm{PR}})$ as defined in \eqref{eq:beta-round}. 

To prove the approximation ratio, we first bound the objective value of \eqref{eq:sign-miqp} at $z^{\mathrm{PR}}$. 

Let $F(z)$ be as in \Cref{lem:concave-swap}. For each selected pair of fractional coordinates, the
term $\sum_{i\in\I}z_i(1-z_i)$ takes the same value at the two
swap endpoints. Thus, choosing the endpoint with the smaller
objective value of \eqref{eq:sign-miqp} also minimizes $F$
over the two endpoints.
By \Cref{lem:concave-swap} and \cite[Theorem~3.2]{harvey2014pipage}, it follows that
\[
   F(z^\mathrm{PR})\le F(z^\mathrm{rel}).
\]
Since $z^\mathrm{PR}$ is binary, the above inequality reduces to
\[
    \|y-X\beta(z^\mathrm{PR})\|_2^2-\|y-X\beta(z^{\mathrm{rel}})\|_2^2\le \lambda_{\max}(X^{\top}X)\Delta^2\sum_{i\in \I}z^\mathrm{rel}_i(1-z^{\mathrm{rel}}_i).
\]

Leveraging the above results, the objective value of \eqref{eq:srr} at $\beta(z^{\mathrm{PR}})$ satisfies 
\begin{align*}
    &\frac{\frac{1}{n}\|y-X\beta(z^{\mathrm{PR}})\|_2^2+\lambda\|\beta(z^{\mathrm{PR}})\|_2^2}{\V_{\mathrm{rel}}(k)}=\frac{\frac{1}{n}\|y-X\beta(z^{\mathrm{PR}})\|_2^2+\lambda\phi(\beta^*;k)}{\frac{1}{n}\|y-X\beta^*\|_2^2+\lambda\phi(\beta^*;k)} \\
    =&1+\frac{1}{n}\frac{\|y-X\beta(z^{\mathrm{PR}})\|_2^2-\|y-X\beta^*\|_2^2}{\frac{1}{n}\|y-X\beta^*\|_2^2+\lambda\phi(\beta^*;k)}\le1+\frac{1}{\lambda n}\frac{\|y-X\beta(z^{\mathrm{PR}})\|_2^2-\|y-X\beta^*\|_2^2}{\phi(\beta^*;k)} \\
    \le&1+\frac{\lambda_{\max}(X^{\top}X)}{\lambda n}\frac{\Delta^2\sum_{i\in \I}z^\mathrm{rel}_i(1-z^{\mathrm{rel}}_i)}{\phi(\beta^*;k)},
\end{align*}
where the first equality follows from \Cref{cor:inv}. Rearranging the above inequality yields
\[
    \frac{1}{n}\|y-X\beta(z^{\mathrm{PR}})\|_2^2+\lambda\|\beta(z^{\mathrm{PR}})\|_2^2\le \left[1+\frac{\lambda_{\max}(X^{\top}X)}{\lambda n}\frac{\Delta^2\sum_{i\in \I}z^\mathrm{rel}_i(1-z^{\mathrm{rel}}_i)}{\phi(\beta^*;k)}\right]\V_{\mathrm{rel}}(k).
\]

We next establish two upper bounds on the term
\[
\frac{\Delta^2\sum_{i\in \I}z^\mathrm{rel}_i(1-z^{\mathrm{rel}}_i)}{\phi(\beta^*;k)}:
\]
\begin{enumerate}[(i)]
    \item Since $|\I|\le r$, it is easy to verify that 
    \[
    \sum_{i\in \I}z^\mathrm{rel}_i(1-z^{\mathrm{rel}}_i)\le\sum_{i\in \I}\frac{1}{4}=\frac{|\I|}{4}\le\frac{r}{4}.
    \]
    Also, \Cref{def:knorm} gives $\phi(\beta^*;k)\ge k\Delta^2$. Consequently,
    \[
        \frac{\Delta^2\sum_{i\in \I}z^\mathrm{rel}_i(1-z^{\mathrm{rel}}_i)}{\phi(\beta^*;k)}\le\frac{r}{4k}.
    \]
    \item For the second bound, we rewrite the term $\Delta^2\sum_{i\in \I}z^\mathrm{rel}_i(1-z^{\mathrm{rel}}_i)$ as
    \begin{align*}
        &\Delta^2\sum_{i\in \I}z^\mathrm{rel}_i(1-z^{\mathrm{rel}}_i) =\Delta^2\sum_{i\in \I}z^{\mathrm{rel}}_i-\sum_{i\in \I}(\Delta z^{\mathrm{rel}}_i)^2 \\
        &=(k-\upsilon)\Delta^2-\sum_{i\in \I}(\beta^*_i)^2 =\phi(\beta^*;k)-\|\beta^*\|_2^2,
    \end{align*}
    where the second equality uses \eqref{eq:z-rel}, and the last equality follows from Part (i) of \Cref{lem:upsilon}. Let $q:=\min\{p,r+k-1\}$. Since $\|\beta^*\|_0\le q$, \Cref{lem:phi_k} gives $\|\beta^*\|_2^2=\phi(\beta^*;q)$. Applying \Cref{lem:phi-mont} with $k_1=k$ and $k_2=q$, we can show that
    \[
    \frac{\|\beta^*\|_2^2}{\phi(\beta^*;k)}=\frac{\phi(\beta^*;q)}{\phi(\beta^*;k)}\ge \frac{k}{q}\ge \frac{k}{r+k-1}.
    \]
    Hence, we get
    \[
    \frac{\Delta^2\sum_{i\in \I}z^\mathrm{rel}_i(1-z^{\mathrm{rel}}_i)}{\phi(\beta^*;k)}=\frac{\phi(\beta^*;k)-\|\beta^*\|_2^2}{\phi(\beta^*;k)}\le\frac{r-1}{r+k-1}.
    \]
\end{enumerate}

Finally, since $\V_{\mathrm{rel}}(k)\le \V_{\mathrm{opt}}(k)$, the result follows immediately. We thus complete the proof.
 \qed
\end{proof}

As seen in \Cref{them:binary-ratio}, the approximation ratio decreases as $k$ or $\lambda$ increases. 
These trends are consistent with the numerical results in \Cref{sec:num}.

Based on the proof of \Cref{them:binary-ratio}, since the output of \Cref{alg:round} is feasible for \eqref{eq:srr}, the following corollary is immediate.

\begin{corollary}\label{cor:rel-gap}
We have
\begin{align*}
 \V_{\mathrm{opt}}(k)
 \leq \left[1+\frac{\lambda_{\max}(X^{\top}X)}{\lambda n}
     \min\left\{\frac{r-1}{k+r-1},\frac{r}{4k}\right\}\right] \V_{\mathrm{rel}}(k).
 \label{eq:miqp-gap-worst}
\end{align*}
\end{corollary}

\section{Numerical Experiments}\label{sec:num} 
In this section, we conduct numerical experiments on both synthetic and real data to evaluate the performance of our proposed algorithms.  

\textbf{Baselines.} We compare our \Cref{alg:round,alg:pipage,alg:reduce} with four baseline approximation algorithms that admit theoretical performance guarantees: the Greedy algorithm and the Randomized algorithm studied in \cite{xie2020scalable}, the Beam Search algorithm  \cite{liu2023okridge}, and the Dual IHT (DIHT) algorithm \cite{yuan2020dual}. 
We present the main comparisons in figures.
The detailed numerical results for all algorithms are presented in the tables of Appendix~\ref{app:tables}. We omit \Cref{alg:reduce} from the figures because its requirement $k\ge r$ limits applicability and it offers
no numerical advantage compared to \Cref{alg:round,alg:pipage}. Its results remain available in Appendix~\ref{app:tables}.

\textbf{Algorithm implementation.} We solve the convex relaxation \eqref{eq:rel} using the proximal gradient method \cite{mcdonald2014spectral}.
Since Greedy, Randomized, and Beam Search return supports,  for each solution returned by \Cref{alg:round,alg:pipage,alg:reduce} or DIHT, we extract its support and then compute the optimal value of the restricted ridge regression problem. This ensures a consistent comparison and may improve the objective value.
For each instance, we run the randomized algorithm 10 times and report the best result obtained. For  Beam Search, we set the beam width to $10$. In addition, we impose each algorithm a time limit of $600$ seconds per instance. If an algorithm reaches this limit, it returns the best feasible solution found.

\textbf{Hardware \& Software.} All the experiments are coded in Python 3.11 with calls to Gurobi 12.0.0 on a PC with 3.3 GHz Intel Core i7 processor and 16GB of memory.

\subsection{Experiments on Synthetic Data}\label{subsec:synthetic}
In this subsection, we conduct numerical experiments on synthetic data generated following  \cite{xie2020scalable}. We note that $X$ is full row rank, so $r=n$. We focus on the low-rank and high-dimensional regime where $r\ll p$. In particular, we set $r\in\{20,40,60\}$, $p\in\{1{,}000,3{,}000,5{,}000\}$, $k\in\{20,30,40,50,60,100,200\}$, and $\lambda\in\{0.1,0.5,1.0\}$.
For each parameter setting $(r,p,k,\lambda)$, we generate $10$ instances using different random seeds and report the average results. Since \Cref{alg:reduce} requires $k\ge r$, we evaluate its performance only for $k\in\{100,200\}$. 

\begin{figure}[!htbp]
\vspace{-1em}
    \centering
\includegraphics[width=0.9\textwidth]{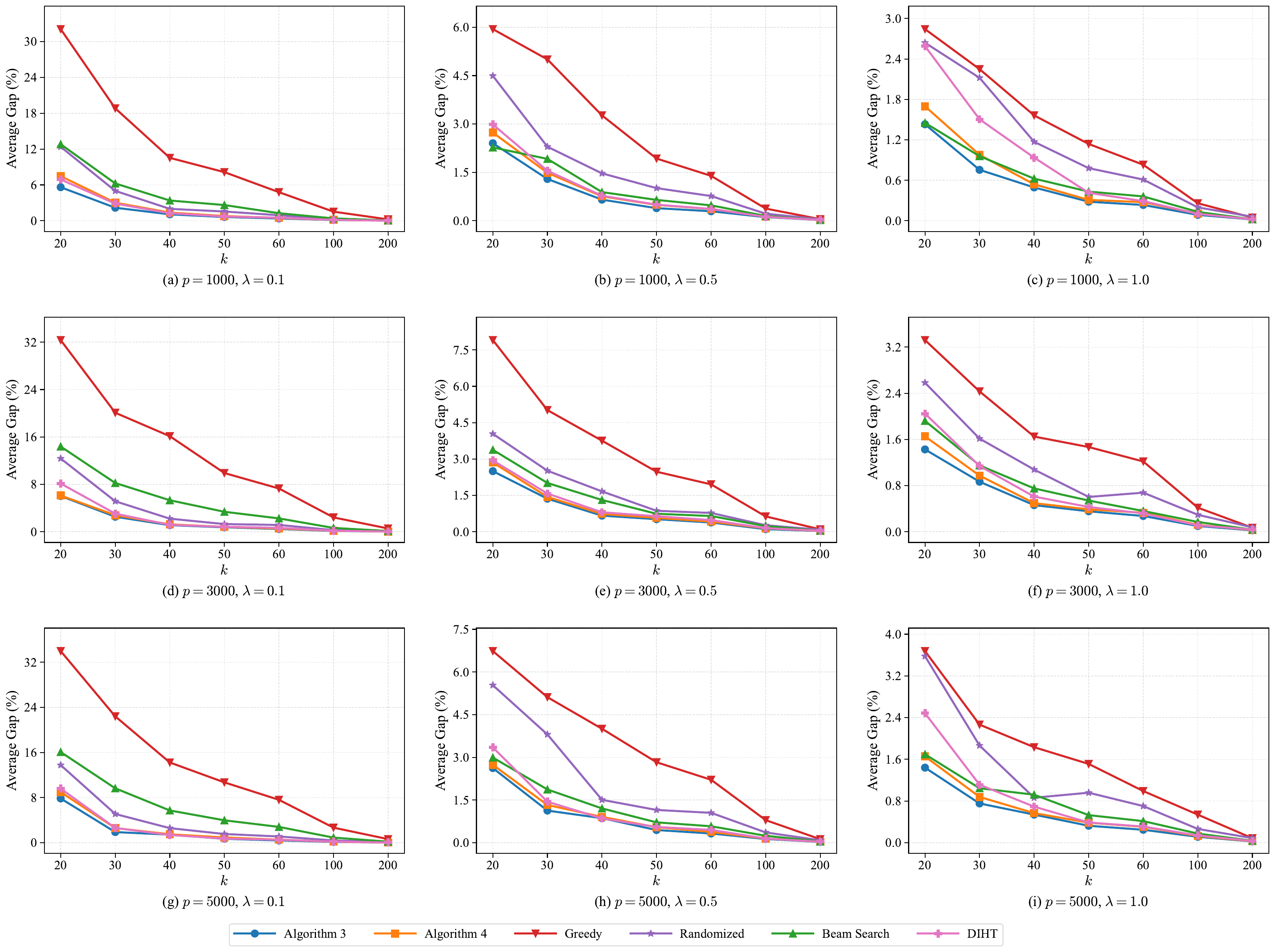}
\vspace{-1em}
    \caption{Average gaps on synthetic data with $r=20$.}
\label{fig:synthetic_r20}
\end{figure}

\begin{figure}[!htbp]
    \centering    \includegraphics[width=0.9\textwidth]{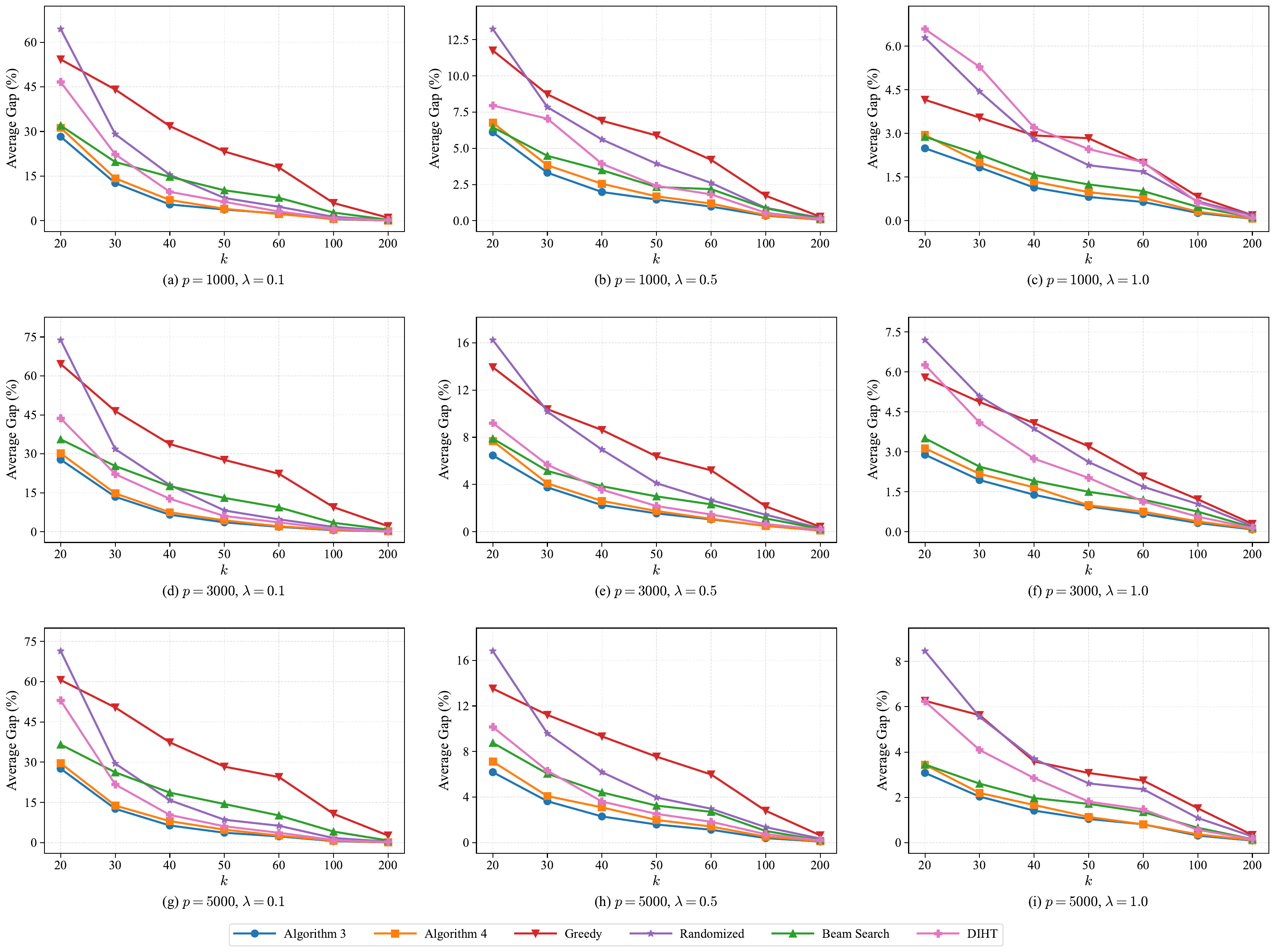}
    \vspace{-1em}
    \caption{Average gaps on synthetic data with $r=40$.}
\label{fig:synthetic_r40}
\end{figure}

\begin{figure}[!htbp]
    \centering
\includegraphics[width=0.9\textwidth]{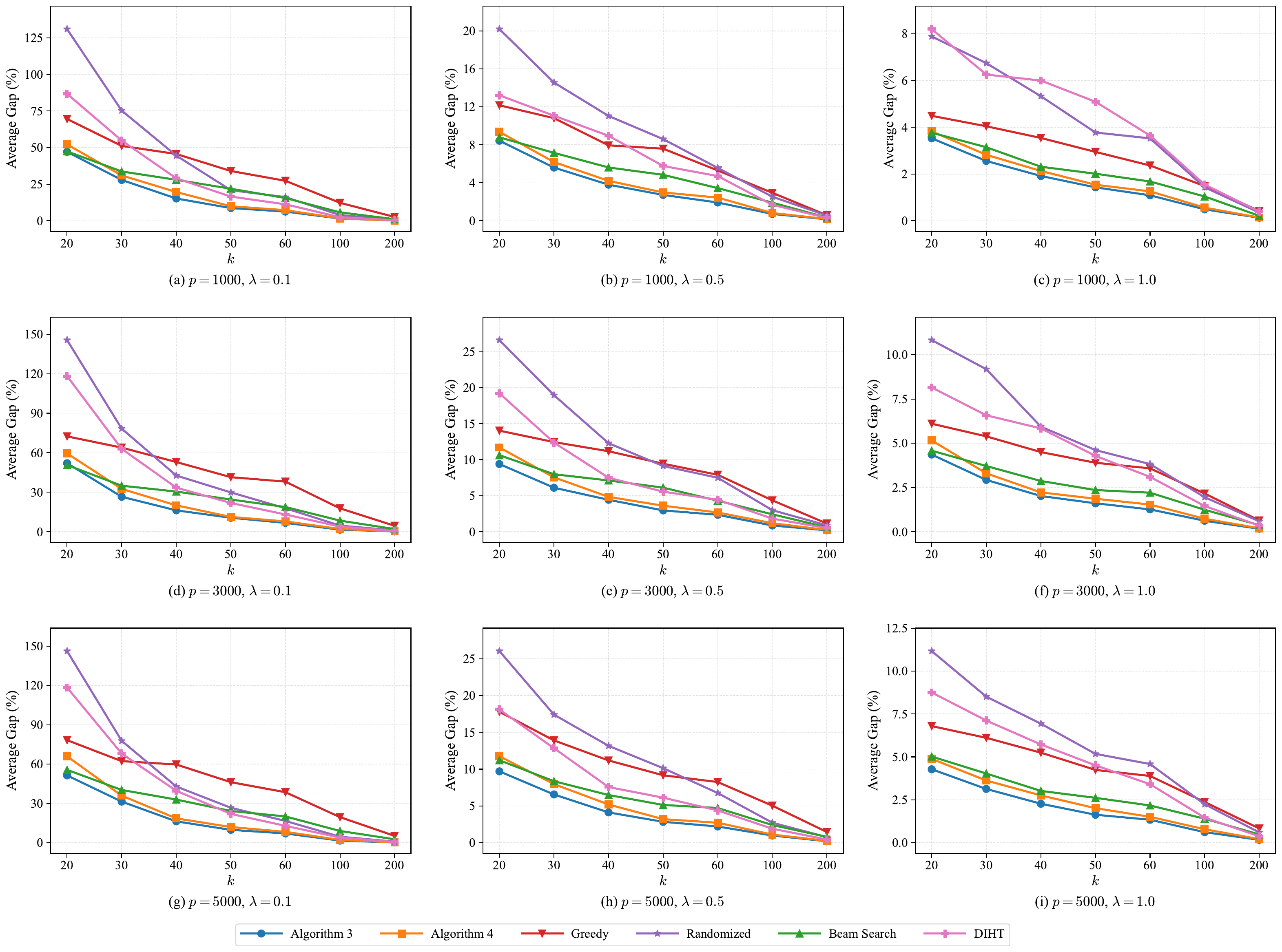}
\vspace{-1em}
    \caption{Average gaps on synthetic data with $r=60$.}
\label{fig:synthetic_r60}
\end{figure}

\Cref{fig:synthetic_r20,fig:synthetic_r40,fig:synthetic_r60} present the average gaps for $r=20,40,60$, respectively. Since the optimal value $\V_{\mathrm{opt}}(k)$ is often difficult to obtain, we instead use $\V_{\mathrm{rel}}(k)$ and compute the gaps as
\[
 \mathrm{Gap}(\%):=\frac{\mathrm{upper\ bound}-\V_{\mathrm{rel}}(k)}{\V_{\mathrm{rel}}(k)}\times 100\%,
\]
where the upper bound is the objective value returned by each approximation algorithm. Across all $189$ parameter settings, \Cref{alg:round,alg:pipage} achieve smaller gaps than Greedy and Randomized. They also outperform DIHT in $187$ and $171$ settings, respectively, and Beam Search in $186$ and $174$ settings, respectively.
Although \Cref{alg:pipage} only approximately solves problem \eqref{eq:sign-miqp}, its gaps are generally only slightly larger than those of \Cref{alg:round}. In a few instances, \Cref{alg:pipage} achieves even smaller gaps, which can be attributed to the final refit on the selected supports. 
\Cref{tab:average_time_synthetic} reports the average running times in seconds for each rank $r$. For our \Cref{alg:round,alg:pipage}, these times include solving \eqref{eq:rel} and applying \Cref{alg:sparsity_reduction} to obtain a solution satisfying the sparsity bound.
As seen, \Cref{alg:pipage} runs substantially faster than Randomized and Beam Search. It significantly reduces the running time of \Cref{alg:round} for larger ranks while maintaining comparable solution quality. Together with the gap comparisons, these results demonstrate that \Cref{alg:pipage} offers a favorable balance between solution quality and computational efficiency.

Next, we would like to analyze the effect of parameters $k$, $\lambda$, $r$, and $p$ on solution quality. 
As shown in \Cref{fig:synthetic_r20,fig:synthetic_r40,fig:synthetic_r60}, the average gaps are more sensitive to the rank $r$ than to the dimension $p$. This suggests that rank has a greater effect on algorithms of \eqref{eq:srr}.
In addition,
the average gaps of \Cref{alg:round,alg:pipage} decrease monotonically as $k$ increases.  As expected, increasing $\lambda$ also reduces the average gaps.  These trends are consistent with the dependence of the approximation ratio in \Cref{them:binary-ratio} on $k$ and $\lambda$.

\begin{table}[!htbp]
\caption{Average running times (seconds) on synthetic data.}
\label{tab:average_time_synthetic}
\centering
\begingroup
\renewcommand{\arraystretch}{1}
\setlength{\tabcolsep}{2.2pt}
\begin{tabular}{ccccccc}
    \toprule
    $r$ & \Cref{alg:round} & \Cref{alg:pipage} & Greedy & Randomized
        & DIHT & Beam Search \\
    \midrule
    $20$ & 0.134 & 0.104 & 0.007 & 0.369 & 0.064 & 0.424 \\
    $40$ & 0.231 & 0.086 & 0.012 & 0.457 & 0.077 & 0.356 \\
    $60$ & 2.799 & 0.089 & 0.017 & 0.545 & 0.097 & 0.390 \\
    \bottomrule
\end{tabular}%
\endgroup
\end{table}

\begin{figure}[!htbp]
    \centering    \includegraphics[width=0.9\textwidth]{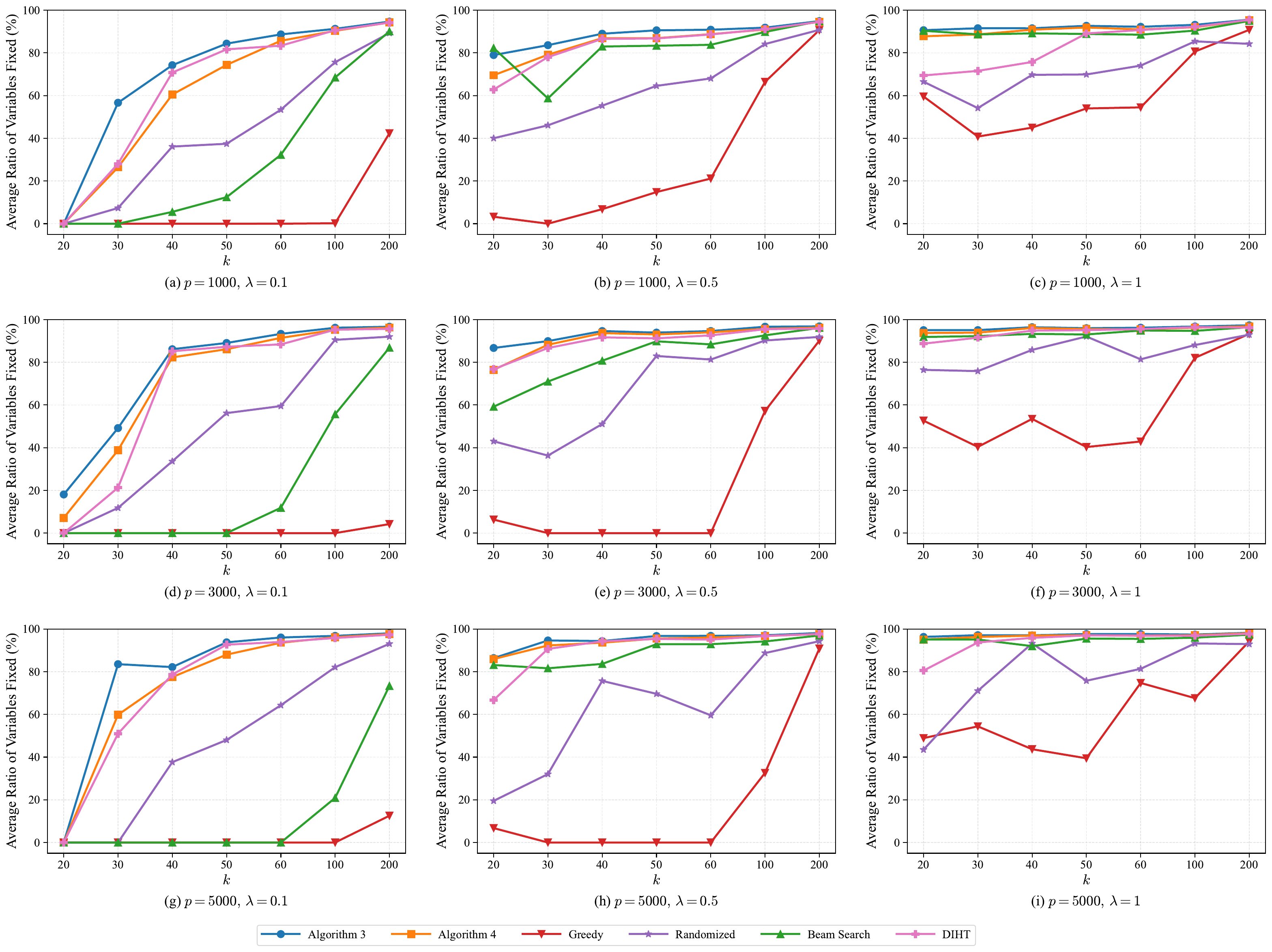}
    \vspace{-1em}
    \caption{Average ratios of variables fixed on synthetic data with $r=20$.}
\label{fig:screening_r20}
\end{figure}

\begin{figure}[!htbp]
    \centering
    \includegraphics[width=0.9\textwidth]{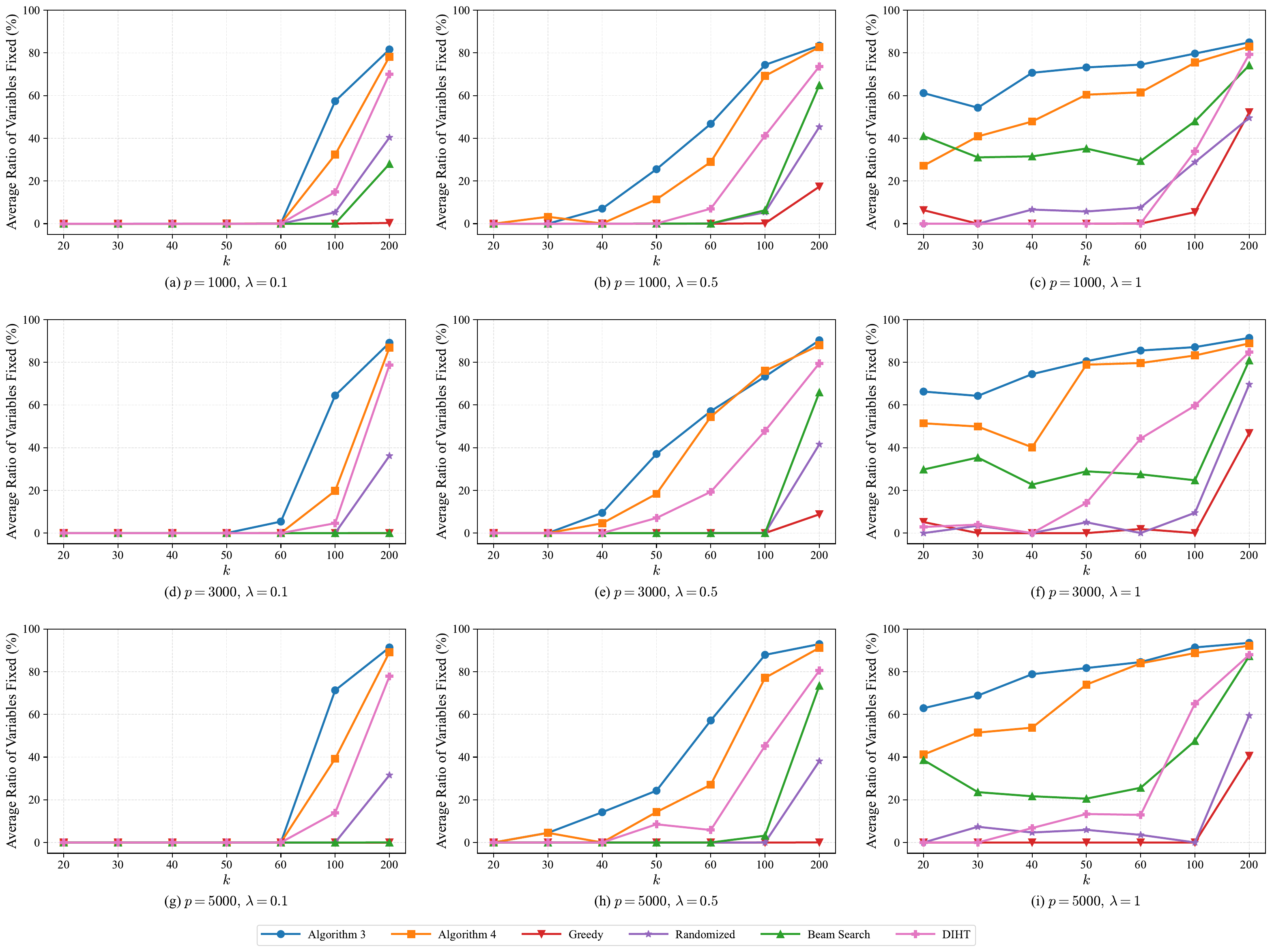}
    \vspace{-1em}
    \caption{Average ratios of variables fixed on synthetic data with $r=40$.}
    \label{fig:screening_r40}
\end{figure}

\begin{figure}[!htbp]
    \centering
\includegraphics[width=0.9\textwidth]{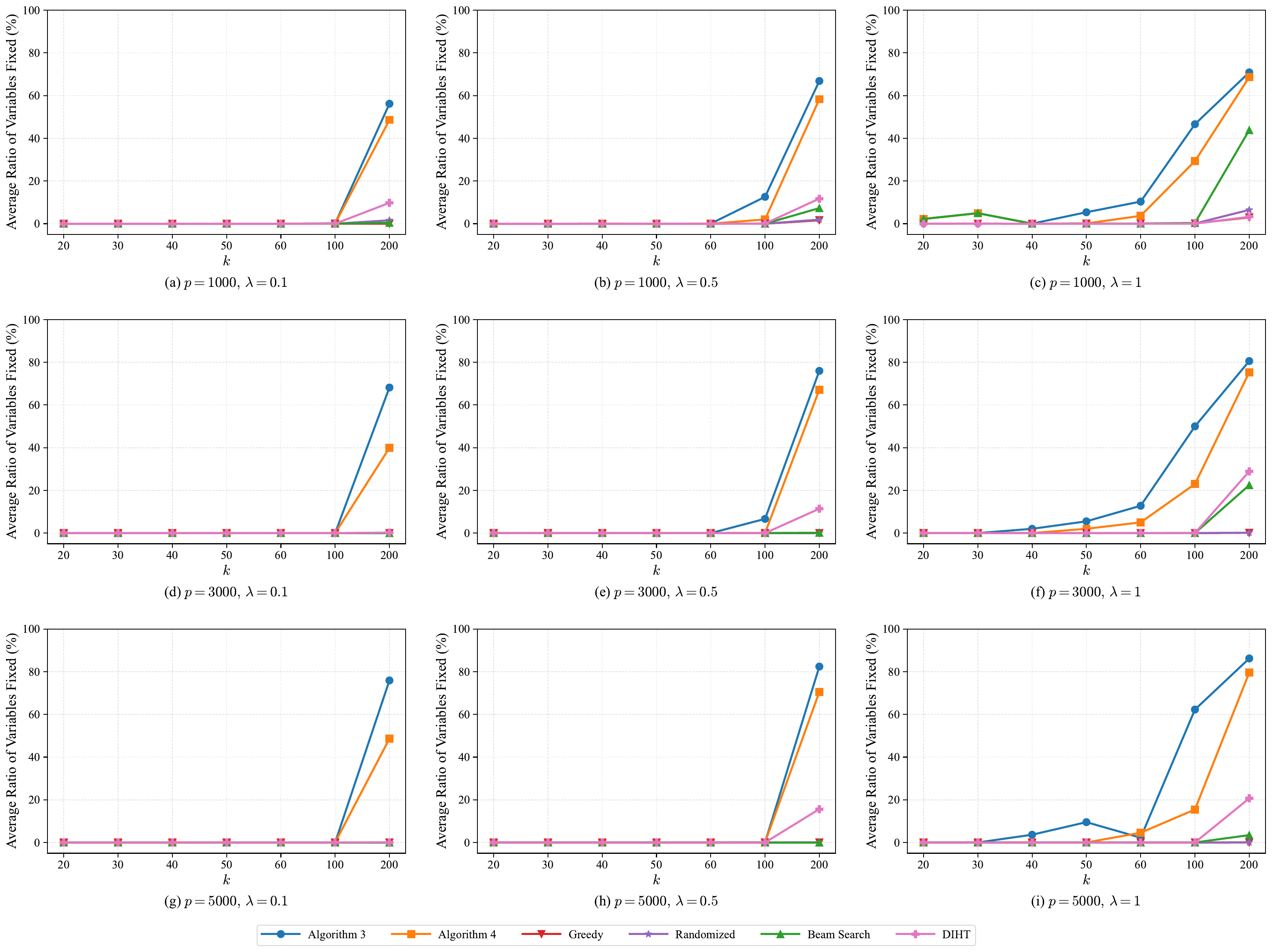}
\vspace{-1em}
    \caption{Average ratios of variables fixed on synthetic data with $r=60$.}
\label{fig:screening_r60}
\end{figure}

We further evaluate the variable fixing capacity of our algorithms by applying the fixing rules in \cite{atamturk2020safe}. A tighter upper bound can enable more variable fixed, which can reduce the problem size and help solve \eqref{eq:srr} to  optimality more efficiently. \Cref{fig:screening_r20,fig:screening_r40,fig:screening_r60} present the  proportions of the $p$ variables fixed by each algorithm. 
\Cref{alg:round,alg:pipage} consistently achieve higher average fixing ratios than the baselines in most parameter settings.  Even small improvements in gaps can lead to great improvements in variable fixing. For example, at $(r,p,k,\lambda)=(40,5{,}000,100,0.5)$, their average gaps are $0.402\%$ and $0.532\%$, compared with $0.748\%$ for DIHT and $1.030\%$ for Beam Search. In contrast, their fixing ratios are $87.884\%$ and $77.128\%$, whereas those of DIHT and Beam Search are only $45.296\%$ and $3.190\%$, respectively. These results demonstrate the enhanced ability of our algorithms to fix variables for \eqref{eq:srr}. Finally, we note that, as both the gaps and fixing ratios are averaged, an algorithm may exhibit a larger average gap yet still achieve a higher average fixing ratio than the other algorithm under the same parameter setting. 

\subsection{Experiments on Real Data}\label{subsec:real}
In this subsection, we conduct numerical experiments on five real-world datasets: 
Breast Cancer ($n=78$, $p=4{,}751$), Lung Cancer ($n=86$, $p=7{,}129$), and AML ($n=116$, $p=6{,}283$), all from \cite{bair2004semi}; Riboflavin ($n=71$, $p=4{,}088$) from \cite{buhlmann2014high}; and Leukemia ($n=72$, $p=3{,}571$) from \cite{dettling2004bagboosting}. For each dataset, the matrix $X$ has full row rank, so $r=n$. 
Then we set $k \in \{20,30,40,50,60,100,200\}$, except for the {AML} dataset, for which we replace $k=100$ with $k=150$ to satisfy the condition $k\ge r$ required by \Cref{alg:reduce}. We set $\lambda\in\{0.1,0.5,1.0\}$, with all three values within the range considered in \cite{atamturk2020safe}. 

\begin{figure}[!htbp]
    \centering
    \includegraphics[width=0.9\textwidth]{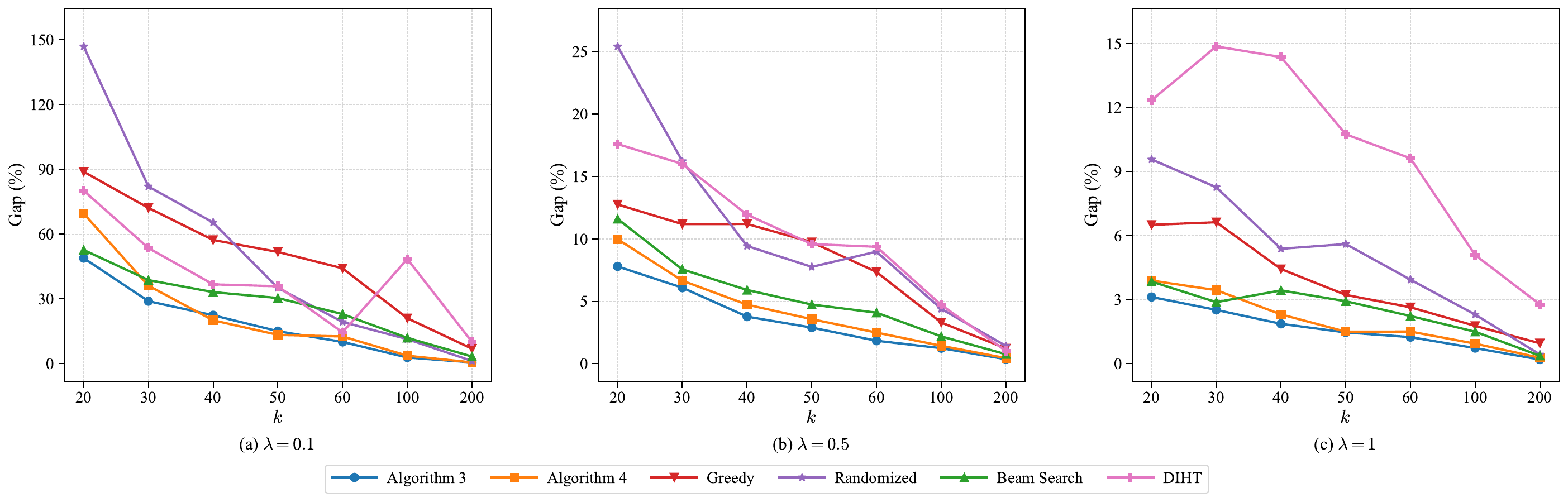}
    \vspace{-1em}
    \caption{Gaps with \textit{Breast Cancer} dataset.}
    \label{fig:real_Breast}
\end{figure}

\begin{figure}[!htbp]
    \centering
    \includegraphics[width=0.9\textwidth]{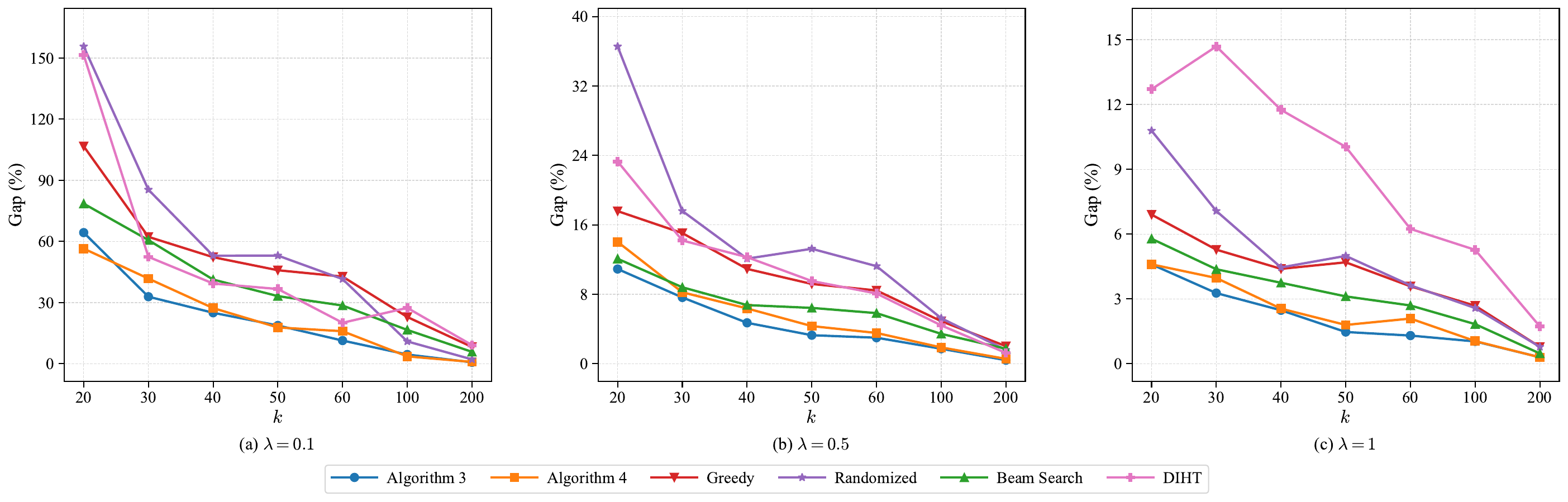}
    \vspace{-1em}
    \caption{Gaps with \textit{Lung Cancer} dataset.}
    \label{fig:real_Lung}
\end{figure}

\begin{figure}[!htbp]
    \centering
    \includegraphics[width=0.9\textwidth]{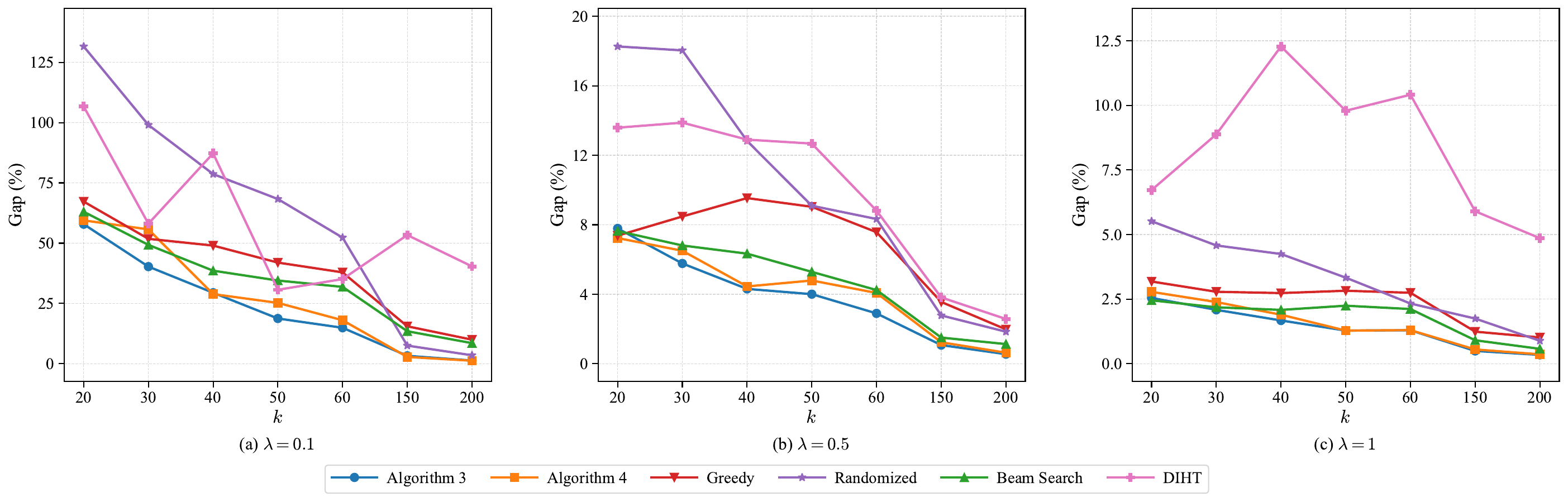}
    \vspace{-1em}
    \caption{Gaps with \textit{AML} dataset.}
    \label{fig:real_AML}
\end{figure}

\begin{figure}[!htbp]
    \centering
    \includegraphics[width=0.9\textwidth]{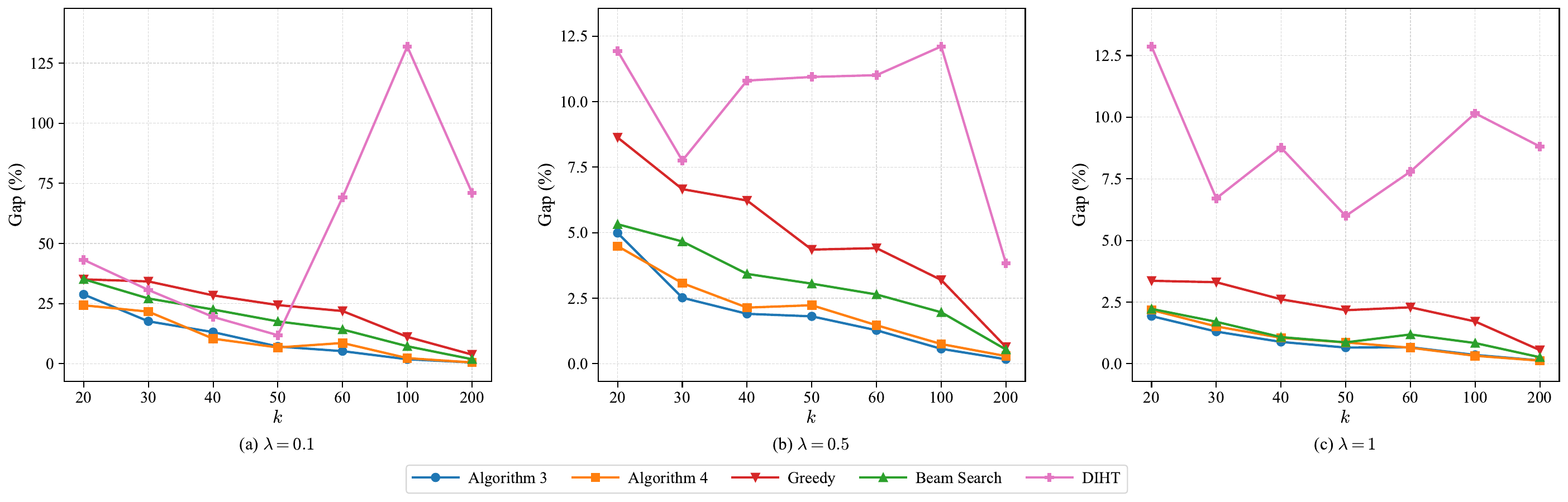}
    \vspace{-1em}
    \caption{Gaps with \textit{Riboflavin} dataset.}
    \label{fig:real_Riboflavin}
\end{figure}

\begin{figure}[!htbp]
    \centering
    \includegraphics[width=0.9\textwidth]{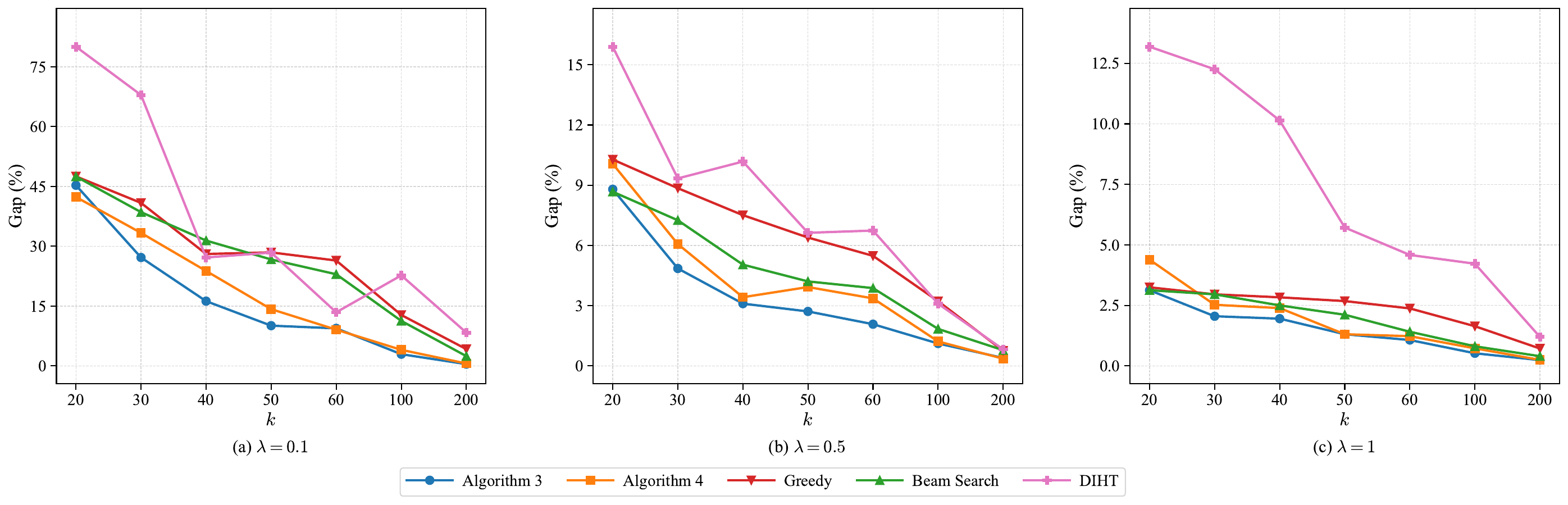}
    \vspace{-1em}
    \caption{Gaps with \textit{Leukemia} dataset.}
    \label{fig:real_leukemia}
\end{figure}

\Cref{fig:real_Breast,fig:real_Lung,fig:real_AML,fig:real_Riboflavin,fig:real_leukemia} display the gaps of each algorithm for each dataset, respectively.
Among the $105$ tested instances, \Cref{alg:round,alg:pipage} achieve smaller gaps than all four baselines in $102$ and $96$ cases, respectively. It is noted that the gaps of DIHT fluctuate substantially as $k$ varies, whereas \Cref{alg:round,alg:pipage} exhibit much more stable performance across different values of $k$. Besides, \Cref{tab:real_time} summarizes average running times across all instances for each dataset. The running time of \Cref{alg:round} generally increases with $r$ and is highest for AML with the largest rank. Its running times on Riboflavin and Leukemia also differ despite their similar ranks. This is because the size $|\I|$ can differ between the two datasets even for the same $k$ and $\lambda$ (see \Cref{tab:real_gap_time} in Appendix \ref{app:tables}). Leukemia has a larger $|\I|$, which makes problem \eqref{eq:sign-miqp} harder to solve and thus takes longer. In addition,  \Cref{alg:pipage} is substantially faster than \Cref{alg:round}, with average running times below $1$ second. It is also generally faster than Randomized and Beam Search. 
For these real instances, \Cref{alg:pipage} overall provides an efficient alternative to \Cref{alg:round} that preserves competitive solution quality.

\begin{table}[!htbp]
\centering
\caption{Average times across all parameter settings for each dataset.}
\label{tab:real_time}
\begingroup
\renewcommand{\arraystretch}{1}
\setlength{\tabcolsep}{2.2pt}
\begin{tabular}{ccccccc}
    \toprule
    Dataset & \Cref{alg:round} & \Cref{alg:pipage} & Greedy
            & Randomized & DIHT & Beam Search \\
    \midrule
    Breast Cancer & 97.114 & 0.505 & 0.033 & 1.005 & 0.187 & 0.541 \\
    Lung Cancer & 215.449 & 0.878 & 0.067 & 1.572 & 0.340 & 0.780 \\
    AML & 384.638 & 0.685 & 0.093 & 1.613 & 0.369 & 0.842 \\
    Riboflavin & 6.314 & 0.478 & 0.027 & 0.848 & 0.155 & 0.476 \\
    Leukemia & 23.533 & 0.247 & 0.026 & 0.829 & 0.133 & 0.443 \\
    \bottomrule
\end{tabular}%
\endgroup
\end{table}

\section{Conclusion}\label{sec:conclude}
We develop two relax-and-sparsify strategies for sparse ridge regression. Parameter reduction admits a data-independent approximation ratio, and binary redistribution outperforms the benchmarks in most numerical experiments. Both strategies rely only on the low rank  of $X$ and properties of $\phi$, rather than the specific least-squares loss. Hence, they can be extended to other sparse learning problems, such as sparse portfolio selection and sparse logistic regression, although the approximation ratios may not transfer.

\newpage
{\let\clearpage\relax
\vspace{10pt}
\renewcommand{\bibname}{\protect\leftline{References }
	
}
\bibliography{Reference}
}
\newpage
\appendix
\renewcommand*{\theHsection}{appendix.\Alph{section}}

\section{Proofs}
\subsection{Proof of \Cref{them:NP}}\label{proof:themNP}
\themNP*
\begin{proof}
In order to prove the result, we reduce the NP-hard decision version of \textit{Exact $k$-item Subset Sum Problem} (E-$k$SSP) (see \cite{caprara2000approximation}) to problem \eqref{eq:srr}. 

    \noindent\textit{E-$k$SSP instance.} Given $m\ge1$ positive integers $a_1,\dots,a_m$, a positive target integer $T$, and an integer $k\in[m]$, the E-$k$SSP asks whether there exists a subset $S\subseteq[m]$ such that 
    \[
        |S|=k, \quad\text{and} \quad \sum_{i\in S}a_i=T.
    \]

    \noindent\textit{Construction of the \eqref{eq:srr} instance.} Set $n:=r$ and $p:=m+r-1$, and consider the instance:
    \[
        X:=\begin{bmatrix}
            a_1 & \dots & a_m & 1 & 0_{r-2}^\top \\
            \frac{2T}{k}-a_1 & \dots & \frac{2T}{k}-a_m & -1 & 0_{r-2}^\top \\
            0_{r-2} & \dots & 0_{r-2} & 0_{r-2} & I_{r-2} 
        \end{bmatrix}
        \in\Re^{r\times p},
        \quad
        y:=\begin{bmatrix}
            1 \\ 1 \\ 0_{r-2}
        \end{bmatrix}
        \in\Re^r,
    \]
    where the last block row and block column of $X$ and the zero block of $y$ are omitted when $r=2$. 
    By \cite[Proposition 5]{xie2020scalable}, the problem \eqref{eq:srr} is equivalent to
    \begin{equation}\label{eq:C-opt-exact}
        \min_{z\in\{0,1\}^p}
        \bigg\{
            \lambda y^\top M(z)^{-1}y \ : \
            \sum_{i\in [p]} z_i\le k
        \bigg\},
        \quad M(z):=r\lambda I_r+X\Diag(z)X^\top.
    \end{equation}

    We divide the remainder of the proof into three steps.
    
    \noindent\textbf{Step I:} compute $y^\top M(z)y$ explicitly. For any feasible solution $z$ of problem \eqref{eq:C-opt-exact}, we have
    \begin{align*}
        y^\top M(z)y&=r\lambda y^\top y+y^\top X\Diag(z)X^\top y=r\lambda y^\top y+\sum_{i\in[p]}z_iy^\top (X_iX_i^\top)y \\
        &=2\lambda r+\sum_{i\in[p]}z_i(X_i^\top y)^2=2\lambda r+\frac{4T^2}{k^2}\sum_{i\in[m]}z_i.
    \end{align*}
    Since $M(z)$ is positive definite, by the Cauchy--Schwarz inequality, the objective value of \eqref{eq:C-opt-exact} is bounded below by
    \begin{equation}\label{ineq:C-opt-obj}
        \lambda y^\top M(z)^{-1}y\ge\lambda\frac{(y^\top y)^2}{y^\top M(z)y}=\frac{4\lambda}{2\lambda r+\frac{4T^2}{k^2}\sum_{i\in[m]}z_i}\ge\frac{2\lambda}{\lambda r+\frac{2T^2}{k}}=\frac{2\lambda k}{\lambda rk+2T^2},
    \end{equation}
    where the last inequality follows from $\sum_{i\in[m]}z_i\le\sum_{i\in[p]}z_i\le k$. 
    
    We now characterize when both inequalities in \eqref{ineq:C-opt-obj} hold with equality. Equality in the first inequality holds if and only if $y$ is an eigenvector of $M(z)$. Since
    \begin{align*}
        M(z)y&=\lambda ry+X\Diag(z)X^\top y=\lambda ry+\sum_{i\in[p]}z_iX_i(X_i^\top y)=\lambda ry+\sum_{i\in[m]}z_iX_i\frac{2T}{k} \\
        &=\lambda ry+\frac{2T}{k}\sum_{i\in[m]}z_i\begin{bmatrix}a_i \\ \frac{2T}{k}-a_i \\0_{r-2}\end{bmatrix}=\lambda ry+\frac{2T}{k}\begin{bmatrix}
            \sum_{i\in[m]}z_ia_i \\ \frac{2T}{k}\sum_{i\in[m]}z_i-\sum_{i\in[m]}z_ia_i \\0_{r-2}
        \end{bmatrix},
    \end{align*}
    this condition holds if and only if $\sum_{i\in[m]}z_ia_i=T/k\sum_{i\in[m]}z_i$. Equality in the second inequality of \eqref{ineq:C-opt-obj} holds if and only if $\sum_{i\in[m]}z_i=k$, which also implies that $z_i=0$ for all $i\in[m+1,p]$.

    \noindent\textbf{Step II:} exact subset sum case. For each $i\in [p]$, let $z_i\in\{0,1\}$ be an indicator variable such that 
    \[
        z_i=\begin{cases}
            1,\quad i\in S, \\
            0,\quad i\notin S.        
        \end{cases}
    \]
    If an exact subset exists, the sum of the chosen $a_i$ is $T$. Let $z^*$ be the solution corresponding to the exact subset sum. It follows that
    \[
        \sum_{i\in[m]}z_i^*a_i=T, \quad z_i^*=0,\ \forall i\in[m+1,p] \quad \text{and} \quad \sum_{i\in[m]}z_i^*=k.
    \]
    Obviously, $T/k\sum_{i\in[m]}z_i^*=T=\sum_{i\in[m]}z_i^*a_i$. Therefore, both inequalities in \eqref{ineq:C-opt-obj} hold with equality, and hence the optimal value $\V_{\mathrm{opt}}(k)$ of \eqref{eq:C-opt-exact} is 
    \[
        \V_{\mathrm{opt}}(k)=\frac{2\lambda k}{\lambda rk+2T^2}.
    \]

    \noindent\textbf{Step III:} non-exact subset sum case. Suppose that no exact subset sum exists. For any feasible solution $z$ of \eqref{eq:C-opt-exact}, we consider two cases. If $\sum_{i\in[m]}z_i<k$, the second inequality in \eqref{ineq:C-opt-obj} is strict, which gives 
    \[
        \lambda y^\top M(z)^{-1}y>\frac{2\lambda k}{\lambda rk+2T^2}.
    \]
    Otherwise, $\sum_{i\in[m]}z_i=k$. In this case, since $\sum_{i\in[m]}z_ia_i\neq T$, it is easy to check that $y$ cannot be an eigenvector of $M(z)$ based on the definition of eigenvectors. Thus, the first inequality in \eqref{ineq:C-opt-obj} is strict, which yields the same conclusion. Combining both cases, we can conclude that the optimal value of \eqref{eq:C-opt-exact} satisfies
    \[
        \V_{\mathrm{opt}}(k)>\frac{2\lambda k}{\lambda rk+2T^2}.
    \]

    Therefore, 
    \[
        \text{E-}k\text{SSP is YES} \ \Longleftrightarrow \ \text{problem \eqref{eq:C-opt-exact} admits an optimal value } \V_{\mathrm{opt}}(k)=\frac{2\lambda k}{\lambda rk+2T^2}.
    \]
    We have a polynomial-time reduction from E-$k$SSP to the decision version of \eqref{eq:C-opt-exact}. Hence, the decision version is NP-hard, and so \eqref{eq:C-opt-exact} is NP-hard. We thus complete the proof. \qed
\end{proof}

\subsection{Proof of \Cref{prop:tight}}\label{proof:proptight}
\proptight*
\begin{proof}
Set $n:=r$ and $p:=r+k-1$, and consider the instance:
\[
  X=\left[\,2\,1_r1_{k-1}^{\top}\ \middle|\ rIr\,\right]
    \in\Re^{r\times p},
    \quad
    y=1_r \in \Re^r,
\]
where the first block of $X$ is omitted when $k=1$.

According to \cite{mcdonald2014spectral}, we can rewrite \eqref{eq:rel} as
\begin{equation}\label{eq:socp}
      \min_{\beta\in\Re^p,\ z\in[0,1]^p}
    \bigg\{
        \frac{1}{r}\|y-X\beta\|_2^2
+\lambda\sum_{i\in [p]}\frac{\beta_i^2}{z_i} \ : \ \sum_{i\in [p]} z_i \le k 
    \bigg\}, 
\end{equation}
where $0/0=0$ and $\beta_i^2/0=+\infty$ if $\beta_i\ne0$.
For each fixed feasible  $z$,  the subproblem of \eqref{eq:socp} over $\beta$ has the unique minimizer (\cite{xie2020scalable}): 
\begin{equation*}
    \beta(z)=\Diag(z)X^\top M(z)^{-1}y,
    \quad
    M(z):=r\lambda I_r+X\Diag(z)X^\top.
\end{equation*} Substituting $\beta(z)$ into the objective of \eqref{eq:socp} yields the known
reformulation:
\begin{equation}\label{eq:C-opt}
    \min_{z\in[0,1]^p}
    \bigg\{
        \lambda y^\top M(z)^{-1}y \ : \
\sum_{i\in [p]} z_i\le k
    \bigg\}.
\end{equation}
Every optimal solution of \eqref{eq:rel} is of the form $\beta(z)$
for some optimal solution $z$ of \eqref{eq:C-opt} \cite{xie2020scalable}. Consequently, if \eqref{eq:C-opt} has a unique optimal solution $z^*$, then $\beta(z^*)$ must be unique and optimal
for \eqref{eq:rel}.

We first derive a lower bound  of
\eqref{eq:C-opt}. For the constructed instance,
\[
    M(z)
    =
    r\lambda I_r
    +4  \sum_{i\in [k-1]} z_i \cdot 1_r1_r^\top
    +r^2\Diag(z_{[k,p]}),
\]
and hence
\[    y^\top M(z)y
    =
    r^2\bigg(
        \lambda
        +4  \sum_{i\in [k-1]} z_i
        +  \sum_{i\in [k, p]} z_i
    \bigg)=r^2\bigg(
        \lambda
        +3  \sum_{i\in [k-1]} z_i
        +  \sum_{i\in [p]} z_i
    \bigg).
\]
Since $M(z)$ is positive definite, the Cauchy--Schwarz inequality yields
\begin{align}\label{ineq:obj}
     \lambda y^\top M(z)^{-1}y
    \ge      \lambda
    \frac{(y^{\top}y)^2}{y^\top M(z)y}
    =
    \frac{\lambda}{
  \lambda
        +3  \sum_{i\in [k-1]} z_i
        +  \sum_{i\in [p]} z_i}
    \ge
    \frac{\lambda}{\lambda+4k-3},   
\end{align}
where the  last inequality follows from  $\sum_{i\in[k-1]}z_i\le k-1$ and
$\sum_{i\in[p]}z_i\le k$.

We now characterize when both inequalities in
\eqref{ineq:obj} hold with equality. Equality in the first
inequality holds if and only if $y=1_r$ is an eigenvector of $M(z)$.
Since
\[
    M(z)y
    =
r\bigg(\lambda+4\sum_{i\in[k-1]}z_i\bigg)1_r
    +r^2z_{[k,p]},
\]
this condition is equivalent to all entries of $z_{[k,p]}$
being equal. Equality in the second inequality of \eqref{ineq:obj} holds if and only if $z_i=1$ for all $i\in [k-1]$ and $\sum_{i\in[p]}z_i=k$.
Together, these conditions uniquely determine a vector
\[
    z_i^*= 1, \ \forall i\in [k-1], \quad z_i^*=\frac 1 r, \ \forall i\in [k,p].
\]
This vector $z^*$ is feasible and attains the lower bound in
\eqref{ineq:obj}.  Hence,
$z^*$ is the unique optimal solution of \eqref{eq:C-opt}.
The corresponding unique optimal solution of \eqref{eq:rel}  is 
\[\beta^*:=\beta(z^*)
    =\Diag(z^*)X^\top M(z^*)^{-1}y
    =\frac{1}{\lambda+4k-3}
    \begin{bmatrix}
        2\,1_{k-1}\\[2pt]
        \frac{1}{r}\,1_r
    \end{bmatrix}.
\]
Since every entry of $\beta^*$ is
positive,
\[ \|\beta^*\|_0=p=r+k-1,
\]
which establishes the tightness of both sparsity bounds. \qed
\end{proof}

\section{Detailed Computational Results}
\label{app:tables}

\begin{table}[!htbp]
\caption{Average gaps and average times on synthetic data with $r=20$.}
\label{tab:synthetic_r20}
\centering
\begingroup
\scriptsize
\renewcommand{\arraystretch}{1.15}
\setlength{\tabcolsep}{2.2pt}
\resizebox{1.0\textwidth}{!}{%
\begin{tabular}{cccc*{7}{cc}}
    \toprule
    \multirow{2}{*}{$r$}
        & \multirow{2}{*}{$p$}
        & \multirow{2}{*}{$\lambda$}
        & \multirow{2}{*}{$k$}
        & \multicolumn{2}{c}{\Cref{alg:reduce}}
        & \multicolumn{2}{c}{\Cref{alg:round}}
        & \multicolumn{2}{c}{\Cref{alg:pipage}}
        & \multicolumn{2}{c}{Greedy}
        & \multicolumn{2}{c}{Randomized}
        & \multicolumn{2}{c}{DIHT}
        & \multicolumn{2}{c}{Beam Search} \\
    \cmidrule(lr){5-6}
    \cmidrule(lr){7-8}
    \cmidrule(lr){9-10}
    \cmidrule(lr){11-12}
    \cmidrule(lr){13-14}
    \cmidrule(lr){15-16}
    \cmidrule(lr){17-18}
    & & & &
    \makecell{Avg.\\Gap} & \makecell{Avg.\\Time} &
    \makecell{Avg.\\Gap} & \makecell{Avg.\\Time} &
    \makecell{Avg.\\Gap} & \makecell{Avg.\\Time} &
    \makecell{Avg.\\Gap} & \makecell{Avg.\\Time} &
    \makecell{Avg.\\Gap} & \makecell{Avg.\\Time} &
    \makecell{Avg.\\Gap} & \makecell{Avg.\\Time} &
    \makecell{Avg.\\Gap} & \makecell{Avg.\\Time} \\
    \midrule

    \multirow{63}{*}{$20$} & \multirow{21}{*}{$1{,}000$} & \multirow{7}{*}{$0.1$} & 20 & - & - & \textbf{5.612} & 0.059 & 7.466 & 0.030 & 32.092 & 0.001 & 12.380 & 0.108 & 6.912 & 0.048 & 12.813 & 0.054 \\
    & & & 30 & - & - & \textbf{2.178} & 0.065 & 3.028 & 0.032 & 18.823 & 0.001 & 4.999 & 0.123 & 2.870 & 0.021 & 6.244 & 0.078 \\
    & & & 40 & - & - & \textbf{1.070} & 0.050 & 1.369 & 0.027 & 10.531 & 0.002 & 1.990 & 0.104 & 1.256 & 0.022 & 3.386 & 0.101 \\
    & & & 50 & - & - & \textbf{0.623} & 0.062 & 0.788 & 0.031 & 8.140 & 0.002 & 1.546 & 0.119 & 0.702 & 0.022 & 2.629 & 0.132 \\
    & & & 60 & - & - & \textbf{0.381} & 0.061 & 0.470 & 0.034 & 4.782 & 0.002 & 0.921 & 0.119 & 0.515 & 0.027 & 1.258 & 0.156 \\
    & & & 100 & 0.355 & 0.039 & \textbf{0.124} & 0.065 & 0.139 & 0.033 & 1.518 & 0.003 & 0.295 & 0.128 & 0.133 & 0.022 & 0.390 & 0.270 \\
    & & & 200 & 0.061 & 0.056 & \textbf{0.023} & 0.062 & 0.026 & 0.032 & 0.227 & 0.006 & 0.054 & 0.125 & 0.027 & 0.023 & 0.047 & 0.532 \\
    \cmidrule(lr){3-18}

    & & \multirow{7}{*}{$0.5$} & 20 & - & - & 2.403 & 0.053 & 2.734 & 0.012 & 5.942 & 0.001 & 4.494 & 0.121 & 2.986 & 0.047 & \textbf{2.273} & 0.051 \\
    & & & 30 & - & - & \textbf{1.295} & 0.036 & 1.486 & 0.013 & 5.003 & 0.001 & 2.295 & 0.126 & 1.550 & 0.020 & 1.918 & 0.075 \\
    & & & 40 & - & - & \textbf{0.654} & 0.034 & 0.757 & 0.012 & 3.267 & 0.002 & 1.469 & 0.116 & 0.775 & 0.023 & 0.888 & 0.103 \\
    & & & 50 & - & - & \textbf{0.390} & 0.035 & 0.491 & 0.015 & 1.929 & 0.002 & 1.011 & 0.125 & 0.499 & 0.021 & 0.646 & 0.132 \\
    & & & 60 & - & - & \textbf{0.294} & 0.039 & 0.363 & 0.018 & 1.392 & 0.002 & 0.766 & 0.123 & 0.365 & 0.021 & 0.477 & 0.151 \\
    & & & 100 & 0.287 & 0.025 & \textbf{0.112} & 0.043 & 0.125 & 0.013 & 0.377 & 0.004 & 0.219 & 0.126 & 0.122 & 0.021 & 0.148 & 0.266 \\
    & & & 200 & 0.044 & 0.024 & 0.023 & 0.038 & 0.025 & 0.014 & 0.049 & 0.007 & 0.046 & 0.128 & 0.025 & 0.022 & \textbf{0.022} & 0.546 \\
    \cmidrule(lr){3-18}

    & & \multirow{7}{*}{$1.0$} & 20 & - & - & \textbf{1.431} & 0.037 & 1.698 & 0.007 & 2.843 & 0.001 & 2.641 & 0.121 & 2.590 & 0.028 & 1.452 & 0.056 \\
    & & & 30 & - & - & \textbf{0.755} & 0.031 & 0.978 & 0.009 & 2.249 & 0.001 & 2.120 & 0.109 & 1.505 & 0.021 & 0.958 & 0.076 \\
    & & & 40 & - & - & \textbf{0.494} & 0.030 & 0.540 & 0.009 & 1.563 & 0.002 & 1.171 & 0.115 & 0.936 & 0.022 & 0.625 & 0.098 \\
    & & & 50 & - & - & \textbf{0.282} & 0.028 & 0.309 & 0.010 & 1.139 & 0.002 & 0.780 & 0.114 & 0.414 & 0.021 & 0.432 & 0.133 \\
    & & & 60 & - & - & \textbf{0.234} & 0.028 & 0.276 & 0.009 & 0.831 & 0.002 & 0.611 & 0.116 & 0.291 & 0.020 & 0.361 & 0.149 \\
    & & & 100 & 0.257 & 0.017 & \textbf{0.087} & 0.035 & 0.105 & 0.010 & 0.259 & 0.005 & 0.198 & 0.127 & 0.102 & 0.021 & 0.131 & 0.281 \\
    & & & 200 & 0.035 & 0.017 & \textbf{0.019} & 0.028 & 0.022 & 0.011 & 0.046 & 0.007 & 0.060 & 0.110 & 0.021 & 0.022 & 0.021 & 0.489 \\
    \cmidrule(lr){2-18}

    & \multirow{21}{*}{$3{,}000$} & \multirow{7}{*}{$0.1$} & 20 & - & - & \textbf{6.036} & 0.165 & 6.112 & 0.127 & 32.352 & 0.002 & 12.359 & 0.332 & 8.134 & 0.058 & 14.380 & 0.087 \\
    & & & 30 & - & - & \textbf{2.527} & 0.169 & 2.784 & 0.137 & 20.088 & 0.003 & 5.145 & 0.333 & 3.051 & 0.066 & 8.229 & 0.127 \\
    & & & 40 & - & - & \textbf{1.104} & 0.172 & 1.244 & 0.146 & 16.128 & 0.004 & 2.202 & 0.388 & 1.158 & 0.061 & 5.314 & 0.173 \\
    & & & 50 & - & - & \textbf{0.750} & 0.209 & 0.821 & 0.183 & 9.909 & 0.004 & 1.289 & 0.389 & 0.836 & 0.070 & 3.366 & 0.221 \\
    & & & 60 & - & - & \textbf{0.463} & 0.233 & 0.548 & 0.192 & 7.299 & 0.005 & 1.144 & 0.371 & 0.657 & 0.067 & 2.248 & 0.294 \\
    & & & 100 & 0.456 & 0.196 & \textbf{0.125} & 0.214 & 0.159 & 0.182 & 2.414 & 0.008 & 0.269 & 0.395 & 0.157 & 0.063 & 0.639 & 0.620 \\
    & & & 200 & 0.081 & 0.196 & \textbf{0.033} & 0.221 & 0.038 & 0.201 & 0.547 & 0.016 & 0.082 & 0.390 & 0.045 & 0.067 & 0.122 & 1.434 \\
    \cmidrule(lr){3-18}

    & & \multirow{7}{*}{$0.5$} & 20 & - & - & \textbf{2.499} & 0.077 & 2.856 & 0.045 & 7.906 & 0.002 & 4.039 & 0.323 & 2.942 & 0.055 & 3.379 & 0.087 \\
    & & & 30 & - & - & \textbf{1.362} & 0.085 & 1.436 & 0.058 & 5.018 & 0.003 & 2.517 & 0.296 & 1.575 & 0.058 & 2.011 & 0.127 \\
    & & & 40 & - & - & \textbf{0.663} & 0.085 & 0.734 & 0.058 & 3.756 & 0.004 & 1.664 & 0.366 & 0.802 & 0.055 & 1.310 & 0.174 \\
    & & & 50 & - & - & \textbf{0.517} & 0.105 & 0.571 & 0.076 & 2.476 & 0.004 & 0.863 & 0.377 & 0.639 & 0.060 & 0.740 & 0.217 \\
    & & & 60 & - & - & \textbf{0.380} & 0.108 & 0.417 & 0.076 & 1.954 & 0.005 & 0.775 & 0.364 & 0.480 & 0.059 & 0.648 & 0.285 \\
    & & & 100 & 0.356 & 0.087 & \textbf{0.106} & 0.106 & 0.141 & 0.084 & 0.632 & 0.008 & 0.266 & 0.390 & 0.146 & 0.059 & 0.232 & 0.622 \\
    & & & 200 & 0.080 & 0.100 & \textbf{0.030} & 0.116 & 0.035 & 0.088 & 0.098 & 0.016 & 0.083 & 0.372 & 0.042 & 0.055 & 0.039 & 1.371 \\
    \cmidrule(lr){3-18}

    & & \multirow{7}{*}{$1.0$} & 20 & - & - & \textbf{1.427} & 0.058 & 1.653 & 0.031 & 3.322 & 0.002 & 2.583 & 0.322 & 2.047 & 0.053 & 1.922 & 0.091 \\
    & & & 30 & - & - & \textbf{0.867} & 0.072 & 0.973 & 0.040 & 2.433 & 0.003 & 1.614 & 0.328 & 1.139 & 0.058 & 1.150 & 0.126 \\
    & & & 40 & - & - & \textbf{0.461} & 0.066 & 0.495 & 0.043 & 1.650 & 0.004 & 1.077 & 0.344 & 0.609 & 0.054 & 0.751 & 0.177 \\
    & & & 50 & - & - & \textbf{0.354} & 0.081 & 0.387 & 0.055 & 1.468 & 0.004 & 0.603 & 0.359 & 0.429 & 0.061 & 0.540 & 0.220 \\
    & & & 60 & - & - & \textbf{0.273} & 0.076 & 0.329 & 0.051 & 1.218 & 0.005 & 0.677 & 0.359 & 0.317 & 0.056 & 0.359 & 0.300 \\
    & & & 100 & 0.274 & 0.067 & \textbf{0.099} & 0.071 & 0.118 & 0.049 & 0.417 & 0.009 & 0.294 & 0.386 & 0.111 & 0.058 & 0.169 & 0.618 \\
    & & & 200 & 0.065 & 0.072 & \textbf{0.024} & 0.088 & 0.030 & 0.063 & 0.069 & 0.016 & 0.075 & 0.340 & 0.037 & 0.055 & 0.035 & 1.489 \\
    \cmidrule(lr){2-18}

    & \multirow{21}{*}{$5{,}000$} & \multirow{7}{*}{$0.1$} & 20 & - & - & \textbf{7.847} & 0.336 & 8.971 & 0.309 & 33.989 & 0.003 & 13.767 & 0.550 & 9.599 & 0.105 & 16.084 & 0.123 \\
    & & & 30 & - & - & \textbf{1.886} & 0.339 & 2.598 & 0.305 & 22.410 & 0.005 & 5.096 & 0.578 & 2.615 & 0.108 & 9.647 & 0.184 \\
    & & & 40 & - & - & 1.432 & 0.354 & 1.495 & 0.316 & 14.230 & 0.006 & 2.583 & 0.553 & \textbf{1.361} & 0.119 & 5.719 & 0.233 \\
    & & & 50 & - & - & \textbf{0.691} & 0.331 & 0.936 & 0.301 & 10.702 & 0.008 & 1.544 & 0.636 & 0.733 & 0.110 & 3.972 & 0.343 \\
    & & & 60 & - & - & \textbf{0.399} & 0.372 & 0.535 & 0.341 & 7.610 & 0.009 & 1.104 & 0.625 & 0.522 & 0.118 & 2.803 & 0.494 \\
    & & & 100 & 0.461 & 0.406 & \textbf{0.143} & 0.419 & 0.171 & 0.382 & 2.673 & 0.014 & 0.399 & 0.650 & 0.185 & 0.111 & 0.907 & 0.948 \\
    & & & 200 & 0.101 & 0.417 & \textbf{0.030} & 0.441 & 0.036 & 0.412 & 0.635 & 0.032 & 0.095 & 0.733 & 0.040 & 0.123 & 0.184 & 2.177 \\
    \cmidrule(lr){3-18}

    & & \multirow{7}{*}{$0.5$} & 20 & - & - & \textbf{2.623} & 0.161 & 2.739 & 0.129 & 6.735 & 0.004 & 5.536 & 0.524 & 3.350 & 0.098 & 2.993 & 0.116 \\
    & & & 30 & - & - & \textbf{1.129} & 0.222 & 1.326 & 0.129 & 5.115 & 0.007 & 3.804 & 0.726 & 1.442 & 0.106 & 1.871 & 0.277 \\
    & & & 40 & - & - & 0.862 & 0.199 & 0.908 & 0.133 & 4.005 & 0.007 & 1.507 & 0.684 & \textbf{0.852} & 0.107 & 1.205 & 0.387 \\
    & & & 50 & - & - & \textbf{0.445} & 0.170 & 0.544 & 0.142 & 2.828 & 0.008 & 1.151 & 0.604 & 0.553 & 0.104 & 0.715 & 0.380 \\
    & & & 60 & - & - & \textbf{0.324} & 0.182 & 0.380 & 0.151 & 2.212 & 0.009 & 1.050 & 0.626 & 0.451 & 0.110 & 0.578 & 0.472 \\
    & & & 100 & 0.412 & 0.213 & \textbf{0.130} & 0.200 & 0.146 & 0.171 & 0.794 & 0.015 & 0.361 & 0.668 & 0.146 & 0.105 & 0.243 & 0.963 \\
    & & & 200 & 0.107 & 0.204 & \textbf{0.027} & 0.215 & 0.034 & 0.190 & 0.123 & 0.032 & 0.084 & 0.692 & 0.035 & 0.122 & 0.044 & 2.214 \\
    \cmidrule(lr){3-18}

    & & \multirow{7}{*}{$1.0$} & 20 & - & - & \textbf{1.440} & 0.121 & 1.655 & 0.096 & 3.676 & 0.004 & 3.578 & 0.544 & 2.484 & 0.100 & 1.696 & 0.123 \\
    & & & 30 & - & - & \textbf{0.753} & 0.134 & 0.878 & 0.087 & 2.265 & 0.005 & 1.864 & 0.606 & 1.112 & 0.097 & 1.045 & 0.199 \\
    & & & 40 & - & - & \textbf{0.539} & 0.180 & 0.568 & 0.097 & 1.832 & 0.008 & 0.863 & 0.772 & 0.692 & 0.104 & 0.921 & 0.441 \\
    & & & 50 & - & - & \textbf{0.325} & 0.126 & 0.384 & 0.103 & 1.511 & 0.008 & 0.959 & 0.580 & 0.387 & 0.101 & 0.529 & 0.384 \\
    & & & 60 & - & - & \textbf{0.248} & 0.140 & 0.304 & 0.116 & 0.991 & 0.009 & 0.703 & 0.612 & 0.310 & 0.106 & 0.414 & 0.485 \\
    & & & 100 & 0.341 & 0.158 & \textbf{0.113} & 0.159 & 0.128 & 0.140 & 0.538 & 0.014 & 0.264 & 0.622 & 0.139 & 0.100 & 0.175 & 0.918 \\
    & & & 200 & 0.099 & 0.153 & \textbf{0.024} & 0.167 & 0.028 & 0.144 & 0.085 & 0.028 & 0.088 & 0.629 & 0.033 & 0.119 & 0.037 & 2.030 \\

    \bottomrule
\end{tabular}
}
\endgroup
\end{table}

\begin{table}[!htbp]
\caption{Average gaps and average times on synthetic data with $r=40$.}
\label{tab:synthetic_r40}
\centering
\begingroup
\scriptsize
\renewcommand{\arraystretch}{1.15}
\setlength{\tabcolsep}{2.2pt}
\resizebox{1.0\textwidth}{!}{%
\begin{tabular}{cccc*{7}{cc}}
    \toprule
    \multirow{2}{*}{$r$}
        & \multirow{2}{*}{$p$}
        & \multirow{2}{*}{$\lambda$}
        & \multirow{2}{*}{$k$}
        & \multicolumn{2}{c}{\Cref{alg:reduce}}
        & \multicolumn{2}{c}{\Cref{alg:round}}
        & \multicolumn{2}{c}{\Cref{alg:pipage}}
        & \multicolumn{2}{c}{Greedy}
        & \multicolumn{2}{c}{Randomized}
        & \multicolumn{2}{c}{DIHT}
        & \multicolumn{2}{c}{Beam Search} \\
    \cmidrule(lr){5-6}
    \cmidrule(lr){7-8}
    \cmidrule(lr){9-10}
    \cmidrule(lr){11-12}
    \cmidrule(lr){13-14}
    \cmidrule(lr){15-16}
    \cmidrule(lr){17-18}
    & & & &
    \makecell{Avg.\\Gap} & \makecell{Avg.\\Time} &
    \makecell{Avg.\\Gap} & \makecell{Avg.\\Time} &
    \makecell{Avg.\\Gap} & \makecell{Avg.\\Time} &
    \makecell{Avg.\\Gap} & \makecell{Avg.\\Time} &
    \makecell{Avg.\\Gap} & \makecell{Avg.\\Time} &
    \makecell{Avg.\\Gap} & \makecell{Avg.\\Time} &
    \makecell{Avg.\\Gap} & \makecell{Avg.\\Time} \\
    \midrule

    \multirow{63}{*}{$40$} & \multirow{21}{*}{$1{,}000$} & \multirow{7}{*}{$0.1$} & 20 & - & - & \textbf{28.251} & 0.249 & 31.242 & 0.022 & 54.221 & 0.001 & 64.410 & 0.127 & 46.618 & 0.025 & 31.983 & 0.058 \\
    & & & 30 & - & - & \textbf{12.656} & 0.283 & 14.248 & 0.024 & 44.099 & 0.002 & 29.113 & 0.137 & 22.323 & 0.029 & 19.742 & 0.087 \\
    & & & 40 & - & - & \textbf{5.466} & 0.197 & 6.983 & 0.026 & 31.808 & 0.002 & 15.452 & 0.133 & 9.711 & 0.029 & 14.766 & 0.120 \\
    & & & 50 & - & - & \textbf{3.771} & 0.255 & 4.030 & 0.024 & 23.284 & 0.003 & 7.660 & 0.132 & 6.353 & 0.029 & 10.234 & 0.146 \\
    & & & 60 & - & - & 2.406 & 0.188 & \textbf{2.211} & 0.026 & 17.881 & 0.003 & 4.693 & 0.139 & 3.092 & 0.028 & 7.688 & 0.176 \\
    & & & 100 & 1.538 & 0.052 & \textbf{0.507} & 0.141 & 0.650 & 0.026 & 5.958 & 0.005 & 1.304 & 0.136 & 0.844 & 0.028 & 2.742 & 0.307 \\
    & & & 200 & 0.227 & 0.053 & \textbf{0.088} & 0.111 & 0.106 & 0.028 & 1.024 & 0.010 & 0.256 & 0.137 & 0.144 & 0.029 & 0.336 & 0.625 \\
    \cmidrule(lr){3-18}

    & & \multirow{7}{*}{$0.5$} & 20 & - & - & \textbf{6.109} & 0.125 & 6.763 & 0.010 & 11.750 & 0.001 & 13.231 & 0.145 & 7.952 & 0.024 & 6.450 & 0.057 \\
    & & & 30 & - & - & \textbf{3.315} & 0.119 & 3.830 & 0.011 & 8.731 & 0.002 & 7.842 & 0.139 & 7.045 & 0.027 & 4.481 & 0.085 \\
    & & & 40 & - & - & \textbf{1.983} & 0.125 & 2.554 & 0.011 & 6.903 & 0.002 & 5.607 & 0.140 & 3.930 & 0.028 & 3.486 & 0.115 \\
    & & & 50 & - & - & \textbf{1.463} & 0.103 & 1.684 & 0.011 & 5.888 & 0.003 & 3.935 & 0.142 & 2.413 & 0.027 & 2.328 & 0.144 \\
    & & & 60 & - & - & \textbf{0.977} & 0.110 & 1.187 & 0.013 & 4.205 & 0.004 & 2.611 & 0.144 & 1.821 & 0.027 & 2.185 & 0.175 \\
    & & & 100 & 1.199 & 0.041 & \textbf{0.348} & 0.080 & 0.390 & 0.011 & 1.730 & 0.006 & 0.892 & 0.139 & 0.549 & 0.028 & 0.858 & 0.305 \\
    & & & 200 & 0.211 & 0.040 & \textbf{0.076} & 0.087 & 0.080 & 0.012 & 0.278 & 0.014 & 0.207 & 0.138 & 0.121 & 0.029 & 0.154 & 0.624 \\
    \cmidrule(lr){3-18}

    & & \multirow{7}{*}{$1.0$} & 20 & - & - & \textbf{2.489} & 0.047 & 2.944 & 0.007 & 4.156 & 0.001 & 6.287 & 0.140 & 6.585 & 0.027 & 2.883 & 0.056 \\
    & & & 30 & - & - & \textbf{1.831} & 0.075 & 2.001 & 0.007 & 3.545 & 0.002 & 4.439 & 0.142 & 5.286 & 0.027 & 2.272 & 0.084 \\
    & & & 40 & - & - & \textbf{1.139} & 0.079 & 1.344 & 0.007 & 2.932 & 0.003 & 2.801 & 0.154 & 3.200 & 0.027 & 1.572 & 0.119 \\
    & & & 50 & - & - & \textbf{0.818} & 0.077 & 0.982 & 0.009 & 2.835 & 0.003 & 1.907 & 0.140 & 2.462 & 0.026 & 1.247 & 0.144 \\
    & & & 60 & - & - & \textbf{0.645} & 0.091 & 0.783 & 0.009 & 1.989 & 0.004 & 1.687 & 0.146 & 2.013 & 0.026 & 1.017 & 0.169 \\
    & & & 100 & 0.786 & 0.033 & \textbf{0.267} & 0.085 & 0.312 & 0.009 & 0.825 & 0.006 & 0.663 & 0.155 & 0.622 & 0.026 & 0.475 & 0.300 \\
    & & & 200 & 0.170 & 0.032 & \textbf{0.068} & 0.060 & 0.077 & 0.009 & 0.192 & 0.010 & 0.190 & 0.149 & 0.096 & 0.027 & 0.117 & 0.628 \\
    \cmidrule(lr){2-18}

    & \multirow{21}{*}{$3{,}000$} & \multirow{7}{*}{$0.1$} & 20 & - & - & \textbf{27.769} & 0.361 & 30.172 & 0.122 & 64.592 & 0.004 & 73.801 & 0.415 & 43.669 & 0.071 & 35.576 & 0.096 \\
    & & & 30 & - & - & \textbf{13.443} & 0.388 & 14.754 & 0.102 & 46.461 & 0.005 & 31.824 & 0.440 & 22.157 & 0.065 & 25.302 & 0.148 \\
    & & & 40 & - & - & \textbf{6.547} & 0.387 & 7.485 & 0.147 & 33.759 & 0.007 & 17.972 & 0.460 & 12.686 & 0.074 & 17.532 & 0.197 \\
    & & & 50 & - & - & \textbf{3.611} & 0.316 & 4.379 & 0.121 & 27.653 & 0.008 & 8.183 & 0.458 & 6.048 & 0.070 & 13.013 & 0.247 \\
    & & & 60 & - & - & \textbf{1.886} & 0.277 & 2.089 & 0.141 & 22.325 & 0.009 & 4.755 & 0.464 & 3.601 & 0.073 & 9.362 & 0.292 \\
    & & & 100 & 2.290 & 0.153 & \textbf{0.589} & 0.387 & 0.803 & 0.141 & 9.455 & 0.016 & 1.870 & 0.480 & 1.092 & 0.073 & 3.435 & 0.506 \\
    & & & 200 & 0.327 & 0.174 & \textbf{0.112} & 0.269 & 0.131 & 0.135 & 2.159 & 0.035 & 0.316 & 0.489 & 0.183 & 0.079 & 0.803 & 1.027 \\
    \cmidrule(lr){3-18}

    & & \multirow{7}{*}{$0.5$} & 20 & - & - & \textbf{6.460} & 0.165 & 7.683 & 0.046 & 13.924 & 0.004 & 16.228 & 0.403 & 9.205 & 0.065 & 7.880 & 0.099 \\
    & & & 30 & - & - & \textbf{3.755} & 0.254 & 4.090 & 0.052 & 10.375 & 0.005 & 10.171 & 0.432 & 5.672 & 0.061 & 5.151 & 0.149 \\
    & & & 40 & - & - & \textbf{2.251} & 0.227 & 2.604 & 0.055 & 8.615 & 0.007 & 6.952 & 0.423 & 3.559 & 0.067 & 3.854 & 0.197 \\
    & & & 50 & - & - & \textbf{1.550} & 0.193 & 1.755 & 0.054 & 6.381 & 0.008 & 4.108 & 0.430 & 2.173 & 0.064 & 2.993 & 0.230 \\
    & & & 60 & - & - & \textbf{1.046} & 0.194 & 1.093 & 0.051 & 5.203 & 0.010 & 2.668 & 0.442 & 1.460 & 0.066 & 2.318 & 0.288 \\
    & & & 100 & 1.849 & 0.087 & 0.495 & 0.272 & \textbf{0.482} & 0.061 & 2.153 & 0.016 & 1.439 & 0.486 & 0.657 & 0.069 & 1.123 & 0.499 \\
    & & & 200 & 0.297 & 0.103 & \textbf{0.099} & 0.168 & 0.117 & 0.061 & 0.419 & 0.031 & 0.300 & 0.477 & 0.177 & 0.073 & 0.241 & 1.051 \\
    \cmidrule(lr){3-18}

    & & \multirow{7}{*}{$1.0$} & 20 & - & - & \textbf{2.891} & 0.078 & 3.131 & 0.027 & 5.794 & 0.004 & 7.196 & 0.403 & 6.265 & 0.063 & 3.510 & 0.094 \\
    & & & 30 & - & - & \textbf{1.940} & 0.162 & 2.174 & 0.036 & 4.871 & 0.005 & 5.085 & 0.436 & 4.097 & 0.060 & 2.443 & 0.147 \\
    & & & 40 & - & - & \textbf{1.388} & 0.183 & 1.667 & 0.041 & 4.077 & 0.007 & 3.865 & 0.433 & 2.742 & 0.066 & 1.908 & 0.195 \\
    & & & 50 & - & - & \textbf{0.950} & 0.103 & 0.994 & 0.032 & 3.207 & 0.008 & 2.617 & 0.462 & 2.023 & 0.065 & 1.498 & 0.226 \\
    & & & 60 & - & - & \textbf{0.664} & 0.122 & 0.753 & 0.035 & 2.078 & 0.010 & 1.688 & 0.474 & 1.135 & 0.068 & 1.203 & 0.297 \\
    & & & 100 & 1.186 & 0.067 & \textbf{0.328} & 0.177 & 0.387 & 0.037 & 1.216 & 0.016 & 1.039 & 0.466 & 0.564 & 0.065 & 0.754 & 0.488 \\
    & & & 200 & 0.248 & 0.081 & \textbf{0.086} & 0.200 & 0.109 & 0.048 & 0.295 & 0.031 & 0.215 & 0.486 & 0.138 & 0.073 & 0.165 & 1.023 \\
    \cmidrule(lr){2-18}

    & \multirow{21}{*}{$5{,}000$} & \multirow{7}{*}{$0.1$} & 20 & - & - & \textbf{27.585} & 0.499 & 29.594 & 0.248 & 60.618 & 0.006 & 71.396 & 0.693 & 52.968 & 0.134 & 36.573 & 0.127 \\
    & & & 30 & - & - & \textbf{12.625} & 0.572 & 13.846 & 0.283 & 50.369 & 0.009 & 29.492 & 0.714 & 21.703 & 0.146 & 26.215 & 0.180 \\
    & & & 40 & - & - & \textbf{6.410} & 0.571 & 8.028 & 0.327 & 37.392 & 0.011 & 15.824 & 0.761 & 10.336 & 0.150 & 18.649 & 0.286 \\
    & & & 50 & - & - & \textbf{3.715} & 0.468 & 4.864 & 0.271 & 28.334 & 0.014 & 8.559 & 0.761 & 6.119 & 0.142 & 14.423 & 0.358 \\
    & & & 60 & - & - & \textbf{2.353} & 0.415 & 2.665 & 0.265 & 24.464 & 0.016 & 6.308 & 0.778 & 3.704 & 0.144 & 10.147 & 0.428 \\
    & & & 100 & 2.461 & 0.362 & \textbf{0.597} & 0.445 & 0.755 & 0.300 & 10.763 & 0.027 & 1.731 & 0.898 & 1.024 & 0.125 & 4.146 & 0.700 \\
    & & & 200 & 0.373 & 0.393 & \textbf{0.124} & 0.552 & 0.149 & 0.356 & 2.673 & 0.054 & 0.433 & 0.891 & 0.221 & 0.133 & 0.831 & 1.374 \\
    \cmidrule(lr){3-18}

    & & \multirow{7}{*}{$0.5$} & 20 & - & - & \textbf{6.188} & 0.263 & 7.114 & 0.106 & 13.522 & 0.006 & 16.814 & 0.749 & 10.149 & 0.128 & 8.760 & 0.129 \\
    & & & 30 & - & - & \textbf{3.638} & 0.276 & 4.090 & 0.117 & 11.215 & 0.009 & 9.583 & 0.724 & 6.312 & 0.130 & 6.051 & 0.184 \\
    & & & 40 & - & - & \textbf{2.304} & 0.287 & 3.099 & 0.124 & 9.326 & 0.012 & 6.192 & 0.784 & 3.610 & 0.139 & 4.415 & 0.295 \\
    & & & 50 & - & - & \textbf{1.597} & 0.273 & 1.995 & 0.128 & 7.552 & 0.014 & 3.962 & 0.754 & 2.525 & 0.130 & 3.252 & 0.363 \\
    & & & 60 & - & - & \textbf{1.138} & 0.219 & 1.418 & 0.113 & 5.990 & 0.016 & 2.979 & 0.734 & 1.837 & 0.134 & 2.709 & 0.419 \\
    & & & 100 & 1.571 & 0.194 & \textbf{0.402} & 0.308 & 0.532 & 0.151 & 2.805 & 0.026 & 1.363 & 0.825 & 0.748 & 0.129 & 1.030 & 0.695 \\
    & & & 200 & 0.347 & 0.217 & \textbf{0.103} & 0.362 & 0.124 & 0.165 & 0.645 & 0.051 & 0.362 & 0.885 & 0.205 & 0.138 & 0.221 & 1.395 \\
    \cmidrule(lr){3-18}

    & & \multirow{7}{*}{$1.0$} & 20 & - & - & \textbf{3.080} & 0.164 & 3.440 & 0.066 & 6.259 & 0.006 & 8.455 & 0.725 & 6.232 & 0.125 & 3.449 & 0.122 \\
    & & & 30 & - & - & \textbf{2.034} & 0.178 & 2.195 & 0.074 & 5.628 & 0.009 & 5.550 & 0.696 & 4.092 & 0.129 & 2.606 & 0.187 \\
    & & & 40 & - & - & \textbf{1.416} & 0.259 & 1.664 & 0.105 & 3.587 & 0.011 & 3.694 & 0.729 & 2.850 & 0.141 & 1.966 & 0.288 \\
    & & & 50 & - & - & \textbf{1.047} & 0.189 & 1.136 & 0.085 & 3.076 & 0.014 & 2.616 & 0.813 & 1.814 & 0.137 & 1.715 & 0.386 \\
    & & & 60 & - & - & 0.812 & 0.176 & \textbf{0.809} & 0.078 & 2.744 & 0.019 & 2.359 & 0.795 & 1.471 & 0.141 & 1.351 & 0.417 \\
    & & & 100 & 1.022 & 0.159 & \textbf{0.313} & 0.239 & 0.371 & 0.093 & 1.519 & 0.026 & 1.089 & 0.804 & 0.566 & 0.120 & 0.658 & 0.692 \\
    & & & 200 & 0.286 & 0.174 & \textbf{0.094} & 0.261 & 0.110 & 0.106 & 0.341 & 0.051 & 0.278 & 0.886 & 0.156 & 0.133 & 0.158 & 1.385 \\

    \bottomrule
\end{tabular}
}
\endgroup
\end{table}

\begin{table}[!htbp]
\caption{Average gaps and average times on synthetic data with $r=60$.}
\label{tab:synthetic_r60}
\centering
\begingroup
\scriptsize
\renewcommand{\arraystretch}{1.15}
\setlength{\tabcolsep}{2.2pt}
\resizebox{1.0\textwidth}{!}{%
\begin{tabular}{cccc*{7}{cc}}
    \toprule
    \multirow{2}{*}{$r$}
        & \multirow{2}{*}{$p$}
        & \multirow{2}{*}{$\lambda$}
        & \multirow{2}{*}{$k$}
        & \multicolumn{2}{c}{\Cref{alg:reduce}}
        & \multicolumn{2}{c}{\Cref{alg:round}}
        & \multicolumn{2}{c}{\Cref{alg:pipage}}
        & \multicolumn{2}{c}{Greedy}
        & \multicolumn{2}{c}{Randomized}
        & \multicolumn{2}{c}{DIHT}
        & \multicolumn{2}{c}{Beam Search} \\
    \cmidrule(lr){5-6}
    \cmidrule(lr){7-8}
    \cmidrule(lr){9-10}
    \cmidrule(lr){11-12}
    \cmidrule(lr){13-14}
    \cmidrule(lr){15-16}
    \cmidrule(lr){17-18}
    & & & &
    \makecell{Avg.\\Gap} & \makecell{Avg.\\Time} &
    \makecell{Avg.\\Gap} & \makecell{Avg.\\Time} &
    \makecell{Avg.\\Gap} & \makecell{Avg.\\Time} &
    \makecell{Avg.\\Gap} & \makecell{Avg.\\Time} &
    \makecell{Avg.\\Gap} & \makecell{Avg.\\Time} &
    \makecell{Avg.\\Gap} & \makecell{Avg.\\Time} &
    \makecell{Avg.\\Gap} & \makecell{Avg.\\Time} \\
    \midrule

    \multirow{63}{*}{$60$} & \multirow{21}{*}{$1{,}000$} & \multirow{7}{*}{$0.1$} & 20 & - & - & \textbf{47.095} & 1.843 & 52.120 & 0.017 & 69.529 & 0.002 & 130.949 & 0.167 & 86.689 & 0.031 & 47.393 & 0.058 \\
    & & & 30 & - & - & \textbf{27.880} & 4.292 & 30.984 & 0.020 & 51.137 & 0.003 & 75.290 & 0.169 & 54.805 & 0.031 & 33.681 & 0.095 \\
    & & & 40 & - & - & \textbf{15.204} & 5.431 & 19.632 & 0.021 & 45.616 & 0.004 & 44.510 & 0.190 & 29.024 & 0.030 & 27.969 & 0.229 \\
    & & & 50 & - & - & \textbf{8.654} & 3.455 & 9.879 & 0.023 & 33.994 & 0.004 & 20.883 & 0.160 & 16.499 & 0.031 & 21.866 & 0.135 \\
    & & & 60 & - & - & \textbf{6.257} & 3.031 & 7.237 & 0.022 & 27.294 & 0.005 & 16.141 & 0.156 & 11.220 & 0.032 & 15.504 & 0.192 \\
    & & & 100 & 5.458 & 0.067 & \textbf{1.487} & 1.600 & 1.589 & 0.023 & 12.248 & 0.007 & 3.882 & 0.170 & 2.420 & 0.032 & 5.733 & 0.336 \\
    & & & 200 & 0.700 & 0.077 & \textbf{0.195} & 1.302 & 0.215 & 0.025 & 2.593 & 0.014 & 0.627 & 0.176 & 0.362 & 0.032 & 1.058 & 0.793 \\
    \cmidrule(lr){3-18}

    & & \multirow{7}{*}{$0.5$} & 20 & - & - & \textbf{8.437} & 0.494 & 9.371 & 0.008 & 12.168 & 0.002 & 20.196 & 0.155 & 13.221 & 0.029 & 8.789 & 0.059 \\
    & & & 30 & - & - & \textbf{5.603} & 1.250 & 6.184 & 0.010 & 10.813 & 0.003 & 14.566 & 0.169 & 11.069 & 0.030 & 7.156 & 0.095 \\
    & & & 40 & - & - & \textbf{3.800} & 1.350 & 4.204 & 0.010 & 7.947 & 0.004 & 11.044 & 0.189 & 8.954 & 0.029 & 5.616 & 0.216 \\
    & & & 50 & - & - & \textbf{2.722} & 1.651 & 2.981 & 0.011 & 7.594 & 0.004 & 8.602 & 0.160 & 5.771 & 0.031 & 4.839 & 0.137 \\
    & & & 60 & - & - & \textbf{1.922} & 1.488 & 2.438 & 0.010 & 5.299 & 0.005 & 5.565 & 0.159 & 4.709 & 0.031 & 3.440 & 0.196 \\
    & & & 100 & 2.495 & 0.063 & \textbf{0.717} & 0.930 & 0.816 & 0.011 & 2.942 & 0.011 & 2.532 & 0.170 & 1.677 & 0.030 & 1.908 & 0.341 \\
    & & & 200 & 0.467 & 0.065 & \textbf{0.151} & 0.478 & 0.180 & 0.012 & 0.583 & 0.019 & 0.549 & 0.173 & 0.306 & 0.031 & 0.334 & 0.782 \\
    \cmidrule(lr){3-18}

    & & \multirow{7}{*}{$1.0$} & 20 & - & - & \textbf{3.529} & 0.146 & 3.826 & 0.006 & 4.492 & 0.002 & 7.892 & 0.164 & 8.205 & 0.029 & 3.769 & 0.058 \\
    & & & 30 & - & - & \textbf{2.560} & 0.473 & 2.823 & 0.007 & 4.045 & 0.003 & 6.747 & 0.175 & 6.259 & 0.029 & 3.148 & 0.125 \\
    & & & 40 & - & - & \textbf{1.917} & 0.544 & 2.132 & 0.007 & 3.544 & 0.003 & 5.332 & 0.168 & 6.003 & 0.029 & 2.312 & 0.106 \\
    & & & 50 & - & - & \textbf{1.430} & 0.604 & 1.541 & 0.007 & 2.947 & 0.004 & 3.777 & 0.162 & 5.094 & 0.030 & 2.014 & 0.136 \\
    & & & 60 & - & - & \textbf{1.089} & 0.562 & 1.258 & 0.008 & 2.366 & 0.005 & 3.528 & 0.157 & 3.647 & 0.030 & 1.680 & 0.196 \\
    & & & 100 & 1.632 & 0.055 & \textbf{0.488} & 0.662 & 0.560 & 0.008 & 1.485 & 0.008 & 1.431 & 0.177 & 1.548 & 0.030 & 1.039 & 0.341 \\
    & & & 200 & 0.328 & 0.058 & \textbf{0.128} & 0.487 & 0.137 & 0.009 & 0.420 & 0.018 & 0.350 & 0.167 & 0.396 & 0.031 & 0.211 & 0.773 \\
    \cmidrule(lr){2-18}

    & \multirow{21}{*}{$3{,}000$} & \multirow{7}{*}{$0.1$} & 20 & - & - & 52.014 & 4.893 & 59.448 & 0.097 & 72.399 & 0.005 & 145.513 & 0.507 & 118.159 & 0.082 & \textbf{50.621} & 0.098 \\
    & & & 30 & - & - & \textbf{26.534} & 8.008 & 32.524 & 0.121 & 63.867 & 0.008 & 78.280 & 0.492 & 62.914 & 0.088 & 34.947 & 0.142 \\
    & & & 40 & - & - & \textbf{16.238} & 8.536 & 20.022 & 0.124 & 52.924 & 0.010 & 42.669 & 0.511 & 33.484 & 0.087 & 30.508 & 0.200 \\
    & & & 50 & - & - & \textbf{10.587} & 6.585 & 11.246 & 0.117 & 41.426 & 0.012 & 29.860 & 0.533 & 21.827 & 0.086 & 24.451 & 0.259 \\
    & & & 60 & - & - & \textbf{6.719} & 5.807 & 7.833 & 0.131 & 38.032 & 0.014 & 17.812 & 0.560 & 13.101 & 0.087 & 18.767 & 0.311 \\
    & & & 100 & 7.063 & 0.216 & \textbf{1.560} & 3.880 & 2.071 & 0.135 & 17.746 & 0.023 & 4.645 & 0.563 & 3.867 & 0.082 & 8.430 & 0.542 \\
    & & & 200 & 0.860 & 0.230 & \textbf{0.241} & 2.613 & 0.322 & 0.130 & 4.509 & 0.045 & 0.995 & 0.592 & 0.555 & 0.088 & 2.021 & 1.218 \\
    \cmidrule(lr){3-18}

    & & \multirow{7}{*}{$0.5$} & 20 & - & - & \textbf{9.403} & 1.352 & 11.698 & 0.041 & 14.035 & 0.005 & 26.610 & 0.536 & 19.223 & 0.072 & 10.628 & 0.100 \\
    & & & 30 & - & - & \textbf{6.110} & 2.098 & 7.549 & 0.043 & 12.448 & 0.008 & 18.960 & 0.514 & 12.374 & 0.081 & 7.984 & 0.147 \\
    & & & 40 & - & - & \textbf{4.433} & 2.851 & 4.840 & 0.046 & 11.175 & 0.010 & 12.290 & 0.551 & 7.493 & 0.080 & 7.117 & 0.207 \\
    & & & 50 & - & - & \textbf{2.959} & 2.472 & 3.612 & 0.048 & 9.443 & 0.012 & 9.125 & 0.528 & 5.562 & 0.082 & 6.123 & 0.259 \\
    & & & 60 & - & - & \textbf{2.348} & 2.657 & 2.673 & 0.048 & 7.907 & 0.014 & 7.494 & 0.535 & 4.444 & 0.077 & 4.305 & 0.309 \\
    & & & 100 & 3.455 & 0.143 & \textbf{0.863} & 2.409 & 1.183 & 0.052 & 4.357 & 0.023 & 2.989 & 0.568 & 1.823 & 0.080 & 2.419 & 0.548 \\
    & & & 200 & 0.596 & 0.169 & \textbf{0.198} & 2.463 & 0.238 & 0.054 & 1.135 & 0.045 & 0.881 & 0.603 & 0.500 & 0.078 & 0.618 & 1.193 \\
    \cmidrule(lr){3-18}

    & & \multirow{7}{*}{$1.0$} & 20 & - & - & \textbf{4.359} & 0.477 & 5.156 & 0.024 & 6.103 & 0.005 & 10.828 & 0.552 & 8.150 & 0.073 & 4.579 & 0.099 \\
    & & & 30 & - & - & \textbf{2.927} & 0.852 & 3.300 & 0.029 & 5.389 & 0.008 & 9.182 & 0.516 & 6.578 & 0.072 & 3.713 & 0.149 \\
    & & & 40 & - & - & \textbf{2.019} & 0.927 & 2.223 & 0.033 & 4.501 & 0.010 & 5.937 & 0.555 & 5.838 & 0.079 & 2.867 & 0.208 \\
    & & & 50 & - & - & \textbf{1.606} & 1.138 & 1.864 & 0.034 & 3.897 & 0.012 & 4.614 & 0.556 & 4.287 & 0.076 & 2.355 & 0.256 \\
    & & & 60 & - & - & \textbf{1.261} & 1.282 & 1.535 & 0.038 & 3.579 & 0.014 & 3.821 & 0.561 & 3.098 & 0.073 & 2.205 & 0.307 \\
    & & & 100 & 2.172 & 0.124 & \textbf{0.620} & 1.626 & 0.733 & 0.036 & 2.157 & 0.023 & 1.950 & 0.591 & 1.453 & 0.075 & 1.241 & 0.561 \\
    & & & 200 & 0.526 & 0.150 & \textbf{0.167} & 1.820 & 0.196 & 0.041 & 0.636 & 0.045 & 0.588 & 0.527 & 0.342 & 0.082 & 0.373 & 1.064 \\
    \cmidrule(lr){2-18}

    & \multirow{21}{*}{$5{,}000$} & \multirow{7}{*}{$0.1$} & 20 & - & - & \textbf{51.350} & 5.519 & 65.903 & 0.257 & 78.241 & 0.009 & 146.195 & 0.858 & 118.354 & 0.179 & 55.582 & 0.126 \\
    & & & 30 & - & - & \textbf{31.353} & 7.060 & 35.928 & 0.282 & 62.230 & 0.012 & 77.727 & 0.898 & 68.022 & 0.195 & 40.255 & 0.197 \\
    & & & 40 & - & - & \textbf{16.290} & 6.508 & 18.549 & 0.324 & 59.673 & 0.016 & 42.884 & 0.870 & 39.598 & 0.192 & 32.919 & 0.312 \\
    & & & 50 & - & - & \textbf{9.814} & 7.697 & 11.952 & 0.387 & 46.181 & 0.020 & 26.654 & 0.888 & 21.993 & 0.191 & 24.195 & 0.357 \\
    & & & 60 & - & - & \textbf{7.065} & 5.620 & 8.287 & 0.300 & 38.525 & 0.023 & 16.601 & 0.960 & 12.863 & 0.187 & 20.044 & 0.443 \\
    & & & 100 & 7.094 & 0.429 & \textbf{1.596} & 6.312 & 2.116 & 0.388 & 19.508 & 0.038 & 4.558 & 1.004 & 3.914 & 0.186 & 8.955 & 0.728 \\
    & & & 200 & 0.979 & 0.524 & \textbf{0.245} & 3.116 & 0.316 & 0.392 & 5.198 & 0.073 & 0.839 & 0.999 & 0.536 & 0.191 & 2.610 & 1.532 \\
    \cmidrule(lr){3-18}

    & & \multirow{7}{*}{$0.5$} & 20 & - & - & \textbf{9.699} & 1.256 & 11.759 & 0.094 & 17.802 & 0.008 & 26.046 & 0.841 & 18.112 & 0.169 & 11.208 & 0.123 \\
    & & & 30 & - & - & \textbf{6.557} & 2.977 & 7.963 & 0.118 & 13.902 & 0.012 & 17.410 & 0.838 & 12.841 & 0.185 & 8.357 & 0.192 \\
    & & & 40 & - & - & \textbf{4.118} & 4.105 & 5.195 & 0.132 & 11.181 & 0.016 & 13.164 & 0.857 & 7.563 & 0.181 & 6.482 & 0.291 \\
    & & & 50 & - & - & \textbf{2.841} & 5.514 & 3.201 & 0.157 & 9.166 & 0.020 & 10.131 & 0.912 & 6.120 & 0.173 & 5.132 & 0.357 \\
    & & & 60 & - & - & \textbf{2.204} & 3.594 & 2.713 & 0.140 & 8.245 & 0.023 & 6.755 & 0.932 & 4.384 & 0.181 & 4.716 & 0.443 \\
    & & & 100 & 3.809 & 0.248 & \textbf{0.976} & 6.021 & 1.103 & 0.177 & 5.035 & 0.038 & 2.720 & 0.955 & 1.904 & 0.184 & 2.402 & 0.724 \\
    & & & 200 & 0.784 & 0.299 & \textbf{0.201} & 2.037 & 0.257 & 0.156 & 1.451 & 0.074 & 0.758 & 0.994 & 0.435 & 0.175 & 0.792 & 1.519 \\
    \cmidrule(lr){3-18}

    & & \multirow{7}{*}{$1.0$} & 20 & - & - & \textbf{4.279} & 0.461 & 4.885 & 0.060 & 6.795 & 0.009 & 11.166 & 0.856 & 8.753 & 0.167 & 5.020 & 0.129 \\
    & & & 30 & - & - & \textbf{3.133} & 1.330 & 3.623 & 0.076 & 6.111 & 0.012 & 8.509 & 0.869 & 7.117 & 0.184 & 4.031 & 0.190 \\
    & & & 40 & - & - & \textbf{2.271} & 1.591 & 2.757 & 0.082 & 5.246 & 0.016 & 6.935 & 0.921 & 5.725 & 0.176 & 3.018 & 0.276 \\
    & & & 50 & - & - & \textbf{1.629} & 2.485 & 2.011 & 0.098 & 4.242 & 0.020 & 5.176 & 0.911 & 4.508 & 0.172 & 2.614 & 0.356 \\
    & & & 60 & - & - & \textbf{1.338} & 1.846 & 1.508 & 0.094 & 3.889 & 0.023 & 4.584 & 0.991 & 3.409 & 0.177 & 2.170 & 0.436 \\
    & & & 100 & 2.053 & 0.214 & \textbf{0.615} & 3.115 & 0.784 & 0.106 & 2.373 & 0.037 & 2.242 & 0.975 & 1.460 & 0.179 & 1.398 & 0.720 \\
    & & & 200 & 0.628 & 0.247 & \textbf{0.166} & 1.368 & 0.209 & 0.110 & 0.830 & 0.082 & 0.601 & 1.004 & 0.375 & 0.176 & 0.480 & 1.545 \\
    \bottomrule
\end{tabular}
}
\endgroup
\end{table}

\begin{table}[!htbp]
\caption{Gaps and times on real data.}
\label{tab:real_gap_time}
\centering
\begingroup
\scriptsize
\renewcommand{\arraystretch}{1.15}
\setlength{\tabcolsep}{2.2pt}

\resizebox{0.8\textwidth}{!}{%
\begin{tabular}{ccccc*{7}{cc}}
    \toprule
    \multirow{2}{*}{Dataset}
        & \multirow{2}{*}{$\lambda$}
        & \multirow{2}{*}{$k$}
        & \multirow{2}{*}{$\upsilon$}
        & \multirow{2}{*}{$|\I|$}
        & \multicolumn{2}{c}{\Cref{alg:reduce}}
        & \multicolumn{2}{c}{\Cref{alg:round}}
        & \multicolumn{2}{c}{\Cref{alg:pipage}}
        & \multicolumn{2}{c}{Greedy}
        & \multicolumn{2}{c}{Randomized}
        & \multicolumn{2}{c}{DIHT}
        & \multicolumn{2}{c}{Beam Search} \\
    \cmidrule(lr){6-7}
    \cmidrule(lr){8-9}
    \cmidrule(lr){10-11}
    \cmidrule(lr){12-13}
    \cmidrule(lr){14-15}
    \cmidrule(lr){16-17}
    \cmidrule(lr){18-19}
    & & & & &
    Gap & Time &
    Gap & Time &
    Gap & Time &
    Gap & Time &
    Gap & Time &
    Gap & Time &
    Gap & Time \\
    \midrule

    \multirow{21}{*}{\makecell{Breast Cancer \\ ($n=78,p=4{,}751$)}} & \multirow{7}{*}{$0.1$} & 20 & 1 & 68 & - & - & \textbf{48.923} & 72.758 & 69.535 & 0.540 & 88.918 & 0.011 & 146.942 & 1.147 & 80.137 & 0.242 & 52.709 & 0.123 \\
    & & 30 & 5 & 72 & - & - & \textbf{28.983} & 417.290 & 36.155 & 0.932 & 72.178 & 0.014 & 82.095 & 0.919 & 53.605 & 0.183 & 38.756 & 0.199 \\
    & & 40 & 16 & 68 & - & - & 22.449 & 70.364 & \textbf{20.117} & 0.941 & 57.329 & 0.018 & 65.396 & 0.922 & 36.776 & 0.190 & 33.170 & 0.305 \\
    & & 50 & 25 & 68 & - & - & 15.083 & 110.645 & \textbf{13.372} & 0.925 & 51.779 & 0.024 & 35.141 & 0.869 & 35.870 & 0.179 & 30.446 & 0.387 \\
    & & 60 & 33 & 71 & - & - & \textbf{10.110} & 114.122 & 12.642 & 0.929 & 44.194 & 0.029 & 19.358 & 1.144 & 14.553 & 0.174 & 22.985 & 0.441 \\
    & & 100 & 71 & 71 & 11.191 & 0.750 & \textbf{2.893} & 600.000 & 3.737 & 0.930 & 21.082 & 0.048 & 11.401 & 1.057 & 48.477 & 0.167 & 12.006 & 0.769 \\
    & & 200 & 168 & 69 & 1.611 & 1.102 & 0.626 & 168.879 & \textbf{0.607} & 0.942 & 6.963 & 0.094 & 1.273 & 0.952 & 10.162 & 0.211 & 3.333 & 1.670 \\
    \cmidrule(lr){2-19}

    & \multirow{7}{*}{$0.5$} & 20 & 3 & 57 & - & - & \textbf{7.798} & 3.451 & 9.981 & 0.267 & 12.762 & 0.010 & 25.429 & 0.796 & 17.607 & 0.187 & 11.604 & 0.142 \\
    & & 30 & 7 & 61 & - & - & \textbf{6.095} & 17.419 & 6.660 & 0.308 & 11.191 & 0.014 & 16.237 & 0.864 & 16.012 & 0.187 & 7.562 & 0.189 \\
    & & 40 & 18 & 57 & - & - & \textbf{3.776} & 8.315 & 4.736 & 0.270 & 11.201 & 0.018 & 9.436 & 1.098 & 11.945 & 0.185 & 5.916 & 0.321 \\
    & & 50 & 28 & 58 & - & - & \textbf{2.896} & 9.165 & 3.570 & 0.341 & 9.743 & 0.024 & 7.762 & 0.959 & 9.592 & 0.185 & 4.744 & 0.390 \\
    & & 60 & 38 & 58 & - & - & \textbf{1.845} & 7.588 & 2.500 & 0.366 & 7.356 & 0.028 & 8.983 & 1.200 & 9.373 & 0.180 & 4.093 & 0.434 \\
    & & 100 & 67 & 66 & 4.888 & 0.444 & \textbf{1.247} & 340.405 & 1.438 & 0.451 & 3.302 & 0.046 & 4.398 & 1.159 & 4.711 & 0.183 & 2.191 & 0.739 \\
    & & 200 & 171 & 61 & 1.019 & 0.663 & \textbf{0.364} & 30.646 & 0.451 & 0.631 & 1.229 & 0.092 & 1.429 & 1.315 & 1.025 & 0.191 & 0.763 & 1.664 \\
    \cmidrule(lr){2-19}

    & \multirow{7}{*}{$1.0$} & 20 & 5 & 52 & - & - & \textbf{3.128} & 1.247 & 3.898 & 0.192 & 6.508 & 0.011 & 9.572 & 0.906 & 12.345 & 0.177 & 3.850 & 0.111 \\
    & & 30 & 9 & 52 & - & - & \textbf{2.523} & 1.834 & 3.448 & 0.203 & 6.630 & 0.014 & 8.269 & 0.866 & 14.866 & 0.179 & 2.886 & 0.177 \\
    & & 40 & 24 & 49 & - & - & \textbf{1.870} & 1.601 & 2.312 & 0.294 & 4.437 & 0.020 & 5.387 & 0.896 & 14.370 & 0.182 & 3.437 & 0.248 \\
    & & 50 & 29 & 52 & - & - & \textbf{1.468} & 3.040 & 1.502 & 0.218 & 3.228 & 0.024 & 5.609 & 0.906 & 10.751 & 0.179 & 2.929 & 0.323 \\
    & & 60 & 37 & 54 & - & - & \textbf{1.248} & 7.165 & 1.508 & 0.245 & 2.646 & 0.029 & 3.940 & 0.907 & 9.623 & 0.177 & 2.236 & 0.375 \\
    & & 100 & 72 & 61 & 2.666 & 0.739 & \textbf{0.733} & 49.536 & 0.941 & 0.383 & 1.775 & 0.045 & 2.304 & 1.153 & 5.093 & 0.202 & 1.500 & 0.722 \\
    & & 200 & 173 & 53 & 0.690 & 0.552 & \textbf{0.198} & 3.934 & 0.282 & 0.297 & 0.959 & 0.088 & 0.434 & 1.063 & 2.780 & 0.180 & 0.374 & 1.625 \\
    \midrule

    \multirow{21}{*}{\makecell{Lung Cancer\\ ($n=86,p=7{,}129$)}} & \multirow{7}{*}{$0.1$} & 20 & 0 & 77 & - & - & 64.342 & 185.525 & \textbf{56.493} & 1.255 & 106.697 & 0.023 & 155.672 & 1.325 & 151.491 & 0.368 & 78.537 & 0.190 \\
    & & 30 & 6 & 76 & - & - & \textbf{32.894} & 308.505 & 41.781 & 1.341 & 62.263 & 0.030 & 85.390 & 1.334 & 52.446 & 0.398 & 60.666 & 0.268 \\
    & & 40 & 19 & 75 & - & - & \textbf{25.020} & 167.726 & 27.304 & 1.469 & 52.290 & 0.040 & 53.014 & 1.441 & 39.373 & 0.346 & 41.305 & 0.434 \\
    & & 50 & 22 & 78 & - & - & 18.765 & 600.000 & \textbf{17.808} & 1.561 & 45.881 & 0.046 & 53.078 & 1.551 & 36.672 & 0.348 & 33.181 & 0.555 \\
    & & 60 & 25 & 83 & - & - & \textbf{11.379} & 600.000 & 15.968 & 1.614 & 42.764 & 0.060 & 41.550 & 1.426 & 20.109 & 0.361 & 28.564 & 0.658 \\
    & & 100 & 60 & 77 & 15.948 & 1.357 & 4.523 & 600.000 & \textbf{3.614} & 1.553 & 22.946 & 0.096 & 10.963 & 1.561 & 27.414 & 0.386 & 16.592 & 1.034 \\
    & & 200 & 164 & 79 & 3.729 & 1.914 & \textbf{0.705} & 600.000 & 0.902 & 1.556 & 8.244 & 0.177 & 2.137 & 1.619 & 9.142 & 0.378 & 5.835 & 2.300 \\
    \cmidrule(lr){2-19}

    & \multirow{7}{*}{$0.5$} & 20 & 3 & 57 & - & - & \textbf{10.916} & 5.631 & 14.037 & 0.552 & 17.584 & 0.023 & 36.555 & 1.360 & 23.299 & 0.304 & 12.094 & 0.173 \\
    & & 30 & 11 & 60 & - & - & \textbf{7.637} & 23.530 & 8.232 & 0.610 & 15.051 & 0.032 & 17.600 & 1.433 & 14.198 & 0.322 & 8.810 & 0.312 \\
    & & 40 & 14 & 64 & - & - & \textbf{4.716} & 32.273 & 6.372 & 0.524 & 10.927 & 0.038 & 12.110 & 1.870 & 12.281 & 0.312 & 6.759 & 0.427 \\
    & & 50 & 26 & 65 & - & - & \textbf{3.288} & 27.899 & 4.339 & 0.654 & 9.191 & 0.049 & 13.245 & 1.771 & 9.532 & 0.318 & 6.435 & 0.537 \\
    & & 60 & 31 & 70 & - & - & \textbf{3.001} & 559.650 & 3.548 & 0.874 & 8.446 & 0.055 & 11.242 & 1.477 & 8.094 & 0.324 & 5.833 & 0.641 \\
    & & 100 & 67 & 69 & 5.187 & 0.854 & \textbf{1.729} & 460.412 & 1.881 & 0.788 & 4.904 & 0.093 & 5.243 & 1.691 & 4.421 & 0.339 & 3.444 & 1.045 \\
    & & 200 & 170 & 70 & 1.411 & 1.252 & \textbf{0.415} & 209.193 & 0.562 & 0.866 & 2.027 & 0.172 & 1.596 & 1.822 & 1.248 & 0.377 & 1.727 & 2.313 \\
    \cmidrule(lr){2-19}

    & \multirow{7}{*}{$1.0$} & 20 & 3 & 49 & - & - & \textbf{4.599} & 2.307 & \textbf{4.599} & 0.345 & 6.901 & 0.025 & 10.785 & 1.389 & 12.703 & 0.318 & 5.785 & 0.185 \\
    & & 30 & 12 & 53 & - & - & \textbf{3.270} & 6.918 & 3.969 & 0.379 & 5.280 & 0.033 & 7.070 & 1.562 & 14.683 & 0.294 & 4.373 & 0.316 \\
    & & 40 & 18 & 53 & - & - & \textbf{2.475} & 4.841 & 2.550 & 0.357 & 4.392 & 0.036 & 4.457 & 1.452 & 11.749 & 0.347 & 3.747 & 0.504 \\
    & & 50 & 28 & 53 & - & - & \textbf{1.470} & 2.786 & 1.796 & 0.441 & 4.694 & 0.047 & 4.991 & 1.517 & 10.037 & 0.310 & 3.117 & 0.547 \\
    & & 60 & 33 & 58 & - & - & \textbf{1.305} & 9.455 & 2.091 & 0.548 & 3.584 & 0.063 & 3.627 & 1.519 & 6.239 & 0.326 & 2.697 & 0.684 \\
    & & 100 & 74 & 66 & 3.165 & 0.755 & \textbf{1.033} & 85.601 & 1.050 & 0.512 & 2.685 & 0.096 & 2.577 & 1.855 & 5.271 & 0.322 & 1.833 & 1.040 \\
    & & 200 & 173 & 63 & 0.937 & 1.127 & 0.303 & 32.169 & \textbf{0.297} & 0.634 & 0.781 & 0.168 & 0.780 & 2.045 & 1.724 & 0.341 & 0.475 & 2.221 \\
    \midrule

    \multirow{21}{*}{\makecell{AML\\ ($n=116,p=6{,}283$)}} & \multirow{7}{*}{$0.1$} & 20 & 1 & 93 & - & - & \textbf{57.862} & 600.000 & 59.455 & 0.776 & 67.240 & 0.029 & 131.572 & 1.505 & 106.693 & 0.382 & 63.028 & 0.163 \\
    & & 30 & 6 & 93 & - & - & \textbf{40.284} & 600.000 & 55.714 & 0.789 & 51.806 & 0.040 & 99.023 & 1.844 & 58.079 & 0.351 & 49.286 & 0.258 \\
    & & 40 & 10 & 98 & - & - & 29.445 & 600.000 & \textbf{28.881} & 1.193 & 49.002 & 0.049 & 78.657 & 1.485 & 87.238 & 0.335 & 38.590 & 0.380 \\
    & & 50 & 17 & 100 & - & - & \textbf{18.768} & 600.000 & 25.240 & 1.342 & 41.950 & 0.061 & 68.308 & 1.848 & 30.665 & 0.342 & 34.511 & 0.484 \\
    & & 60 & 26 & 99 & - & - & \textbf{14.924} & 600.000 & 18.041 & 1.114 & 37.893 & 0.073 & 52.302 & 1.548 & 35.114 & 0.331 & 31.876 & 0.616 \\
    & & 150 & 103 & 103 & 11.033 & 1.652 & 3.281 & 600.000 & \textbf{2.816} & 1.612 & 15.532 & 0.160 & 7.527 & 1.334 & 53.370 & 0.332 & 13.489 & 1.494 \\
    & & 200 & 153 & 99 & 4.958 & 2.118 & 1.338 & 600.000 & \textbf{1.257} & 1.677 & 9.979 & 0.219 & 3.516 & 1.683 & 40.371 & 0.353 & 8.551 & 2.369 \\
    \cmidrule(lr){2-19}

    & \multirow{7}{*}{$0.5$} & 20 & 5 & 57 & - & - & 7.788 & 6.453 & \textbf{7.233} & 0.288 & 7.371 & 0.029 & 18.262 & 1.878 & 13.587 & 0.482 & 7.637 & 0.163 \\
    & & 30 & 6 & 70 & - & - & \textbf{5.772} & 103.548 & 6.503 & 0.339 & 8.476 & 0.036 & 18.036 & 1.425 & 13.876 & 0.372 & 6.808 & 0.293 \\
    & & 40 & 13 & 72 & - & - & \textbf{4.308} & 79.801 & 4.450 & 0.413 & 9.535 & 0.046 & 12.827 & 1.776 & 12.900 & 0.368 & 6.336 & 0.407 \\
    & & 50 & 21 & 82 & - & - & \textbf{4.005} & 600.000 & 4.792 & 0.515 & 9.034 & 0.055 & 9.094 & 1.652 & 12.671 & 0.365 & 5.289 & 0.592 \\
    & & 60 & 29 & 84 & - & - & \textbf{2.899} & 600.000 & 4.075 & 0.535 & 7.573 & 0.078 & 8.331 & 1.553 & 8.813 & 0.369 & 4.241 & 0.547 \\
    & & 150 & 108 & 92 & 2.965 & 1.357 & \textbf{1.075} & 600.000 & 1.236 & 0.697 & 3.557 & 0.169 & 2.781 & 1.376 & 3.816 & 0.347 & 1.505 & 1.539 \\
    & & 200 & 157 & 89 & 1.823 & 1.445 & \textbf{0.552} & 600.000 & 0.636 & 0.729 & 1.979 & 0.229 & 1.844 & 1.650 & 2.575 & 0.344 & 1.132 & 2.437 \\
    \cmidrule(lr){2-19}

    & \multirow{7}{*}{$1.0$} & 20 & 6 & 45 & - & - & 2.548 & 0.586 & 2.782 & 0.230 & 3.184 & 0.028 & 5.516 & 1.417 & 6.724 & 0.382 & \textbf{2.455} & 0.160 \\
    & & 30 & 10 & 55 & - & - & \textbf{2.087} & 4.316 & 2.389 & 0.268 & 2.786 & 0.039 & 4.580 & 1.430 & 8.879 & 0.370 & 2.184 & 0.251 \\
    & & 40 & 19 & 52 & - & - & \textbf{1.676} & 6.815 & 1.897 & 0.292 & 2.736 & 0.102 & 4.248 & 1.836 & 12.282 & 0.406 & 2.085 & 0.402 \\
    & & 50 & 27 & 58 & - & - & \textbf{1.284} & 4.834 & \textbf{1.284} & 0.321 & 2.826 & 0.059 & 3.333 & 1.705 & 9.793 & 0.376 & 2.247 & 0.475 \\
    & & 60 & 32 & 67 & - & - & \textbf{1.296} & 71.050 & 1.301 & 0.322 & 2.748 & 0.071 & 2.331 & 1.549 & 10.413 & 0.386 & 2.118 & 0.645 \\
    & & 150 & 111 & 83 & 1.829 & 1.306 & \textbf{0.496} & 600.000 & 0.553 & 0.508 & 1.243 & 0.162 & 1.750 & 1.735 & 5.902 & 0.370 & 0.909 & 1.612 \\
    & & 200 & 163 & 79 & 1.298 & 1.313 & \textbf{0.344} & 600.000 & 0.365 & 0.429 & 1.018 & 0.225 & 0.889 & 1.640 & 4.861 & 0.378 & 0.577 & 2.386 \\
    \midrule

    \multirow{21}{*}{\makecell{Riboflavin\\ ($n=71,p=4{,}088$)}} & \multirow{7}{*}{$0.1$} & 20 & 6 & 58 & - & - & 28.794 & 7.849 & \textbf{24.293} & 0.764 & 35.071 & 0.009 & 89.643 & 0.812 & 43.198 & 0.145 & 35.144 & 0.116 \\
    & & 30 & 8 & 61 & - & - & \textbf{17.638} & 24.627 & 21.662 & 0.739 & 34.227 & 0.013 & 34.280 & 0.720 & 30.659 & 0.151 & 27.138 & 0.171 \\
    & & 40 & 16 & 58 & - & - & 13.152 & 4.588 & \textbf{10.432} & 0.751 & 28.427 & 0.018 & 31.808 & 0.745 & 19.470 & 0.156 & 22.554 & 0.237 \\
    & & 50 & 29 & 57 & - & - & 7.164 & 6.126 & \textbf{6.762} & 0.755 & 24.397 & 0.020 & 22.945 & 0.766 & 11.861 & 0.149 & 17.597 & 0.309 \\
    & & 60 & 34 & 59 & - & - & \textbf{5.237} & 45.155 & 8.638 & 0.760 & 21.909 & 0.027 & 10.517 & 0.800 & 69.115 & 0.155 & 14.248 & 0.337 \\
    & & 100 & 71 & 61 & 5.657 & 0.940 & \textbf{1.868} & 16.609 & 2.341 & 0.760 & 11.187 & 0.039 & 3.767 & 0.852 & 131.891 & 0.167 & 7.277 & 0.608 \\
    & & 200 & 175 & 54 & 1.269 & 0.925 & 0.545 & 9.076 & \textbf{0.542} & 0.763 & 3.817 & 0.073 & 0.795 & 0.927 & 71.078 & 0.161 & 1.887 & 1.515 \\
    \cmidrule(lr){2-19}

    & \multirow{7}{*}{$0.5$} & 20 & 8 & 36 & - & - & 4.989 & 0.449 & \textbf{4.486} & 0.368 & 8.630 & 0.009 & 9.258 & 0.822 & 11.931 & 0.170 & 5.323 & 0.148 \\
    & & 30 & 12 & 48 & - & - & \textbf{2.518} & 1.179 & 3.075 & 0.393 & 6.661 & 0.012 & 7.860 & 0.780 & 7.751 & 0.151 & 4.664 & 0.179 \\
    & & 40 & 23 & 47 & - & - & \textbf{1.904} & 1.513 & 2.142 & 0.433 & 6.226 & 0.016 & 5.189 & 0.886 & 10.803 & 0.143 & 3.428 & 0.265 \\
    & & 50 & 32 & 46 & - & - & \textbf{1.809} & 1.534 & 2.232 & 0.385 & 4.354 & 0.020 & 4.463 & 0.761 & 10.941 & 0.164 & 3.055 & 0.315 \\
    & & 60 & 36 & 50 & - & - & \textbf{1.277} & 3.425 & 1.470 & 0.419 & 4.411 & 0.024 & 3.445 & 0.827 & 11.009 & 0.170 & 2.644 & 0.397 \\
    & & 100 & 77 & 49 & 4.001 & 0.575 & \textbf{0.570} & 1.313 & 0.752 & 0.388 & 3.196 & 0.037 & 1.800 & 0.960 & 12.104 & 0.151 & 1.964 & 0.637 \\
    & & 200 & 174 & 51 & 0.589 & 0.662 & \textbf{0.172} & 3.766 & 0.292 & 0.498 & 0.638 & 0.070 & 0.581 & 1.036 & 3.832 & 0.165 & 0.543 & 1.406 \\
    \cmidrule(lr){2-19}

    & \multirow{7}{*}{$1.0$} & 20 & 9 & 29 & - & - & \textbf{1.932} & 0.294 & 2.181 & 0.220 & 3.365 & 0.009 & 4.988 & 0.780 & 12.869 & 0.147 & 2.225 & 0.119 \\
    & & 30 & 19 & 32 & - & - & \textbf{1.305} & 0.365 & 1.514 & 0.250 & 3.307 & 0.013 & 4.338 & 0.767 & 6.708 & 0.152 & 1.702 & 0.176 \\
    & & 40 & 25 & 38 & - & - & \textbf{0.886} & 0.500 & 1.053 & 0.370 & 2.616 & 0.017 & 2.267 & 0.901 & 8.764 & 0.143 & 1.081 & 0.253 \\
    & & 50 & 37 & 34 & - & - & \textbf{0.656} & 0.360 & 0.866 & 0.213 & 2.172 & 0.020 & 2.840 & 0.761 & 5.994 & 0.154 & 0.868 & 0.310 \\
    & & 60 & 42 & 43 & - & - & 0.669 & 0.979 & \textbf{0.649} & 0.282 & 2.289 & 0.023 & 3.888 & 0.740 & 7.795 & 0.150 & 1.188 & 0.332 \\
    & & 100 & 77 & 44 & 1.498 & 0.435 & 0.357 & 1.914 & \textbf{0.323} & 0.247 & 1.715 & 0.036 & 0.836 & 1.166 & 10.157 & 0.163 & 0.840 & 0.653 \\
    & & 200 & 181 & 44 & 0.222 & 0.502 & 0.134 & 0.981 & \textbf{0.125} & 0.272 & 0.552 & 0.071 & 0.710 & 1.009 & 8.811 & 0.156 & 0.264 & 1.523 \\
    \midrule

    \multirow{21}{*}{\makecell{Leukemia\\ ($n=72,p=3{,}571$)}} & \multirow{7}{*}{$0.1$} & 20 & 4 & 62 & - & - & 45.297 & 5.856 & \textbf{42.356} & 0.339 & 47.525 & 0.009 & 102.254 & 0.756 & 80.037 & 0.116 & 47.494 & 0.110 \\
    & & 30 & 8 & 70 & - & - & \textbf{27.185} & 122.433 & 33.358 & 0.563 & 40.879 & 0.011 & 82.552 & 0.731 & 67.935 & 0.123 & 38.558 & 0.165 \\
    & & 40 & 14 & 67 & - & - & \textbf{16.243} & 102.761 & 23.816 & 0.473 & 28.071 & 0.014 & 47.940 & 0.713 & 27.188 & 0.135 & 31.437 & 0.293 \\
    & & 50 & 21 & 68 & - & - & \textbf{10.096} & 74.981 & 14.235 & 0.439 & 28.470 & 0.017 & 36.772 & 0.730 & 28.418 & 0.157 & 26.699 & 0.282 \\
    & & 60 & 33 & 66 & - & - & 9.424 & 62.761 & \textbf{9.091} & 0.445 & 26.426 & 0.020 & 22.684 & 0.610 & 13.445 & 0.127 & 22.986 & 0.346 \\
    & & 100 & 70 & 62 & 9.113 & 0.605 & \textbf{2.941} & 26.023 & 4.032 & 0.363 & 12.772 & 0.034 & 6.814 & 0.840 & 22.702 & 0.124 & 11.291 & 0.606 \\
    & & 200 & 170 & 58 & 1.869 & 0.615 & \textbf{0.454} & 15.331 & 0.607 & 0.376 & 4.210 & 0.108 & 1.337 & 0.820 & 8.351 & 0.134 & 2.451 & 1.436 \\
    \cmidrule(lr){2-19}

    & \multirow{7}{*}{$0.5$} & 20 & 6 & 44 & - & - & 8.798 & 0.542 & 10.062 & 0.113 & 10.279 & 0.007 & 13.217 & 0.694 & 15.897 & 0.152 & \textbf{8.674} & 0.108 \\
    & & 30 & 12 & 53 & - & - & \textbf{4.854} & 2.463 & 6.069 & 0.150 & 8.852 & 0.011 & 13.171 & 0.786 & 9.338 & 0.136 & 7.254 & 0.171 \\
    & & 40 & 17 & 58 & - & - & \textbf{3.102} & 4.721 & 3.423 & 0.207 & 7.509 & 0.013 & 11.813 & 0.782 & 10.179 & 0.138 & 5.044 & 0.223 \\
    & & 50 & 25 & 57 & - & - & \textbf{2.714} & 8.069 & 3.929 & 0.179 & 6.394 & 0.018 & 10.844 & 1.235 & 6.631 & 0.145 & 4.209 & 0.312 \\
    & & 60 & 36 & 57 & - & - & \textbf{2.078} & 6.388 & 3.363 & 0.227 & 5.482 & 0.022 & 6.937 & 0.802 & 6.741 & 0.137 & 3.876 & 0.338 \\
    & & 100 & 73 & 58 & 3.311 & 0.332 & \textbf{1.117} & 8.463 & 1.221 & 0.233 & 3.224 & 0.033 & 2.984 & 1.031 & 3.099 & 0.131 & 1.843 & 0.579 \\
    & & 200 & 173 & 57 & 0.688 & 0.380 & 0.397 & 37.053 & \textbf{0.359} & 0.230 & 0.756 & 0.063 & 0.780 & 0.978 & 0.823 & 0.104 & 0.784 & 1.319 \\
    \cmidrule(lr){2-19}

    & \multirow{7}{*}{$1.0$} & 20 & 8 & 32 & - & - & \textbf{3.127} & 0.161 & 4.389 & 0.083 & 3.248 & 0.009 & 5.336 & 0.757 & 13.187 & 0.137 & 3.129 & 0.114 \\
    & & 30 & 16 & 35 & - & - & \textbf{2.053} & 0.192 & 2.525 & 0.121 & 2.958 & 0.011 & 7.453 & 0.684 & 12.255 & 0.134 & 2.958 & 0.152 \\
    & & 40 & 21 & 42 & - & - & \textbf{1.953} & 0.571 & 2.390 & 0.107 & 2.834 & 0.015 & 3.914 & 0.829 & 10.139 & 0.143 & 2.500 & 0.221 \\
    & & 50 & 29 & 45 & - & - & \textbf{1.309} & 0.723 & \textbf{1.309} & 0.122 & 2.679 & 0.018 & 2.618 & 0.969 & 5.717 & 0.138 & 2.114 & 0.284 \\
    & & 60 & 37 & 46 & - & - & \textbf{1.070} & 1.127 & 1.229 & 0.103 & 2.376 & 0.019 & 2.667 & 0.827 & 4.585 & 0.138 & 1.412 & 0.344 \\
    & & 100 & 75 & 53 & 1.064 & 0.277 & \textbf{0.523} & 3.069 & 0.725 & 0.162 & 1.640 & 0.034 & 1.761 & 0.766 & 4.222 & 0.110 & 0.811 & 0.587 \\
    & & 200 & 174 & 54 & 0.404 & 0.343 & \textbf{0.240} & 10.504 & 0.250 & 0.151 & 0.719 & 0.062 & 0.931 & 1.060 & 1.197 & 0.136 & 0.401 & 1.313 \\

    \bottomrule
\end{tabular}
}
\endgroup
\end{table}

\end{document}